\documentclass[ejs]{imsart}

\RequirePackage{amsthm,amsmath}
\RequirePackage[numbers]{natbib}
\RequirePackage[colorlinks,citecolor=blue,urlcolor=blue]{hyperref}
\RequirePackage{graphicx}
\usepackage{bm,amsfonts}
\usepackage{multirow}

\RequirePackage[resetlabels]{multibib}
\newcites{supp}{Supplementary References}
\allowdisplaybreaks

\pubyear{2020}
\volume{0}
\issue{0}
\firstpage{1}
\lastpage{8}

\startlocaldefs
\numberwithin{equation}{section}
\theoremstyle{plain}

\newtheorem{theorem}{Theorem}[section]
\newtheorem{proposition}[theorem]{Proposition}
\newtheorem{lemma}[theorem]{Lemma}
\newtheorem{corollary}[theorem]{Corollary}

\theoremstyle{remark}

\newtheorem{assumption}{Assumption}
\newtheorem{example}{Example}[section]

\newtheorem{remark}{Remark}[section]

\endlocaldefs

\begin{document}

\def\T{{ \mathrm{\scriptscriptstyle T} }}
\def\v{{\varepsilon}}
\def\R{\mathbb R}
\def\E{\mathbb E}
\renewcommand{\P}{\mathbb P}
\def\d{{\rm d}}
\renewcommand{\hat}{\widehat}
\renewcommand{\bar}{\overline}
\newcommand{\vv}[1]{\boldsymbol{#1}}
\newcommand{\scr}[1]{\mathscr #1}
\def\N{\mathbb N}
\def\Z{\mathbb Z}
\def\tld{\widetilde}

\begin{frontmatter}
\title{Many-sample tests for the equality and the proportionality hypotheses between large covariance matrices}
\runtitle{Many-sample tests}

\begin{aug}
    \author[A]{\fnms{Tianxing}~\snm{Mei}\ead[label=e1]{meitx@connect.hku.hk}},
    \author[A]{\fnms{Chen}~\snm{Wang}\thanksref{t1}\ead[label=e2]{stacw@hku.hk}}
    \and
    \author[B]{\fnms{Jianfeng} \snm{Yao~}\ead[label=e3]{jeffyao@cuhk.edu.cn}}
   
    \thankstext{t1}{Corresponding author.}
    \address[A]{
      Department of Statistics and Actuarial Science,
      The University of Hong Kong\printead[presep={,\ }]{e1,e2}}
      
      \address[B]{
      School of Data Science,
      The Chinese University of Hong Kong  (Shenzhen)\printead[presep={,\ }]{e3}}
         \runauthor{T. Mei et al.}
\end{aug}







\begin{abstract}
This paper proposes procedures for testing the equality hypothesis and the proportionality hypothesis involving a large number of $q$ covariance matrices of dimension $p\times p$.   
Under a limiting scheme where $p$,
$q$ and the sample sizes from the $q$ populations grow to infinity in a proper manner, the proposed test statistics are shown to be asymptotically normal. Simulation results show that finite sample properties of the test procedures are satisfactory under both the null and alternatives.
As an application, we derive a test procedure for the Kronecker product covariance specification for transposable data. Empirical analysis of datasets from the Mouse Aging Project and the 1000 Genomes Project (phase 3) is also conducted. 
\end{abstract}



\tableofcontents
\end{frontmatter}

\section{Introduction}

Large sample covariance matrices 
serve as a building block of multivariate analysis 
such as principal component analysis, data classification, factor modeling, or regression analysis. However, in such a high-dimensional setting, many traditional statistical methods no longer work, or show very poor performance, assuming a fixed dimension $p$.  Therefore, inference on large sample covariance matrices is of crucial importance. This paper examines some hypothesis testing problems in this regard. 

In the one-sample situation, consider a $p$-dimensional population $\vv{x}\in\mathbf{R}^p$ with covariance matrix $\vv{\Sigma}$. Given a sample from the population, some typical hypotheses about $\vv{\Sigma}$ include:  
 (a) the identity hypothesis:   ``$\vv{\Sigma} = \vv{I}_p$'';             
(b) the sphericity hypothesis: ``$\vv{\Sigma} = c \vv{I}_p$'',  with $c > 0$ unspecified; 
and (c) the diagonality hypothesis: ``$\vv{\Sigma}$ is diagonal''. 
The two-sample situation involves  two $p$-dimensional populations,
say population $\vv{x}$ as above, and another population  $\vv{y}\in\mathbf{R}^p$
with  covariance matrix $\vv{\Sigma}'$.  Given
samples from both populations, we are interested to test hypotheses on the
relationship between $\vv{\Sigma}$ and $\vv{\Sigma}'$.  Typical hypotheses include:
 (d)  the equality  hypothesis:   ``$\vv{\Sigma} = \vv{\Sigma}'$~'';
 (e)   the proportionality hypothesis:  ``$\vv{\Sigma} = c \vv{\Sigma}'$~'' , with $c>0$ unspecified;
 and (f)  the over-dispersion hypothesis: ``$\vv{\Sigma}-\vv{\Sigma}'$ is non-negative definite''.
 
Over the last two decades, tremendous attention has been paid to testing
large covariance matrices where the dimension $p$ is allowed to be large in comparison with the  
sample sizes.  Among the vast literature, most work constructs the test statistic based on some distance 
function between the sample covariance matrices and their targets under the null hypothesis. 
We give below a
few representative papers and class them into three
categories according to the 
mathematical tools employed in deriving the asymptotic null distribution: (a) work  based on  the theory of $U$-statistics or standard
  central limit theorems:
  \cite{Ledoit20021081},
  \cite{Schott2006827,Schott20076535},
  \cite{Chen2010810},
  \cite{Fisher20102554},
  \cite{Srivastava20111090,Srivastava2014289},
  \cite{Li2012908},
  \cite{Qiu20121285},
  \cite{Jiang20122241},
  \cite{Coelho2012627},
  \cite{Ahmad2017500},   \cite{Ishii201999},
  \cite{Zheng2020},
  \cite{HanWu2021}
  and
  \cite{BaiHuWangZhang2021};
(b) work  based on random matrix theory:
  \cite{BJYZ09},
  \cite{WangYao13},
  \cite{LiYao16},
  \cite{Liu2013293},\cite{Tian2015217}
  and \cite{Zhang20183161};
(c) work based on nonparametric statistics:
  \cite{Butucea2017557},
  \cite{Friendly2020144} and \cite{Liao20212775}.
However, almost all of this literature deals with the one-sample  or the
two-sample situations. Only a few  go beyond two samples and
consider several populations. Among  the above long list, only 
\cite{Coelho2012627},
\cite{Ahmad2017500},
\cite{BaiHuWangZhang2021} and 
\cite{Zhang20183161} 
indeed consider several populations, which show that
research on such {
several-sample testing} for large covariance matrices is rather scarce. 
 
The present paper is devoted to the hypothesis testing problem for a
large number of populations. Consider $q$
$p$-dimensional populations  $\vv{x}_1,  \ldots,\vv{x}_q$, with
respective mean $\vv{m}_i$ and covariance matrix $\vv{\Sigma}_i$, $i=1,\ldots,q$.   Given $q$ samples from the populations, each with sample size $n_i$, we aim to test hypotheses about the
relationship among the $q$ covariance matrices  $\vv{\Sigma}_i$.
The key novelty is that  the number of populations  $q$
is allowed to be large in comparison with the dimension $p$ and the sample sizes $n_i$.

Our study is motivated by a few modern genetics data analysis problems where the number of the populations is indeed large.  Two such datasets are analyzed in the paper, namely the {
Mouse Aging Project} dataset \citep{TA2016,TA2019} and a subset of the well-known {
1000 Genomes Project phase 3} dataset \citep{1000gene}. The first data set has $q=16$ populations,  $p=46$ variables, and the 16 sample sizes are  
$n_1=\cdots=n_{16}=40$; and the second (sub)dataset  has $q=26$ populations,  $p=112515$ genes and 
the 26 sample sizes $n_i$'s vary between 64 and 108 (more details on the datasets are given in Sections~\ref{sec:5.1} and \ref{sec:5.2}, respectively). For the first dataset, it is important to know whether the 16 population covariance matrices of dimension 46$\times$46 are proportional, while for the second dataset, an important hypothesis to check is whether the 26  covariance matrices of dimension 112515$\times$112515 are identical.   In both situations, the number of
populations $q$ is not small 
in comparison with the sample sizes, so the existing 
procedures for several-sample testing {
cannot be applied}. 
Addressing such many-population problems is the motivation of this research.  
For definiteness,   we call this new framework  
{
many-sample} hypothesis testing for large covariance matrices.  

{To the best of our knowledge, we are unaware of any existing method for such a many-sample testing problem on large covariance matrices, allowing the number of populations to grow with the dimension and sample sizes}.  We will consider the following hypotheses: 
(g)  the many-sample equality hypothesis:   $\vv{\Sigma}_1=\cdots=\vv{\Sigma}_q$; and
(h) the many-sample proportionality hypothesis: $\vv{\Sigma}_i=a_{ij} \vv{\Sigma}_j$, for $i\ne j$, with $a_{ij}>0$ unspecified.   

In both scenarios, we construct a generalized $U$-statistics involving the $q$  sample covariance matrices to estimate a distance among the $q$ population covariance matrices. The distance is zero under the null hypothesis and increases along the deviation from the null.   Under a proper limiting scheme,
we derive an asymptotic normal distribution for the test statistic under both the null hypothesis and the alternative. The power analysis of the test procedures is also studied under several representative examples.  Simulation results show satisfactory properties of the test procedures in finite samples under both the null and alternatives.

From a technical point of view, the proposed statistics have a complex structure. First of all, the data vectors are high-dimensional with a growing dimension so we do not have a single kernel function but a sequence of kernel functions depending on the growing dimension (in some literature this is called a high-dimensional  $U$-statistic).  Secondly, the statistics involve $q$ samples with different sample sizes $\{n_i\}$, and both the number of populations $q$ and the sample sizes $\{n_i\}$ grow to infinity. That is, the asymptotic analysis is developed under a simultaneous limiting scheme where all the parameters $p$, $q$, and $\{n_i\}$ tending to infinity, which adds great difficulty.

We also apply our procedures to the two genomics datasets mentioned above. For the {
Mouse Aging Project}, our test strongly rejects the hypothesis of proportionality of the 16 covariance matrices (the normalized test statistic is 13.999 to be compared with the standard normal under the null).  The same dataset has also been analyzed in~\cite{TA2019}, where the proportionality hypothesis is accepted which seemingly contradicts our finding. A closer look reveals that their procedure in fact tests an independence hypothesis among the columns of the data matrix.  This independence hypothesis  is 
equivalent to the proportionality hypothesis  {
if} the dataset has a Kronecker product covariance structure. Combining their result and ours thus reveals that the $q$ columns of the dataset are independent but their covariance matrices are not proportional. So there is no contradiction; rather,  the Kronecker product covariance structure assumed in~\cite{TA2019} is unlikely satisfied by this dataset. 

The other application to the {
1000 Genomes Project phase 3} dataset \citep{1000gene}
checks whether the within-genes covariance matrices are identical across the 26 populations, which is strongly rejected by our procedure.  Note that a widely used model for such gene expression datasets is ANOVA which assumes the equality of these covariance matrices. Our test result implies that conclusions from such ANOVA analysis might be questionable.   
 
The rest of the paper is organized as follows. In Section \ref{sec:2}, we develop our many-sample testing procedure for the proportionality hypothesis. 
In Section \ref{sec:kron}, we apply the procedure to general transposable data possessing a Kronecker product covariance structure and propose a specification test.  
In Section \ref{sec:3}, we establish a many-sample testing procedure for the equality hypothesis 
and apply it to the 1000 Genomes Project phase 3 dataset. Detailed proofs of theorems and lemmas and information about the dataset are provided in the supplementary material.

\section{Testing of proportionality}\label{sec:2}

\subsection{Basic settings and assumptions}
Let $\vv{x}_1,\ldots,\vv{x}_q$ be $q$  $p$-dimensional populations  with covariance matrix $\vv{\Sigma}_i=\vv{\Sigma}_{ip}={\rm Cov}(\vv{x}_i)$, $i=1,\ldots,q$. Consider testing  the following hypothesis:
\begin{equation*}\label{eq:2.1}
\begin{split}
    H_0:~~\vv{\Sigma}_i=w_{i}\vv{\Sigma}~~~i=1,\ldots, q.
\end{split}
\end{equation*}  
for some non-specified matrix $\vv{\Sigma}$ and factors~$w_{i}> 0$.
The following  distance is adopted to characterize the proportionality of two matrices: 
 for two $p\times p$ matrices $\vv{A}$ and $\vv{B}$, 
\begin{equation}\label{eq:d_prop}
d_{\rm prop}(\vv{A},\vv{B}):=\frac1{p^2} {\rm tr}
 \left[{\rm tr} (\vv{B})\vv{A} -  {\rm tr} (\vv{A})\vv{B}\right]^2.
\end{equation}
It is easy to see that  $d_{\rm prop}(\vv{A},\vv{B})=0$ if and only if $\vv{A}=w\vv{B}$ for some constant $w$.
{
Furthermore, we define
\begin{equation}\label{eq:M_p}
    M_p = \binom{q}{2}^{-1}\sum_{i<j} d_{\rm prop}(\vv{\Sigma}_i,\vv{\Sigma}_j).
\end{equation}
Note that $M_p$ is non-negative, and $M_p$ vanishes if and only if all $\vv{\Sigma}_j$'s are mutually proportional. In other words, the testing problem can be equivalently reformulated as  
\begin{equation}\label{eq:2.3}
H_0:~~~M_p = 0~~~~\hbox{versus}~~~H_1:~~~M_p>0.
\end{equation}
Inspired by this observation, we propose a test statistic for $H_0$ based on an unbiased estimator 
$U_p$ in (\ref{eq:U_p})
of $M_p$ by samples, of which the key step is to derive an unbiased estimation for the distance $d_{\rm prop}(\vv{\Sigma}_i,\vv{\Sigma}_j)$.
We will elaborate on this topic in the following Section \ref{subsec:2.1}.
}

For each population $\vv{x}_i$, $i=1,\ldots,q$, suppose that we have a sample of size $n_i$, denoted by $\vv{x}_{i,1},\vv{x}_{i,2},\ldots,\vv{x}_{i,n_i}$. Let  $\vv{X}_{ip}=(\vv{x}_{i,1},\ldots,\vv{x}_{i,n_i})$ be the $p\times n_i$ data matrix from the $i$th population. The corresponding sample covariance matrix is
\begin{equation*}\label{eq:2.4}
\vv{S}_{i,p}=\frac{1}{n_i}\vv{X}_{ip}\vv{X}^\T_{ip}=\frac{1}{n_i}\sum_{k=1}^{n_i}\vv{x}_{i,k}\vv{x}^\T_{i,k}.
\end{equation*}
Moreover, the $q$ samples are assumed to be independent of each other. 

The following assumptions are made.

\begin{assumption}\label{assm:2.1}
 For $i=1,\ldots,q$, $\vv{x}_i=\vv{\Sigma}_i^{1/2} \vv{z}_i$, where $\vv{z}_i$ has independent and identically distributed (i.i.d.) entries with zero mean, unit variance, {fourth moment $\nu_{4,i}$}, and finite eighth moment.   In addition, {the $q$ fourth moments $\{\nu_{4,i}\}$ are bounded uniformly for all $q$, that is, there exists a positive constant $\nu_0$ such that $\sup_p \max_{1\leq i\leq q} \nu_{4,i} \leq \nu_0$}.
 
\end{assumption}

\begin{assumption}\label{assm:2.2}
 There exist two positive constants $C$ and $c$ independent of $p$ and $q$  such that 
 \begin{equation}\label{eq:assm2.2}
 \sup_{p} 
 \max_{1\leq i\leq q}
 \|\vv{\Sigma}_i\|\leq C
 ~~~\hbox{and}~~~ 
 \sup_{p}\frac{1}{pq}\sum_{j=1}^q {\rm tr}(\vv{\Sigma}_j) \geq c, 
 \end{equation}
 where $\|\cdot\|$ denotes the operator norm of matrices.
\end{assumption}

\begin{assumption}\label{assm:2.3} 
For each $i=1,\ldots,q$,  $n_i=n_i(p)\to\infty$ as $p\to\infty$ such that $c_{ip}:=p/n_i\to c_i\in(c_0,C_0)$ for some positive constants $c_0$ and $C_0$. 
\end{assumption}

\begin{assumption}\label{assm:2.4} 
{
The number of populations  $q=q(p)\to\infty$ as $p\to\infty$.
}
\end{assumption}

{
The independent component structure in Assumption \ref{assm:2.1}, as a natural extension of multivariate normal distribution, has been commonly adopted in large random matrix theory and related statistical problems; see, for example, \cite{BS10,Chen2010810,yao15}. 
The uniform boundedness condition in Assumption \ref{assm:2.2} helps control the order of the remainder terms in variance calculations. The other condition in (\ref{eq:assm2.2}) ensures the asymptotic variances of our test statistics do not degenerate. These assumptions are technical and can potentially be relaxed in future work.
The scheme given in 
Assumption \ref{assm:2.3} is always referred to as the {\em large-dimensional asymptotics}.  According to the discussion in \cite{ledoit2017numerical,yao15}, a statistical problem with a dimension-to-size ratio between $0.1$ and $10$ is usually viewed as a large-dimensional statistical problem. Assumption \ref{assm:2.4} specifies our many-sample setting, where the number of populations is growing rather than being fixed.
}

\subsection{Unbiased estimation of the distance} 
\label{subsec:2.1}
{
This section is devoted to developing the unbiased estimation of $d_{\rm prop}(\vv{\Sigma}_i,\vv{\Sigma}_j)$ for any $i\neq j$. We denote
\begin{equation*}\label{eq:mu}
    \mu_{i,1} = \frac{1}{p}{\rm tr}(\vv{\Sigma}_i),~~~~\mu_{i,2} = \frac{1}{p}{\rm tr}(\vv{\Sigma}_i^2).
\end{equation*}
and
\begin{equation*}\label{eq:gamma}
    \gamma_{ij} = \gamma_{ji} = \left(\frac{1}{p}{\rm tr}(\vv{\Sigma}_i)\right)\left(\frac{1}{p}{\rm tr}(\vv{\Sigma}_i\vv{\Sigma}_j)\right)\left(\frac{1}{p}{\rm tr}(\vv{\Sigma}_j)\right).
\end{equation*}
Recall the definition (\ref{eq:d_prop}). It holds that
\begin{equation}\label{eq:d_prop_expansion}
    \begin{split}
    d_{\rm prop}(\vv{\Sigma}_i,\vv{\Sigma}_j)
    &=  
    p\left\{\mu_{i,2}\mu_{j,1}^2 + \mu_{j,2}\mu_{i,1}^2 - 2 \gamma_{ij}\right\}.
    \end{split}
\end{equation}
We construct an unbiased estimator for the distance by estimating terms $\mu_{i,2}$, $\mu_{j,2}$, $\mu_{i,1}^2$, $\mu_{j,1}^2$ and $\gamma_{ij}$, respectively.
The following proposition provides the detailed construction of the unbiased estimators.

\begin{proposition}\label{prop:ub_est}
    For any $i=1,\ldots,q$, let
    \begin{align}
    \vv{R}_{i,12} & = \frac{n_i}{p(n_i-1)}{\rm tr}(\vv{S}_{i,p})\vv{S}_{i,p} - \frac{1}{pn_i(n_i-1)}\vv{X}_{ip}\mathcal{D}(\vv{X}_{ip}^\T\vv{X}_{ip})\vv{X}_{ip}^\T, \label{eq:R12}\\
    \vv{R}_{i,2} & =\frac{n_i}{n_i-1}\vv{S}_{i,p}^2 - \frac{1}{n_i(n_i-1)}\vv{X}_{ip}\mathcal{D}(\vv{X}_{ip}^\T\vv{X}_{ip})\vv{X}_{ip}^\T, \label{eq:R2}
\end{align}
where $\mathcal{D}(\vv{B})$ denotes the diagonal matrix made by the diagonal entries of a matrix $\vv{B}$.
We define, for any $i=1,\ldots,q$,
\begin{equation}\label{eq:ub_est_mu}
    \hat{\mu}_{i,12} = \frac{1}{p}{\rm tr}(\vv{R}_{i,12}),~~~~\hat{\mu}_{i,2} = \frac{1}{p}{\rm tr}(\vv{R}_{i,2}).
\end{equation}
In addition, for $1\leq i\neq j\leq q$, we define
\begin{equation}\label{eq:ub_est_gm}
   \hat{\gamma}_{ij}=\hat{\gamma}_{ji} = \frac{1}{p}{\rm tr}(\vv{R}_{i,12}\vv{R}_{j,12}).
\end{equation}
Then, under Assumption \ref{assm:2.1}, $\hat{\mu}_{i,12}$, $\hat{\mu}_{i,2}$, $\hat{\gamma}_{ij}$
are the respective unbiased estimators for $\mu_{i,1}^2$ $\mu_{i,2}$ and $\gamma_{ij}$.
Consequently, let
\begin{equation}\label{eq:h_ij}
\begin{split}
    h(\vv{X}_{ip},\vv{X}_{jp}) &= p\left\{\hat{\mu}_{i,2}\hat{\mu}_{j,12} + \hat{\mu}_{j,2}\hat{\mu}_{i,12} - 2\hat{\gamma}_{ij}\right\}.
    \end{split}
\end{equation}
Then, $h(\vv{X}_{ip},\vv{X}_{jp})$ is an unbiased estimator for the distance $d_{\rm prop}(\vv{\Sigma}_i,\vv{\Sigma}_j)$. 
\end{proposition}
The rest of this section specifies the idea of constructing the unbiased estimators given in the above proposition. The proof of Proposition \ref{prop:ub_est}
can be directly established based on Lemmas \ref{lem:mom} and \ref{lem:ub_est} below. 

The main difficulty in constructing an unbiased estimation of the distance is to estimate the cross-term $\gamma_{ij}$ in (\ref{eq:d_prop_expansion}) unbiasedly since it mixes the inner product and the normalized traces of $\vv{\Sigma}_i$ and $\vv{\Sigma}_j$. To resolve the issue, we introduce the following notation
\begin{equation*}\label{eq:mu_12}
    \mu_{i,12}(\vv{A}) = \left(\frac{1}{p}{\rm tr}(\vv{\Sigma}_i)\right)\left(\frac{1}{p}{\rm tr}(\vv{\Sigma}_i\vv{A})\right),
\end{equation*}
where $\vv{A}$ stands for any $p\times p$ symmetric matrix. 
Note that when $\vv{A}= \vv{I}_p$, $\mu_{i,12} (\vv{I}_p)= \mu_{i,1}^2$. Furthermore, since samples from different populations are mutually independent, 
if we fix $\vv{\Sigma}_j$ and denote
\begin{equation}\label{eq:tildeSgm}
    \widetilde{\vv{\Sigma}_j} = \left(\frac{1}{p}{\rm tr}(\vv{\Sigma}_j)\right) \vv{\Sigma}_j,
\end{equation}
then the cross-term can be further simplified as 
$\gamma_{ij} =\mu_{i,12}(\widetilde{\vv{\Sigma}_j})$. 
Therefore, we first establish an unbiased estimate for $\mu_{i,12}(\vv{A})$ for any deterministic matrix $\vv{A}$ and then utilize the independence of different populations to build an unbiased estimator for $\gamma_{ij}$. 

Let $\vv{A}$ be a deterministic symmetric matrix of size $p$. The quantity $\mu_{i,12}(\vv{A})$ is closely related to the followings: 
\begin{equation*}\label{eq:mu_24}
     \mu_{i,2}(\vv{A}) = \frac{1}{p} {\rm tr}(\vv{\Sigma}_i^2\vv{A}),~~~~
     \mu_{i,4}(\vv{A}) = \frac{(\nu_{4,i}-3)}{p}{\rm tr}(\mathcal{D}(\boldsymbol{\Sigma}_i)\mathcal{D}(\boldsymbol{\Sigma}_i^{1/2}\boldsymbol{A}\boldsymbol{\Sigma}_i^{1/2})).
\end{equation*}
In particular, when $\vv{A} = \vv{I}_p$, $\mu_{i,2}(\vv{I}_p) = \mu_{i,2}$.
These three quantities are connected through the expectations of the following statistics:
\begin{equation*}\label{eq:sample_nu}
\begin{split}
        \hat{\nu}_{i,12}(\vv{A}) & =  \left(\frac{1}{p} {\rm tr}(\boldsymbol{S}_{i,p})\right) \left(\frac{1}{p} {\rm tr}(\boldsymbol{S}_{i,p}\boldsymbol{A})\right),~~~~
        \hat{\nu}_{i,2}(\vv{A}) = \frac{1}{p}{\rm tr}(\boldsymbol{S}_{i,p}^2\boldsymbol{A}), \\
\hat{\nu}_{i,4}(\vv{A}) & = \frac{1}{pn_i(n_i-1)}\sum_{r<s} (\boldsymbol{x}^\T_{i,r}\boldsymbol{x}_{i,r} - \boldsymbol{x}^\T_{i,s}\boldsymbol{x}_{i,s})(\boldsymbol{x}^\T_{i,r}\vv{A}\boldsymbol{x}_{i,r} - \boldsymbol{x}^\T_{i,s}\vv{A}\boldsymbol{x}_{i,s}). 
    \end{split}
    \end{equation*}
The following lemma provides the details of the relationship, and its proof is provided in Appendix \ref{sec:lem_mom}.
\begin{lemma}\label{lem:mom}
For any $p\times p$ deterministic symmetric matrix $\vv{A}$, under Assumption \ref{assm:2.1}, it holds that
    \begin{equation}\label{eq:lsys}
    \mathbb{E}\begin{pmatrix}
        \hat{\nu}_{i,2}(\vv{A})\\
        \hat{\nu}_{i,12}(\vv{A})\\
        \hat{\nu}_{i,4}(\vv{A})
    \end{pmatrix} = 
    \begin{pmatrix}
        1+\frac{1}{n_i} & c_{ip} & \frac{1}{n_i} \\
        \frac{2}{pn_i} & 1 & \frac{1}{pn_i} \\
        2 & 0 & 1
    \end{pmatrix} 
    \begin{pmatrix}
        \mu_{i,2}(\vv{A})\\
        \mu_{i,12}(\vv{A})\\
        \mu_{i,4}(\vv{A})
    \end{pmatrix} 
\end{equation}
\end{lemma}
From the above lemma, we see that though the statistics $\hat{\nu}_{i,12}(\vv{A})$ and $\hat{\nu}_{i,2}(\vv{A})$ are structurally similar to $\mu_{i,12}(\vv{A})$ and $\mu_{i,2}(\vv{A})$, 
they are not unbiased estimators of $\mu_{i,12}(\vv{A})$ and $\mu_{i,2}(\vv{A})$ due to the high-dimensional effects $p/n_i \to c_i$. 
Fortunately, expectations of these three statistics form a closed liner system with respect to $\mu_{i,12}(\vv{A})$, $\mu_{i,2}(\vv{A})$ and $\mu_{i,4}(\vv{A})$.
In addition, the coefficient matrix in (\ref{eq:lsys}) is non-random, invertible, and depends only on the dimension $p$ and the sample size $n_i$. 
Hence, by solving the system of equations, we obtain unbiased estimators for these three quantities. More importantly, the corresponding unbiased estimators $\hat{\mu}_{i,12}(\vv{A})$, $\hat{\mu}_{i,2}(\vv{A})$ and $\hat{\mu}_{i,4}(\vv{A})$ are matrix inner products of matrices $\vv{A}$ and the corresponding $\vv{R}_{i,12}$, $\vv{R}_{i,2}$ and $\vv{R}_{i,4}$ defined, respectively, in (\ref{eq:R12}), (\ref{eq:R2}) and (\ref{eq:R4}). We summarize the results as the lemma below.
\begin{lemma}\label{lem:ub_est}
    Let $\vv{R}_{i,12}$, $\vv{R}_{i,2}$ be the matrices given by (\ref{eq:R12}) and (\ref{eq:R2}) respectively. In addition, define
    \begin{equation}\label{eq:R4}
    \begin{split}
    \vv{R}_{i,4} & =\frac{n_i+2}{n_i(n_i-1)}\vv{X}_{ip}\mathcal{D}(\vv{X}_{ip}^\T\vv{X}_{ip})\vv{X}_{ip}^\T-\frac{2n_i}{n_i-1}\vv{S}_{i,p}^2 \\
    &- \frac{n_i}{n_i-1}{\rm tr}(\vv{S}_{i,p})\vv{S}_{i,p}.
    \end{split}
\end{equation}
For any $p\times p$ deterministic symmetric matrix $\vv{A}$, we define 
\begin{equation*}\label{eq:ub_est}
    \hat{\mu}_{i,12}(\vv{A}) = \frac{1}{p}{\rm tr}(\vv{A}\vv{R}_{i,12}),~~~\hat{\mu}_{i,2}(\vv{A}) = \frac{1}{p}{\rm tr}(\vv{A}\vv{R}_{i,2}),~~~\hat{\mu}_{i,4}(\vv{A}) = \frac{1}{p}{\rm tr}(\vv{A}\vv{R}_{i,4}).
\end{equation*}
Then, under Assumption \ref{assm:2.1}, $\hat{\mu}_{i,12}(\vv{A})$, $\hat{\mu}_{i,2}(\vv{A})$ and $\hat{\mu}_{i,4}(\vv{A})$ are unbiased estimators for $\mu_{i,12}(\vv{A})$,$\mu_{i,2}(\vv{A})$, and $\mu_{i,4}(\vv{A})$, respectively.
\end{lemma}
The proof of this lemma is provided in Appendix \ref{sec:lem_ub_est}. With the help of the lemma, by letting $\vv{A}=\vv{I}_p$, it is easy for us to verify that $\hat{\mu}_{i,12}$ and $\hat{\mu}_{i,2}$ in (\ref{eq:ub_est_mu}) are indeed unbiased estimators to $\mu_{i,1}^2$ and $\mu_{i,2}$, respectively.

Moreover, with the estimation for $\mu_{i,12}(\vv{A})$ in the above lemma, we next verify that $\hat{\gamma}_{ij}$ is an unbiased estimator for the cross-term $\gamma_{ij}$ in (\ref{eq:d_prop_expansion}). Recall that $\gamma_{ij}=\mu_{i,12}(\widetilde{\vv{\Sigma}_j})$. From the definition of $\widetilde{\vv{\Sigma}_j}$ in (\ref{eq:tildeSgm}), we see that
\begin{align*}
    \hat{\mu}_{i,12}(\widetilde{\vv{\Sigma}_j}) = \left(\frac{1}{p}{\rm tr}(\vv{\Sigma}_j)\right)\left(\frac{1}{p}{\rm tr}(\vv{\Sigma}_j \vv{R}_{i,12})\right) = \mu_{j,12}(\vv{R}_{i,12}).
\end{align*}
Thus, we recall the definition of $\hat{\gamma}_{ij}$ in (\ref{eq:ub_est_gm}) to see that
\begin{equation*}
    \hat{\gamma}_{ij} = \frac{1}{p}{\rm tr}(\vv{R}_{i,12}\vv{R}_{j,12}) = \hat{\mu}_{i,12}(\vv{R}_{j,12}) = \hat{\mu}_{j,12}(\vv{R}_{i,12}).
\end{equation*}
Due to the independence of the $i$th and $j$th population,
it holds that
\begin{align*}
    \E[\hat{\gamma}_{ij}] & = \E[\hat{\mu}_{j,12}(\vv{R}_{i,12})] = \E[\mu_{j,12}(\vv{R}_{i,12})] \\
    & = \E[\hat{\mu}_{i,12}(\widetilde{\vv{\Sigma}_j})] =\mu_{i,12}(\widetilde{\vv{\Sigma}_j}) = \gamma_{ij}.
\end{align*}
In other words, $\hat{\gamma}_{ij}$ is indeed unbiased to $\gamma_{ij}$.

Combining the discussion above, we finally see that $h(\vv{X}_{ip},\vv{X}_{jp})$ in (\ref{eq:h_ij}) is indeed an unbiased estimator for the distance $d_{\rm prop}(\vv{\Sigma}_i,\vv{\Sigma}_j)$.

}

\subsection{Test statistic and its asymptotic normality}\label{subsec:2.2}
The proposed test statistic is defined as follows:
\begin{equation*}\label{eq:U_p}
U_p=\binom{q}{2}^{-1}\sum_{i<j} h(\vv{X}_{ip},\vv{X}_{jp}).
\end{equation*}
As mentioned in the Introduction, the statistic $U_p$ is a high-dimensional $U$-statistic with the data vectors (and thus the kernel function) depending on the growing dimension $p$, which complicates our analysis. Another challenge is that we 
have many samples with different sample sizes,  and both the number of populations $q$ and the sample sizes $\{n_i\}$ also grow to infinity.

{
An immediate consequence of the discussion in the last section is that $U_p$ is unbiased to the quantity $M_p$ defined in (\ref{eq:M_p}), i.e., $\E[U_p] = M_p$. 
In particular,
when $H_0$ holds,  we have $d_{\rm prop}(\vv{\Sigma}_i,\vv{\Sigma}_j)=0$ for all possible pairs of $i$ and $j$ and thus $\E[U_p] = M_p = 0$ in the null case.
}

{
To reduce the computational times and complexities in practice, we utilize Equations (\ref{eq:ub_est_mu}), (\ref{eq:ub_est_gm}) and (\ref{eq:h_ij}) in the last section to provide an explicit expression of our $U_p$:
\begin{equation}
    \begin{split}
        U_p & = \frac{2pq}{q-1}\left\{\left(\frac{1}{q}\sum_{i=1}^q \hat{\mu}_{i,2}\right)\left(\frac{1}{q}\sum_{i=1}^q\hat{\mu}_{i,12}\right) - \frac{1}{p}{\rm tr}\left(\frac{1}{q}\sum_{i=1}^q \vv{R}_{i,12}\right)^2\right\}\\
        & - \frac{2p}{q-1}\left\{\left(\frac{1}{q}\sum_{i=1}^n \hat{\mu}_{i,2}\hat{\mu}_{i,12}\right) - \frac{1}{p}{\rm tr}\left(\frac{1}{q}\sum_{i=1}^q \vv{R}_{i,12}^2\right)\right\}.
    \end{split}
\end{equation}
in which the matrix $\vv{R}_{i,12}$ is defined as in (\ref{eq:R12}). It is worth noting that the above expression only involves the calculation of the mean of $q$ numbers or matrices rather than the summation over two distinct indices. Thus, the computational efficiency is significantly improved.
}

{
Next, we focus on the asymptotic distribution of $U_p$. To elaborate on the conclusion, we denote $\langle \vv{A},\vv{B}\rangle = p^{-1}{\rm tr}(\vv{A}\vv{B}^\T)$, the inner product of $p\times p$ matrices and introduce the notation for a fixed $\vv{\Sigma}_i$:
\begin{equation}\label{eq:inner_prod}
\begin{split}
    \langle \vv{A},\vv{B}\rangle_{\vv{\Sigma}_i} &= 2\langle \vv{\Sigma}_i^{1/2}\vv{A}\vv{\Sigma}_i^{1/2}, \vv{\Sigma}_i^{1/2}\vv{B}\vv{\Sigma}_i^{1/2}\rangle \\
    &+ (\nu_{4,i}-3) \langle \mathcal{D}(\vv{\Sigma}_i^{1/2}\vv{A}\vv{\Sigma}_i^{1/2}), \mathcal{D}(\vv{\Sigma}_i^{1/2}\vv{B}\vv{\Sigma}_i^{1/2})\rangle.
    \end{split}
\end{equation}
It is easy to see that $\langle\cdot,\cdot\rangle_{\vv{\Sigma}_i}$ is also an inner product of matrices and thus is non-negative definite, i.e., $\langle\vv{A},\vv{A}\rangle_{\vv{\Sigma}_i} \geq 0$. The structure represented by the new inner product $\langle\cdot,\cdot\rangle_{\vv{\Sigma}_i}$ frequently appears in the calculation of the expectation and variance of products of the quadratic forms related to $\vv{x}_i$; see Appendix \ref{supp:aux} and \ref{app:a1} for more details. This notation will help compute the asymptotic variance of $U_p$ and simplify its expression.

The following theorem shows the asymptotic normality of our $U$-statistic under both null and alternative cases.
\begin{theorem}\label{thm:2.4}
For any $i = 1,\ldots,q$, denote
\begin{equation}\label{eq:alpha_beta_Lambda}
\begin{split}
    \alpha_{i,p} &= \frac{1}{q-1}\sum_{j:j\neq i} \mu_{j,1}^2,~~~
    \beta_{i,p} = \frac{1}{q-1}\sum_{j:j\neq i}\mu_{j,2},\\
\vv{\Lambda}_{i,p} &= \frac{1}{q-1}\sum_{j:j\neq i} \mu_{j,1} \vv{\Sigma}_j,~~~\kappa_{1,i} = \frac{1}{p}{\rm tr}(\vv{\Lambda}_{i,p}\vv{\Sigma}_i).
\end{split}
\end{equation}
In addition, let
\begin{equation}\label{eq:Gamma_ip}
\vv{\Gamma}_{i,p} = \alpha_{i,p}\vv{\Sigma}_i - \mu_{i,1}\vv{\Lambda}_{i,p} + (\mu_{i,1}\beta_{1,i} - \kappa_{1,i})\vv{I}_p. 
\end{equation}
Then, under Assumptions \ref{assm:2.1}---\ref{assm:2.4}, 
it holds that 
\[
\frac{\sqrt{q}(U_p - M_p)}{\sigma_p} \overset{d.}{\to} \mathcal{N}(0,1),
\]
where the asymptotic variance satisfies $\sigma_p^2 = \frac{4}{q}\sum_{i=1}^q \sigma_{i,p}^2$ with
\begin{equation*}\label{eq:thm2.4-2}
\begin{split}
\sigma_{i,p}^2&= 4c_{ip}^2\alpha_{i,p}^2 \mu_{i,2}^2+ 4 c_{ip}\langle \vv{\Gamma}_{i,p},\vv{\Gamma}_{i,p}\rangle_{\vv{\Sigma}_i}.
\end{split}
\end{equation*}
\end{theorem}
}

Theorem \ref{thm:2.4} is proved by adopting the classical Haj\'{e}k projection to the current specific context, and its proof is provided in Appendix \ref{app:thm2.4}. 

{
According to Theorem \ref{thm:2.4}, the asymptotic mean of our statistic $U_p$ is  $M_p$, and the rate of convergence of $U_p-M_p$ is of order $O_p(q^{-1/2})$. 
Moreover, denote
\begin{equation}\label{eq:var_decomp}
    \sigma_{0,p}^2 = \frac{16}{q}\sum_{i=1}^q c_{i,p}^2 \alpha_{i,p}^2\mu_{i,2}^2,~~~\sigma_{r,p}^2 = \frac{16}{q}\sum_{i=1}^q c_{i,p}\langle\vv{\Gamma}_{i,p},\vv{\Gamma}_{i,p}\rangle_{\vv{\Sigma}_i}.
\end{equation}
Under Assumptions \ref{assm:2.2} and \ref{assm:2.3}, it is easy to see that the first term $\sigma_{0,p}^2$ has a positive lower bound.
Hence, we can decompose the asymptotic variance $\sigma_p^2$ into a positive term $\sigma_{0,p}^2$ plus a non-negative term $\sigma_{r,p}^2$. 
In particular, under the null $H_0$, it can be proved that both the mean $M_p$ and the term $\sigma_{r,p}^2$ in the variance vanish.
Consequently, the null distribution for $U_p$ is given as follows.
\begin{corollary}\label{cor:2.5}
Suppose that $H_0$ in (\ref{eq:2.1}) holds. 
Then, under Assumptions \ref{assm:2.1}---\ref{assm:2.4}, we have 
\[
\frac{\sqrt{q}U_p}{\sigma_{0,p}}\overset{d.}{\to}\mathcal{N}(0,1).
\]
\end{corollary}

The critical step of the proof is to verify that all matrices $\{\vv{\Gamma}_{i,p}\}$ are equal to zero. We leave the detailed calculations in Appendix \ref{sec:cor_prop}.

To estimate for the asymptotic variance $\sigma_{0,p}^2$, observe that
\begin{equation}\label{eq:sigma_0_app}
\sigma_{0,p}^2=16 \bigg\{\frac{1}{q}\sum_{i=1}^q c_{ip}^2\mu_{i,2}^2\bigg\}\bigg\{\frac{1}{q}\sum_{i=1}^q\mu_{i,1}^2\bigg\}^2 + O(\frac{1}{q}).
\end{equation}
Define
\begin{equation*}\label{eq:2.13}
\hat{\sigma}_{p}^2=16
\bigg\{\frac{1}{q}\sum_{i=1}^q c_{ip}^2\hat{\mu}_{i,2}^2\bigg\}\bigg\{\frac{1}{q}\sum_{i=1}^q \hat{\mu}_{i,12}\bigg\}^2,
\end{equation*}
where $\hat{\mu}_{i,12}$ and $\hat{\mu}_{i,2}$ are given in (\ref{eq:ub_est_mu}).
We have the following result.
\begin{theorem}\label{thm:2.5}
Under Assumptions \ref{assm:2.1}---\ref{assm:2.4}, $\hat{\sigma}_{p}^2$ is consistent to $\sigma_{0,p}^2$, i.e., $\hat{\sigma}_p^2 \overset{p.}{\to}\sigma_{0,p}^2$. Consequently, under $H_0$, it holds that
\[
    \frac{\sqrt{q} U_p}{\hat{\sigma}_{p}}\overset{d.}{\to}\mathcal{N}(0,1).
\]
\end{theorem}
}

The proof of Theorem \ref{thm:2.5} is given in 
Appendix 
\ref{app:thm2.5}.

In summary, to test $H_0$ in (\ref{eq:2.1}) with an significant level $\alpha$,  we use the test statistic $U_p$ defined in (\ref{eq:U_p}) and reject $H_0$ if $q^{1/2} U_p/\hat{\sigma}_p > z_\alpha$, where $z_{\alpha}$ is the $\alpha$th upper-quantile of the standard normal. 

\subsection{Power of the proposed test}\label{subsec:2.3}
{

In this section, we study the performance of our proposed test procedure under the alternative $H_1$ in (\ref{eq:2.3}) where population covariance matrices are not proportional to each other.  Recall that Theorem \ref{thm:2.4} suggests the test statistic $U_p$ is asymptotically normal under both $H_0$ and $H_1$ with a mean drift $\sqrt{q}M_p$. Based on these results, the next theorem provides the asymptotic representation of the power function and a sufficient and necessary condition for a strong rejection.
\begin{theorem}\label{thm:generalpower}
    Suppose that Assumptions \ref{assm:2.1}---\ref{assm:2.4} holds. For any given significant level $\alpha\in(0,1)$, let $\Phi$ and $z_\alpha$ be the distribution function and the $\alpha$th upper-quantile of the standard normal, respectively. Then,
    the power function satisfies
    \begin{equation}\label{eq:pf}
        \P_{H_1}\left(\frac{\sqrt{q}U_p}{\hat{\sigma}_{p}} > z_\alpha\right)
        =\Phi\left(\frac{\sqrt{q}M_p}{\sigma_p} - \frac{\sigma_{0,p}}{\sigma_p}z_\alpha\right) + o(1),
    \end{equation}
    where $\sigma_p^2$ is defined in Theorem \ref{thm:2.4} and $\sigma_{0,p}^2$ is given in (\ref{eq:var_decomp}).
    In addition,
    the power function tends to $1$  if and only if  $\sqrt{q}M_p \to \infty$. 
\end{theorem}
The detailed proof of the theorem is given in Appendix \ref{sec:proof_generalpower}. 
\begin{remark}\label{rem:power_nece}
Let $d_{\max,q} = \max_{i<j} d_{\rm prop}(\vv{\Sigma}_i,\vv{\Sigma}_j)$. Under the alternative, $d_{\max,q}$ is always positive when $q$ is fixed, but it might decrease to zero as $q$ varies. Observe that 
$
\sqrt{q} d_{\max,q} \geq \sqrt{q} M_p.
$
We immediately find that our test will lose the power if $d_{\max,q} = o(q^{-1/2})$. In other words, a necessary condition to ensure the power is that the convergence rate of $d_{\max,q}$ to $0$ cannot be faster than $q^{-1/2}$. 
\end{remark}
}

{
In what follows, we exhibit several examples to investigate which factors will influence the mean drift $\sqrt{q} M_p$ and derive corresponding sufficient conditions for a significant power. The first one is a particular alternative that significantly deviates from the null.
\begin{example}[$q$ matrices without proportional pair]\label{example:1}
     Suppose that for any $1\leq i\neq j \leq q$, $\vv{\Sigma}_i$ and $\vv{\Sigma}_j$ are not proportional to each other, i.e., $d_{\rm prop}(\vv{\Sigma}_i,\vv{\Sigma}_j)>0$. In addition, we denote
    $d_{1,q} := \min_{i<j} d_{\rm prop}(\vv{\Sigma}_i,\vv{\Sigma}_j)$.
    The setting in this example ensures that $d_{0,q}$ is positive for any $q$ fixed. But, it may vary to zero as $q$ increases. Hence, the converge rate of $d_{1,q}$ influences the performance of the power.
    
    Observe that the mean drift $\sqrt{q}M_p \geq \sqrt{q}d_{1,q}$. 
    According to Theorem \ref{thm:generalpower}, if $\sqrt{q}d_{1,q}\to \infty$, then the power function tends to $1$.
    \hfill$\Box$
\end{example}

The above example is rare. In fact, if there are two mutually proportional populations, $d_{1,q}$ can always be zero for all $q$. A more likely scenario is when $q$ populations can be divided into several sub-groups, with populations within the same group being proportional and populations in different groups not being proportional.
We will discuss this situation in the next example.
\begin{example}
[Non-proportional subgroups without dominant group]\label{example:2}
    Let $m = m(q)$ be a positive integer smaller than $q$ and may vary as $q$ increases.
    Suppose that $\{1,\ldots,q\}$ can be split into $m$ non-empty disjoint subsets $I_1,\ldots,I_m$.     
    We denote $q_i=|I_i|$, the number of elements in $I_i$, $i=1,\ldots,m$. Under our settings, we allow $q_i$'s to vary as $q$ grows. The ratio $q_i/q$ represents the proportion of the $i$th group $I_i$ in the total $q$ populations. In this example, we assume that
    \begin{equation}\label{eq:theta_0}
        \theta_0 := \sup_{q} \max_{1\leq i\leq m} \frac{q_i}{q} <1.
    \end{equation}
    It implies that none of the $m$ groups will be dominant.  
    
    We assume that populations within the same group are mutually proportional, that is, for $i=1,\ldots,m$, $\{\vv{\Sigma}_r:r\in I_i\}$ are proportional to certain basis matrix $\vv{\Lambda}_i$. Here, by taking $p^{-1}{\rm tr}(\vv{\Lambda}_i) =1$ for $i=1,\ldots,m$, the basis matrix is uniquely determined and we also have $\vv{\Sigma}_r = \mu_{r,1}\vv{\Lambda}_i$ for all $r\in I_i$. We denote, for $i=1,\ldots,m$, 
    $W_{i} = q_i^{-1}\sum_{r\in I_i} \mu_{r,1}^2$,
    and let $W_{\min,q}$ be the lower bound for $W_i$'s.
    Moreover, the $m$ bases $\vv{\Lambda}_1,\ldots,\vv{\Lambda}_m$ are not mutually proportional, that is, $d_{\rm prop}(\vv{\Lambda}_i,\vv{\Lambda}_j) >0$ for any $i\neq j$. We denote $d_{2,q} := \min_{1\leq i<j\leq m} d_{\rm prop}(\vv{\Lambda}_i,\vv{\Lambda}_j)$.
    Note that for any $q$ fixed, $d_{2,q}$  is positive, but it may decrease to zero as $q$ increases.

    Next, we show that the divergence of the mean drift $M_p$ is determined by $d_{2,q}$, $W_{\min,q}$ and $\theta_0$. In fact, we 
    observe that
    \begin{align*}
        M_p & = \frac{1}{q(q-1)} \sum_{r\neq s} d_{\rm prop}(\vv{\Sigma}_r,\vv{\Sigma}_s) \\
        & = \frac{1}{q(q-1)}\sum_{i\neq j} \sum_{r\in I_i,s\in I_j} d_{\rm prop}(\vv{\Sigma}_r,\vv{\Sigma}_s)\\
        & = \frac{1}{q(q-1)}\sum_{i\neq j} d_{\rm prop}(\vv{\Lambda}_i,\vv{\Lambda}_j)q_iq_jW_iW_j\\
        & \geq  W_{\min,q}^2 d_{2,q} \sum_{i=1}^m \frac{q_i}{q}\left(1-\frac{q_i}{q}\right) \\
        & \geq W_{\min,q}^2 d_{2,q} \left(1- \theta_0\right).
    \end{align*}
    Consequently, if, in addition, $W_{\min,q}$ has a positive lower bound, a sufficient condition to ensure the power is $\sqrt{q}d_{2,q} \to \infty$.
    In other words, when there is no dominant subgroup and the first-order spectral moments of populations are not too small, to ensure a significant power, the convergence rate of distances between basis matrices $\{\vv{\Lambda}_j\}$ should be slower than $q^{-1/2}$. 
    \hfill$\Box$
\end{example}
The condition (\ref{eq:theta_0}) in the last example indicates that there will not be any group among the $m$ groups that occupy a significant majority. When the condition fails, i.e., $\theta_0 =1$, there comes a more complicated case: the majority of the $q$ populations are mutually proportional with a small number of outliers, which will be discussed in the next example.

\begin{example}[One dominant group with a few outliers]\label{example:3}
    Suppose that there exists a non-empty subset $I_0 \subset\{1,\ldots,q\}$ and an un-specified basis matrix $\vv{\Sigma}$ such that populations outside $I_0$ is proportional to $\vv{\Sigma}$, but within the subset, populations are not. The populations in $I_0$ will be hereafter referred to as {\em outliers}. We denote the minimum distance between outliers and the basis $\vv{\Sigma}$ as $d_{3,q}= \min_{j\in I_0} d_{\rm prop}(\vv{\Sigma}_j,\vv{\Sigma})$.
    Here, $d_{3,q}$ are positive for any $q$ fixed but may decrease to zero as $p,q$ and $\{n_j\}$ grow.
    
    Let $K = |I_0|$ be the number of elements in $I_0$. In this example, we focus on the following small-$K$ assumption: $K = K(q) = o(q)$, which, in particular, allows $K$ to be finite. A direct consequence of the assumption is that the term $\sigma_{r,p}^2$ in (\ref{eq:var_decomp}) is of order $o(1)$ and thus $\sigma_p^2 = \sigma_{0,p}^2 +o(1)$. In other words, the asymptotic variance $\sigma_p^2$ under a small-$K$ alternative approximates the null variance $\sigma_{0,p}^2$. Therefore, the power function in (\ref{eq:pf}) can be further simplified as
    \begin{equation}\label{eq:alterpf}
        \P_{H_1}\left(\frac{\sqrt{q}U_p}{\hat{\sigma}_{p}} > z_\alpha\right)
        =\Phi\left(\frac{\sqrt{q}M_p}{\sigma_{0,p}} - z_\alpha\right) + o(1).
    \end{equation}
    
    Without loss of generality, we assume $p^{-1}{\rm tr}(\vv{\Sigma})=1$ so that $\vv{\Sigma}_i = \mu_{i,1}\vv{\Sigma}$ for any $i\in I_0^c$. Let
    \begin{equation*}\label{eq:W12}
    W_{1,q}=\frac{1}{q-K}\sum_{i\in I_0^c} \mu_{i,1}^2,~~~~W_{2,q} = \frac{1}{q-K}\sum_{i\in I_0^c} c_{ip}^2\mu_{i,1}^4.
    \end{equation*}
    Under Assumption \ref{assm:2.3} and the small-$K$ condition, it is easy to show that both $W_{1,q}$ and $W_{2,q}$ have positive lower bounds. Further, denote $\mu_2=p^{-1}{\rm tr}(\vv{\Sigma}^2)$. We see from (\ref{eq:var_decomp}) that the asymptotic variance can be simplified as $\sigma_{0,p}^2 = 16 \mu_2^2 W_{1,q}^2 W_{2,q}^{1/2}$.

Observe that
\begin{equation*}
\begin{split}
    M_p &=\frac{2}{q(q-1)} 
    \left\{
    \sum_{\{i,j\}\subset I_0} d_{\rm prop}(\vv{\Sigma}_i,\vv{\Sigma}_j)
    +
    \sum_{i\in I_0^c}\sum_{j\in I_0} 
    d_{\rm prop}(\vv{\Sigma}_i,\vv{\Sigma}_j)
    \right\} \\
    & \geq \frac{2}{q(q-1)}\sum_{i\in I_0^c}\sum_{j\in I_0} 
    d_{\rm prop}(\vv{\Sigma}_i,\vv{\Sigma}_j)\\
    & \geq 
    \frac{2K(q-K)}{q(q-1)} d_{3,q} W_{1,q}. 
\end{split}
\end{equation*} 
Therefore, in (\ref{eq:alterpf}), we have
\begin{equation*}\label{eq:alt_meandrift}
    \frac{\sqrt{q}M_p}{\sigma_{0,p}} \geq \frac{(q-K)Kd_{3,q}}{2(q-1)\mu_2\sqrt{qW_{2,q}}} 
\end{equation*}
Consequently, if  
$Kd_{3,q}/\sqrt{q} \to \infty$,
then $\sqrt{q}M_p\to \infty$ and thus the power will eventually tend to $1$. 
To summarize, under the conditions
\begin{equation}
    K=o(q) \quad\text{and}\quad Kd_{3,q}/\sqrt{q} \to \infty,
\end{equation}
the power of the test will tends to 1.

In particular, when $K$ is finite, the above condition reduces to $d_{3,q}/\sqrt{q} \to \infty$, that is, the distances between outliers and the major group should diverge at a rate not slower than $q^{1/2}$. This is understandable because when the number of outliers is limited, to distinguish $H_1$ from $H_0$, the difference between these two must be significant. In this example, this difference is reflected in the divergence of $d_{3,q}$ at a rate faster than $q^{1/2}$. 
    \hfill$\Box$
\end{example}
}

{
At the end of this section, 
we consider the problem of ``finding a needle in a haystack'', where there is only a single outlier group $(K=1)$ among a growing number of populations. On the one hand, this example itself is quite complicated since, with only one outlier, the alternative and the null are very close to each other. On the other hand, unlike in Example \ref{example:3} where we could only obtain lower bound estimates for the mean drift and the power function when there is only one outlier, we can derive an accurate expression of the mean drift and the power function. This enables us to validate our theoretical results in the simulation part (Section \ref{ssec:2.5}).
\begin{example}[Finding a needle in a
haystack]\label{example:4}
Consider a special case of $H_1$ as follows:
\begin{equation}\label{eq:simu0}
\begin{split}
    H^{(\beta)}_1:~&\vv{\Sigma}_i =
    \begin{cases} 
        w_i \vv{\Sigma}, & \quad i=1,\ldots,q-1; \\
        \vv{\Sigma}+\sqrt{\beta}\vv{\Lambda}, & \quad i = q, 
    \end{cases}  
\end{split}
\end{equation}
where $w_i>0$ are non-specified factors and the matrix $\vv{\Lambda}$ is not proportional to $\vv{\Sigma}$, i.e., $d_{\rm prop}(\vv{\Sigma},\vv{\Lambda})>0$. Here, we introduce a non-negative variable $\beta$ to control the degree of deviation of $H_1^{(\beta)}$ from $H_0$. In fact, when $\beta =0$, $H_1^{(0)} = H_0$; and as $\beta$ increase, $H_1^{(\beta)}$ gradually deviates from $H_0$. Without loss of generality, we assume $p^{-1}{\rm tr}(\vv{\Sigma}) =1$ so that $w_i=\mu_{i,1}$ for $i=1,\ldots,q-1$.

For $\beta>0$, since there is only a single outlier, we have
\begin{align*}
    d_{3,q}& =d_{\rm prop}(\vv{\Sigma},\vv{\Sigma}+\sqrt{\beta}\vv{\Lambda})\\
    & = 
    {\rm tr}\left[
    \vv{\Sigma}
    \left(\frac{1}{p}{\rm tr}(\vv{\Sigma})+\sqrt{\beta}\frac{1}{p}{\rm tr}(\vv{\Lambda})\right) - (\vv{\Sigma}+\sqrt{\beta}\vv{\Lambda})
    \left(\frac{1}{p}{\rm tr}(\vv{\Sigma})\right)\right]^2 \\
    & = \beta d_{\rm prop}(\vv{\Sigma},\vv{\Lambda}).
\end{align*}
In addition, the mean drift
\begin{equation*}
\sqrt{q}M_p = \frac{2\sqrt{q}}{q(q-1)}\sum_{i=1}^{q-1} d_{\rm prop}(\vv{\Sigma}_q,\vv{\Sigma}_i) = \frac{2}{\sqrt{q}}\beta d_{\rm prop}(\vv{\Sigma},\vv{\Lambda}) W_{1,q}
\end{equation*}
Meanwhile, we also have $\sigma_{0,p}^2 = 16\mu_2^2 W_{1,p}^2 \sqrt{W_{2,p}}$ and 
\begin{align*}
    \frac{\sqrt{q}M_p}{\sigma_{0,p}} &= \beta \frac{ d_{\rm prop}(\vv{\Sigma},\vv{\Lambda})}{2\sqrt{qW_{2,q}}\mu_2}. 
\end{align*}
Combining with (\ref{eq:alterpf}), the power function now becomes
\begin{equation}\label{eq:sim1}
    \P_{H_1}\left(\frac{\sqrt{q}U_p}{\hat{\sigma}_{p}} > z_\alpha\right)
        =\Phi\left(\beta \frac{ d_{\rm prop}(\vv{\Sigma},\vv{\Lambda})}{2\sqrt{q W_{2,q}}\mu_2} - z_\alpha\right) + o(1).
\end{equation}
Therefore, if we assume 
\[
    \liminf_{p,q\to\infty} \frac{d_{\rm prop}(\vv{\Sigma},\vv{\Lambda})}{\sqrt{q}}  >0,
\]
then, as $\beta$ increases, the value of the power function will gradually tend to one.
    \hfill$\Box$
\end{example}  
}

\subsection{Simulation results}\label{ssec:2.5}

In this section, we conduct a numerical study to evaluate the finite-sample performance of our test.  The simulation design follows the setting in Example \ref{example:4} with dimensions $p=100$, {the number of populations $q= 50,100,150$, and the sample sizes $\{n_j\}$ randomly picked from $\{50,\ldots,150\}$}. The populations have the structure $\vv{x}_i=\vv{\Sigma}^{1/2}\vv{z}_i$ as in Assumption \ref{assm:2.1}, where the noise $\vv{z}$ has {\color{blue}i.i.d.} entries with {zero mean and unit variance}. Moreover, we consider two different types of noises, namely the standard normal noise $z_{ij}\sim\mathcal{N}(\vv{0},\vv{I}_p)$ and the Gamma noise $ z_{ij}\sim{\sf Gamma}(4,2)-2$.

 { Recall that in Example \ref{example:4}, the last population is a single outlier $(K=1)$. In this section, we set $\vv{\Sigma}_q = \vv{\Sigma}_0 + \sqrt{\beta}\vv{\Lambda}_0$}, and the other populations are all proportional to $\vv{\Sigma}_0$. 
 Here, {the variable $\beta$ will vary from $0$ to a positive number $\beta_{\max}$, which is selected according to different settings to control a proper increasing rate of the power. Given two non-negative definite matrices $\vv{\Lambda}$ and $\vv{\Sigma}$, we define
 \[
 \vv{\Sigma}_0 = \vv{\Sigma}/(p^{-1}{\rm tr}(\vv{\Sigma}))~~~\hbox{and}~~~ \vv{\Lambda}_0=\vv{\Lambda}/\sqrt{p^{-1}d_{\rm prop}(\vv{\Lambda},\vv{\Sigma})},
 \]
 which are the respective normalizations of $\vv{\Sigma}$ and $\vv{\Lambda}$. The aim of the normalization procedure is to ensure $p^{-1}{\rm tr}(\vv{\Sigma}_0) =1$ and $p^{-1}d_{\rm prop}(\vv{\Sigma}_0,\vv{\Lambda}_0) = 1$, which helps reduce the theoretic curve of power in (\ref{eq:sim1}) to be
 \begin{equation}\label{eq:sim1''}
     \P_{H_1}\left(\frac{\sqrt{q}U_p}{\hat{\sigma}_p} > z_\alpha \right) =
     \Phi\left(\frac{\beta p}{2\mu_2\sqrt{q W_{2,q}}}\right) + o(1),
 \end{equation}
 where $\mu_2 = p^{-1}{\rm tr}(\vv{\Sigma}_0^2)$ and $W_{2,q} =(q-1)^{-1}\sum_{i=1}^{q-1} c_{ip}^2\mu_{i,1}^4$.}

The following combinations of $\vv{\Sigma}$ and $\vv{\Lambda}$ are considered: 

(a) $\vv{\Sigma}=\vv{I}_p$, and $\vv{\Lambda}=\vv{U}{\rm diag}\{\sqrt{\beta},\cdots,\sqrt{\beta},0,\ldots,0\}\vv{U}^\T$, where the first $p/2$ entries equal to $\sqrt{\beta}$ and $\vv{U}$ is a randomly chosen orthogonal  matrix; 

(b) $\vv{\Sigma}=\vv{U}_1 \vv{D}_1 \vv{U}_1^\T$ and $\vv{\Lambda}=\vv{U}_2 \vv{D}_2 \vv{U}_2^\T$, where $\vv{U}_1$ and $\vv{U}_2$ are two independently and randomly chosen orthogonal matrices, $\vv{D}_1$ and $\vv{D}_2$ are diagonal matrices with their entries randomly picked from the interval $(e^{-3},e^3)$, respectively.

\begin{figure}
\centering
\includegraphics[width = 43 mm, height = 43 mm]{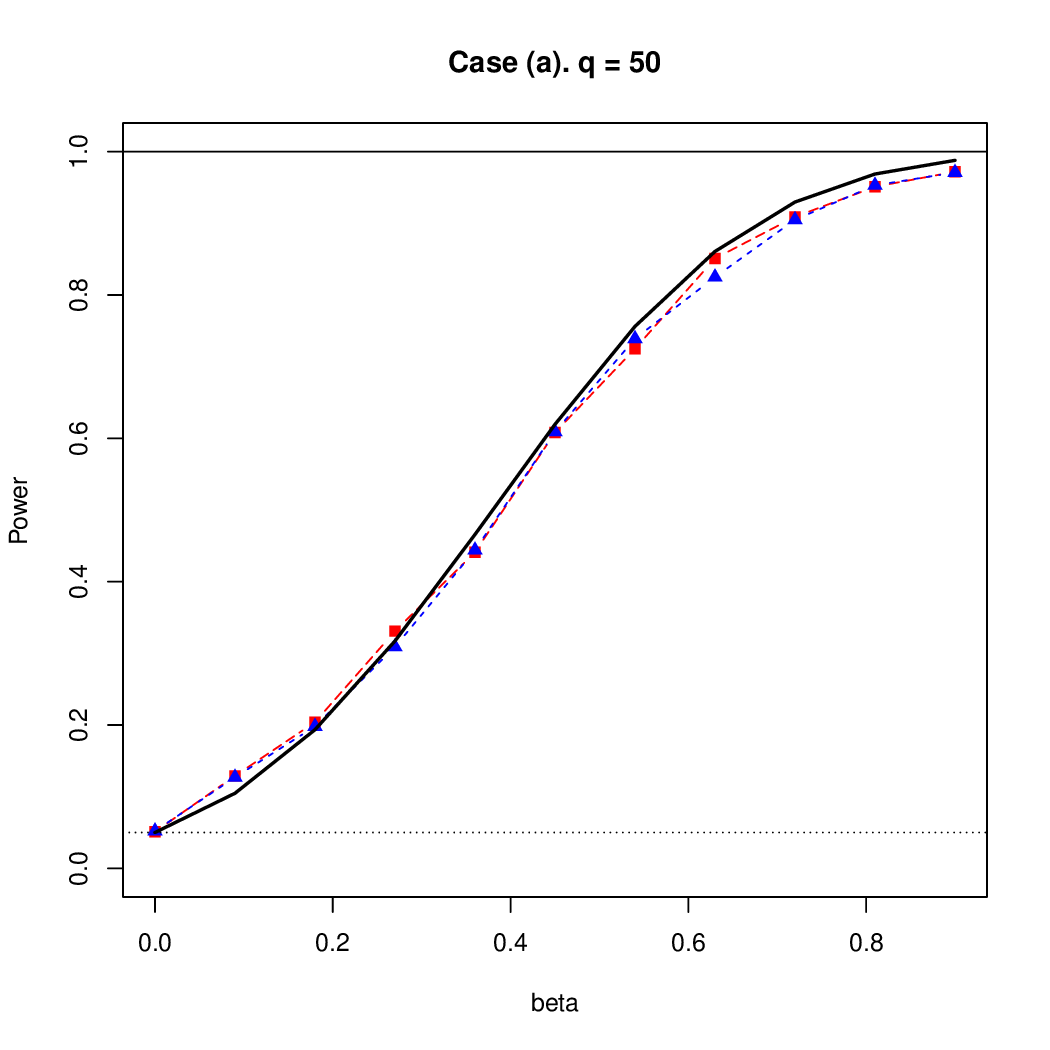}
\includegraphics[width = 43 mm, height = 43 mm]{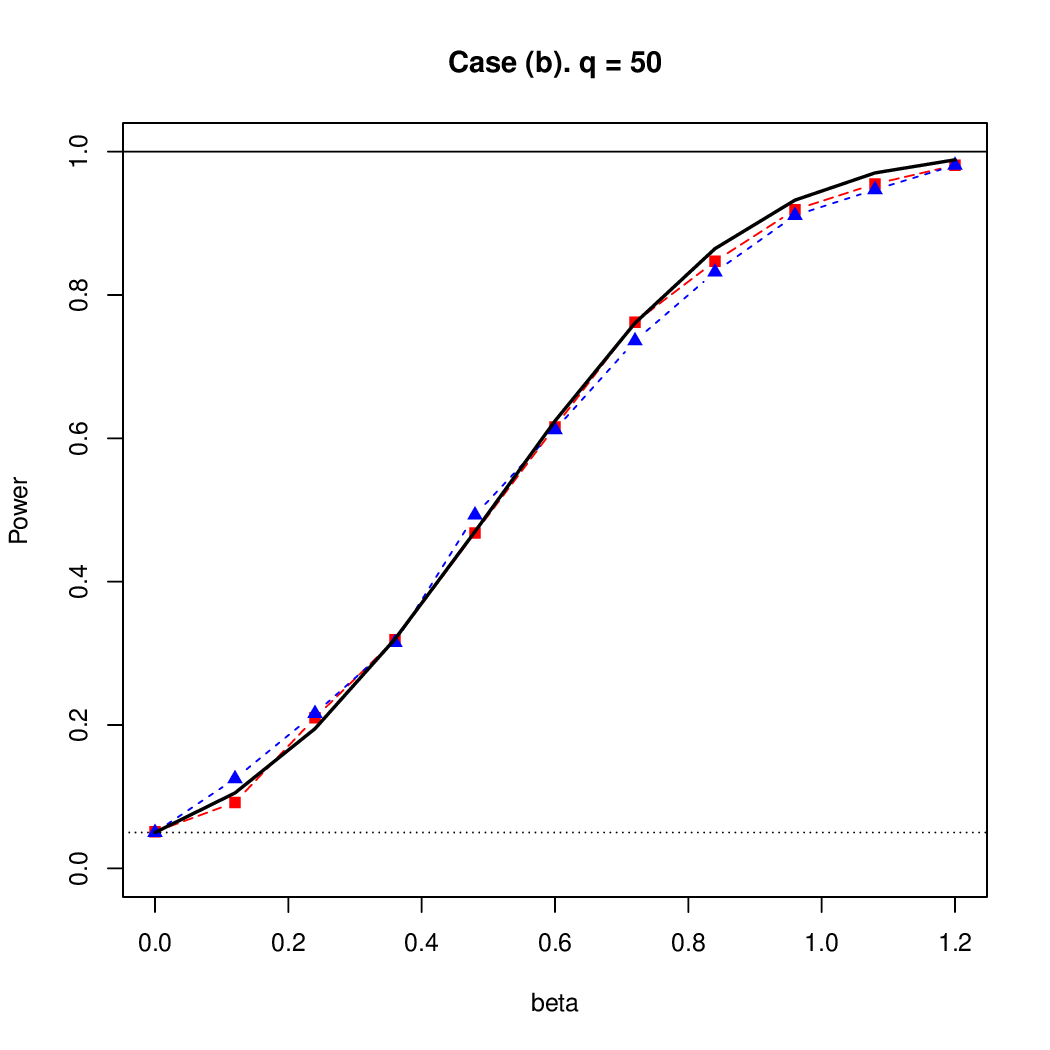}
\includegraphics[width = 43 mm, height = 43 mm]{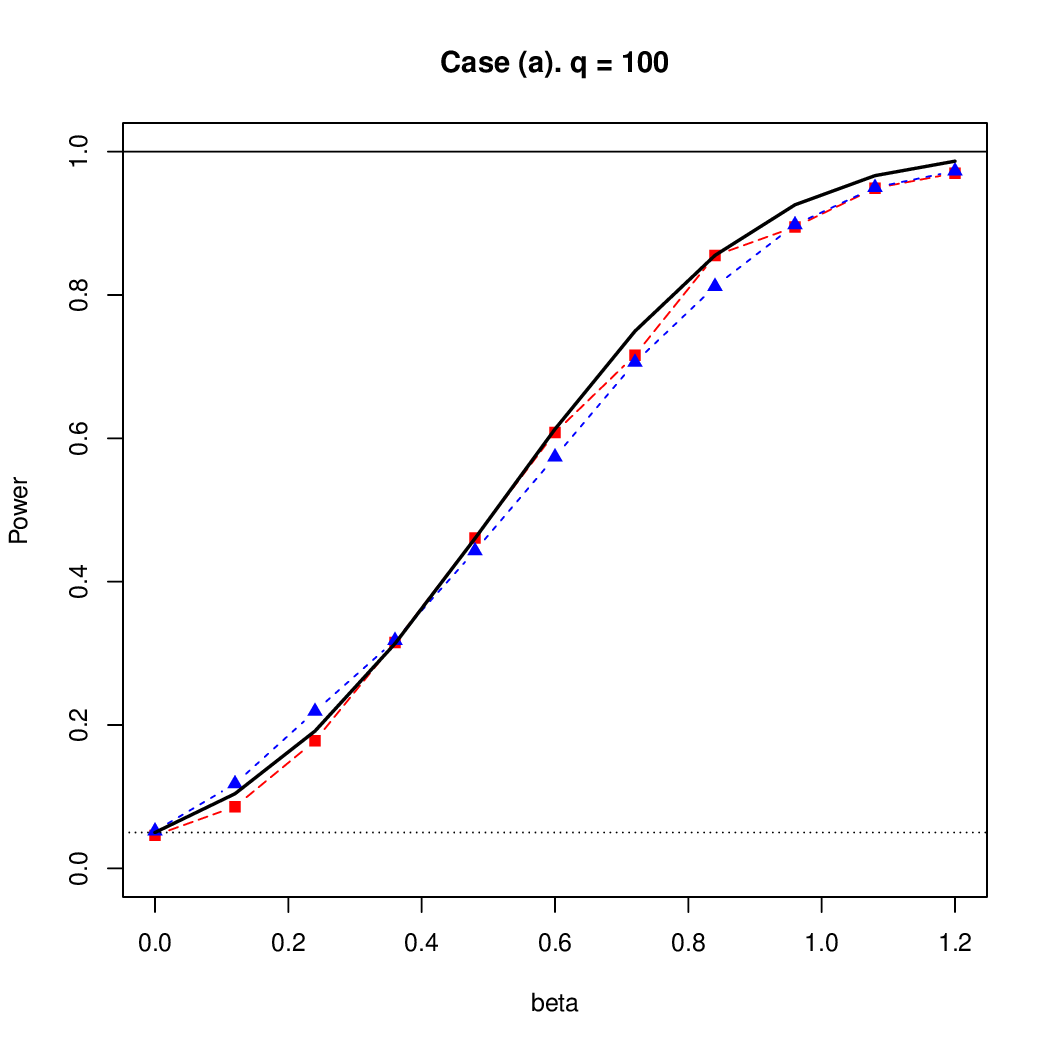}
\includegraphics[width = 43 mm, height = 43 mm]{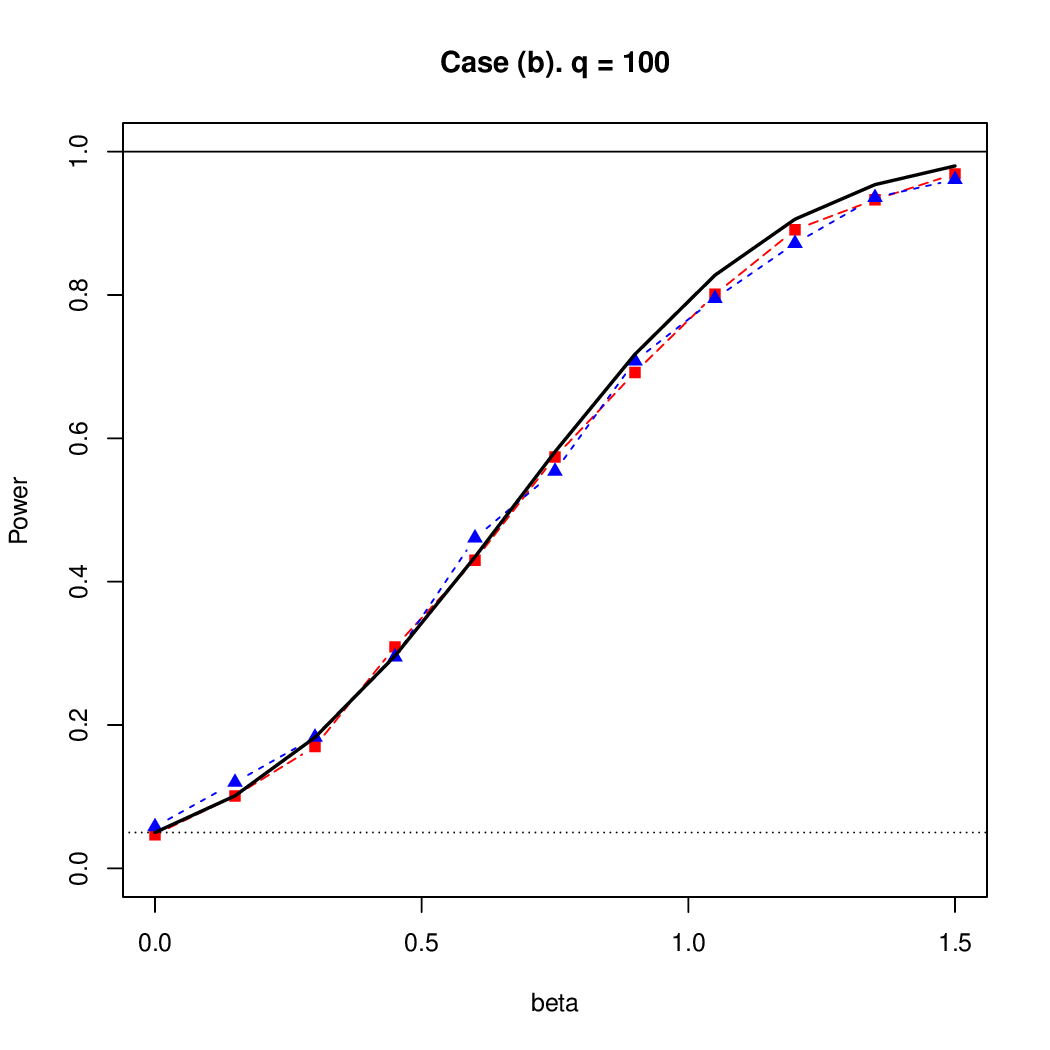}
\includegraphics[width = 43 mm, height = 43 mm]{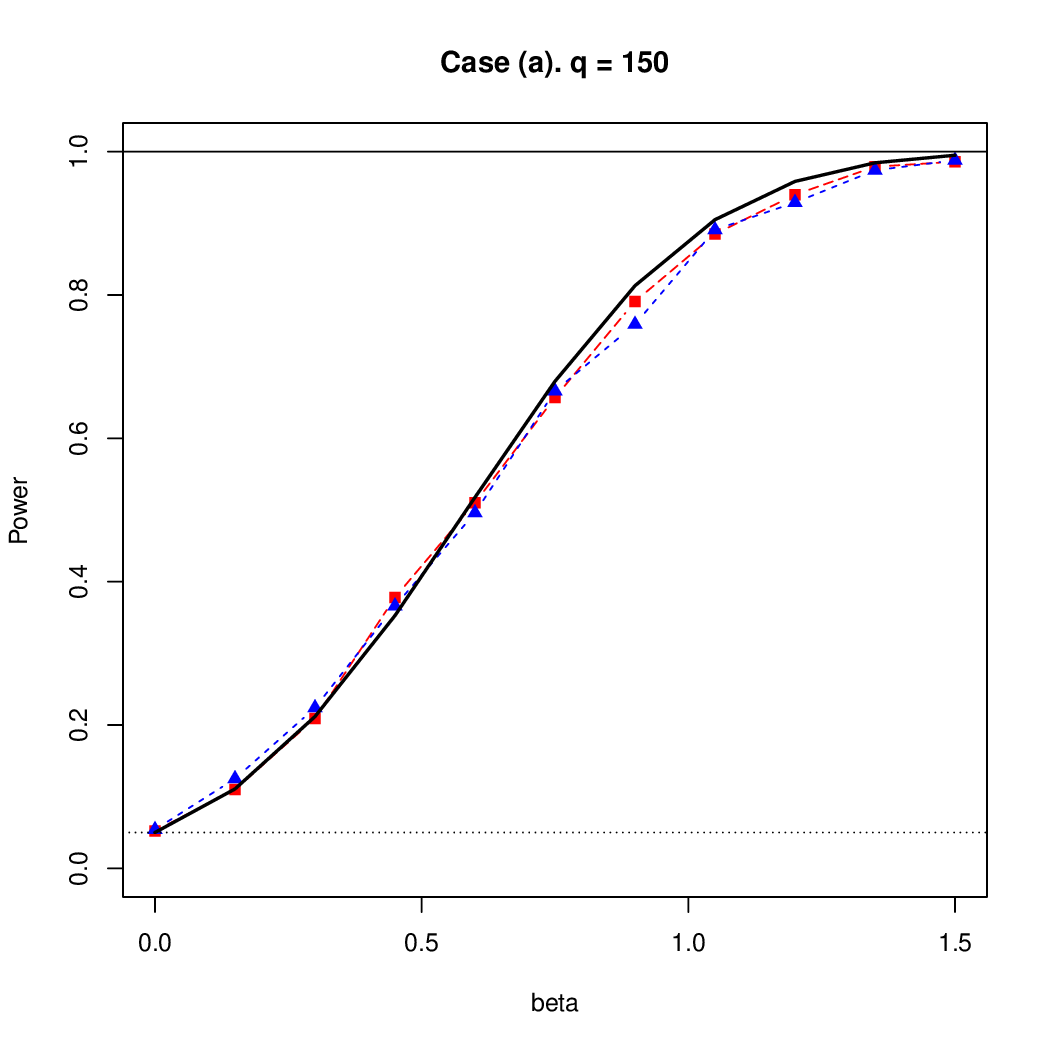}
\includegraphics[width = 43 mm, height = 43 mm]{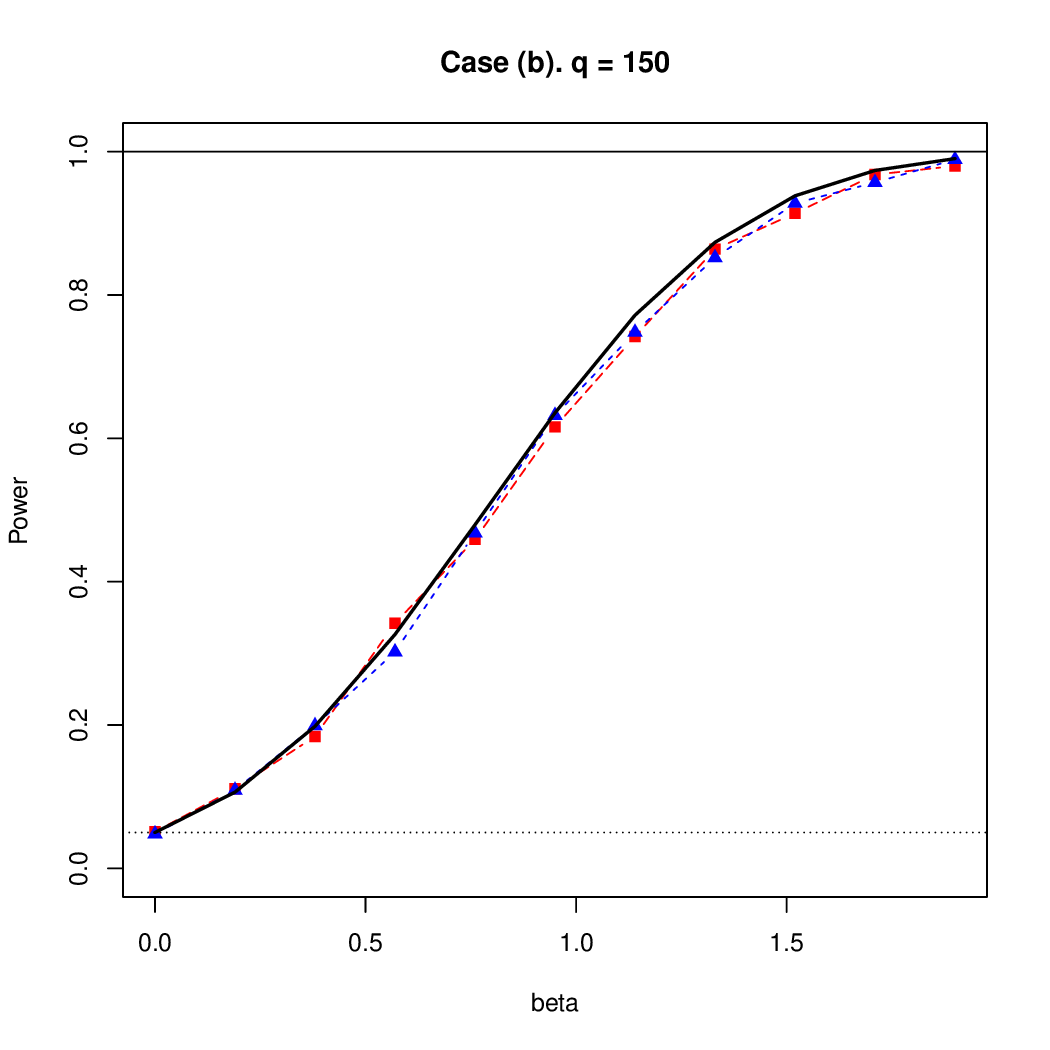}
\caption{Empirical power plot of the many-sample proportionality test. Two cases for the pair $(\vv{\Sigma},\vv{\Lambda})$ with two different noises, normal noise (red square) and Gamma noise (blue triangle), under three different settings: $q=50, 100$, and $150$. The black dash lines stand for the theoretical power curves given in (\ref{eq:sim1''}).}
\label{fig:1}
\end{figure}

{For $i=1,\ldots,q-1$}, we randomly pick a factor $w_i$ uniformly from $(0.5,1.5)$ and set $\vv{\Sigma}_i=w_i\vv{\Sigma}_0$. {For each $\beta$ from $0$ to $\beta_{\max}$ with step size $0.1$,} we repeat simulations $1000$ times for the two cases of $(\vv{\Sigma}, \vv{\Lambda})$,  each with two different noises, and compute the empirical size and power. The empirical sizes under $H_0$ (or equivalently, $\beta=0$) are shown in Table \ref{tab:1}, and the curves of {both theoretical and empirical} powers are displayed in Figure \ref{fig:1}. It can be observed that our proposed testing statistic controls the size well under the null case and effectively rejects $H_0$ when {the outlier covariance matrix $\vv{\Sigma}_0+\sqrt{\beta}\vv{\Lambda}_0$} gradually deviate from $\vv{\Sigma}_0$. {Also, the empirical power curves approximate the theoretical ones closely, which supports our previous theoretical analysis.}

\begin{table}
\caption{Empirical sizes of the many-sample proportionality test.} 
\label{tab:1}
\begin{tabular}{@{}ccccc@{}}
\hline
&\multicolumn{2}{c}{Case (a)}&
\multicolumn{2}{c}{Case (b)}
 \\ [3pt]
\hline
&Normal & Gamma & 
Normal & Gamma
 \\
\hline
$q = 50$& $0.053$ & $0.056$ & 
$0.052$ & $0.054$
 \\
$q=100$ &  $0.046$& $0.048$ & $0.045$ & 
$0.054$
\\
$q=150$ &  $0.055$& $0.053$ & $0.050$ & 
$0.053$
\\
\hline
\end{tabular}
\end{table}


\section{A specification test for transposable data}\label{sec:kron}

\subsection{A specification test for Kronecker product dependence structure }\label{ssec:2.6}
Transposable data is a special class of matrix-valued data 
whose rows and columns represent two different sets of variables \citep{AllenTibsh2010,AllenTibsh2012}. Recent studies of transposable data can be found in various areas such as genomics analysis and longitudinal data in financial markets, for example \cite{Ning2013,TA2016,Leng2012,Lee2013}. 

Modeling dependence between rows and columns has primary importance in transposable data analysis. A widely adopted model for this purpose is the  {\it Kronecker product dependence structure},  where the $p\times q$ transposable data $\tilde{\vv{X}}$ satisfies
\begin{equation}\label{eq:kron}
\tilde{\vv{X}}=\vv{\Sigma}_R^{1/2}\vv{Z}\vv{\Sigma}_C^{1/2},
\end{equation} 
where $\vv{\Sigma}_R$ is a $p\times p$ row covariance matrix, $\vv{\Sigma}_C$ is a $q\times q$ column covariance matrix, and $\vv{Z}$ is a $p\times q$ random matrices having i.i.d. entries with zero mean and unit variance. Note that for $r_1,r_2\in\{1,\ldots,p\}$ and $c_1,c_2\in\{1,\ldots,q\}$, ${\rm cov}(\tilde{\vv{X}}_{r_1,c_1},\tilde{\vv{X}}_{r_2,c_2})=(\vv{\Sigma}_R)_{r_1,r_2}(\vv{\Sigma}_C)_{c_1,c_2}$. Hence,  the covariance matrix of the vectorization $\mathop{\text{vec}}(\tilde{\vv{X}})$ of the transposable data is   $\vv{\Sigma}=\vv{\Sigma}_C\otimes\vv{\Sigma}_R$, namely the Kronecker product of row and column covariance matrices (this covariance structure is also referred as {
separable covariance model} in some other fields). According to \cite{TA2019}, ``... the matrix-variate normal distribution ... is widely used to model high-dimensional transposable data" (p.1310), and this type of distribution does have the Kronecker product covariance structure. Readers are referred to \cite{TA2019} for more justifications of the Kronecker product structure assumption.

There is an unexpectedly interesting relationship between the Kronecker dependence structure and the framework of many-sample tests we consider in the paper. To see this, let us enlarge the model into a generic  {\it
independent component model} (ICM) is defined as follows.   A  $p\times q$ transposable data $\tilde{\vv{X}}$ follows an ICM if $\mathop{\text{vec}}(\tilde{\vv{X}})=\vv{\Gamma}^{1/2}\vv{z}$,
where  $\vv{\Gamma}$ is a $(p\times q)\times (p\times q)$ non-negative definite matrix and $\vv{z}$ is a  $p\times q$-dimensional random vector having i.i.d. entries with zero mean and unit variance.
Clearly, the Kronecker product structure model is a sub-model of ICM with $\vv{z}=\mathop{\text{vec}}(\vv{Z})$ and 
$\vv{\Gamma}=\vv{\Sigma}=\vv{\Sigma}_C\otimes\vv{\Sigma}_R$.

An important problem for transposable data is to test the following  hypothesis
\begin{equation}\label{eq:hykron}
    H^*_0:~\tilde{\vv{X}}~\hbox{satisfies (\ref{eq:kron}) with a diagonal $\vv{\Sigma}_C$}.
\end{equation} 
This is referred to as the diagonal hypothesis for the between-column covariance matrix  $\vv{\Sigma}_C$, see \cite{TA2019}. 
Lemma \ref{lem:kron} below 
establishes the connection between this diagonal hypothesis under Kronecker product structure and the many-sample proportionality test introduced in Section~\ref{sec:2}.   
\begin{lemma}\label{lem:kron}
Suppose that a $p\times q$ transposable data $\tilde{\vv{X}}$ follows ICM. Then, the following three statements are equivalent. 

(i) $\tilde{\vv{X}}$ satisfies (\ref{eq:kron}) with a diagonal $\vv{\Sigma}_C$; 

(ii) $\tilde{\vv{X}}$ satisfies (\ref{eq:kron}) with independent columns;

(iii) The columns of $\tilde{\vv{X}}$ are independent and their covariance matrices are proportional to each other.
\end{lemma}

The proof of the lemma is provided in Appendix 
\ref{app:4}.
By Lemma \ref{lem:kron},  the diagonality hypothesis $H_0^*$ (item (i)) under the Kronecker product structure is equivalent to the proportionality hypothesis $H_0$ we considered in Section \ref{sec:2} for independent columns (item (iii)).
This connection enables the use of our test for the diagonality hypothesis $H_0^*$.
Recall that one advantage of our test is that both dimensions $p$ and $q$ are allowed to grow to infinity with the sample size $n$.

Let $\tilde{\vv{X}}_1,\ldots,\tilde{\vv{X}}_n$ be an i.i.d. sample from a $p\times q$ transposable population $\tilde{\vv{X}}$. By Lemma \ref{lem:kron}, our procedure of testing for the hypothesis $H_0^*$  (\ref{eq:hykron}) is summarized as follows: 

(1) Let $\tilde{\vv{x}}_{i,j}$ be the $i$th column of $\tilde{\vv{X}}_j$, for $i=1,\ldots,q$ and $j=1,\ldots,n$. Define $\vv{X}_{ip}=(\tilde{\vv{x}}_{i,1},\ldots,\tilde{\vv{x}}_{i,n})$, $i=1,\ldots,q$.

(2) Compute the value of the test statistic $U_p$ and 
the asymptotic variance estimate $\hat{{\sigma}}_p^2$ defined in Section \ref{subsec:2.2}.

(3) With significance level $\alpha\in(0,1)$, reject $H_0$ if  $q^{1/2}U_p/\hat{{\sigma}}_p>z_{\alpha}$, where $z_{\alpha}$ is the $\alpha$th upper-quantile of the standard normal.

\subsection{Comparison with test in \cite{TA2019}}\label{subsec:2.4}

The procedure proposed in \cite{TA2019} for the diagonality 
hypothesis  $H_0^*$ (\ref{eq:hykron}) is based on the key 
fact that under $H_0^*$, $\tilde{\vv{X}}$ has independent columns.
Let $\tilde{\vv{x}}_j$ be  the $j$th column of $\tilde{\vv{X}}$. Note that under the Kronecker product structure, 
we have ${\rm cov}(\tilde{\vv{x}}_j,\tilde{\vv{x}}_k)=(\vv{\Sigma}_C)_{jk}\vv{\Sigma}_R$ for any $j,k\in\{1,\ldots,q\}$. Thus, under $H_0^*$,
\[
r_{jk}:=\frac{1}{p}{\rm tr}\{{\rm cov}(\tilde{\vv{x}}_j,\tilde{\vv{x}}_k)\}
=\frac{1}{p}\E[\tilde{\vv{x}}^\T_k\tilde{\vv{x}}_j] 
=\frac{1}{p}(\vv{\Sigma}_C)_{kj} {\rm tr}(\vv{\Sigma}_R) =0,
\]
for any $j\neq k$.
In other words, the between column covariance  matrix $\vv{R}=(r_{jk})_{1\leq j,k\leq q}$ is diagonal.
Based on this observation,
\cite{TA2019} proposed a test statistic $W_n=T_{2n}-T_{3n}$  as an unbiased estimator for the sum of squares of the off-diagonal elements of  $\vv{R}$, namely,  
\[
\E[W_n] = {\rm tr}(\vv{R}^2)- {\rm tr}(\vv{R}\circ \vv{R}) = \sum_{j\neq k } r_{jk}^2,
\]
where $\circ$ stands for the Hadamard product of matrices,
\begin{align*}
T_{2n}&=\frac{1}{p^2n(n-1)}\sum_{i\neq j}{\rm tr}(\tilde{\vv{X}}^\T_i\tilde{\vv{X}}_i\tilde{\vv{X}}^\T_j\tilde{\vv{X}}_j)
-\frac{2}{p^2\binom{n}{3}}\sum_{(i,j,k)}{\rm tr}(\tilde{\vv{X}}^\T_i\tilde{\vv{X}}_i\tilde{\vv{X}}^\T_j\tilde{\vv{X}}_k)\\
&+\frac{1}{p^2\binom{n}{4}}\sum_{(i,j,k,l)}{\rm tr}(\tilde{\vv{X}}^\T_i\tilde{\vv{X}}_j\tilde{\vv{X}}^\T_k\tilde{\vv{X}}_l), \mbox{and}\\
T_{3n}&=\frac{1}{p^2n(n-1)}\sum_{i\neq j}{\rm tr}[(\tilde{\vv{X}}^\T_i\tilde{\vv{X}}_i)\circ(\tilde{\vv{X}}^\T_j\tilde{\vv{X}}_j)]-\frac{2}{p^2\binom{n}{3}}\sum_{(i,j,k)}{\rm tr}[(\tilde{\vv{X}}^\T_i\tilde{\vv{X}}_i)\circ(\vv{X}^\T_j\vv{X}_k)]\\
&+\frac{1}{p^2\binom{n}{4}}\sum_{(i,j,k,l)}{\rm tr}[(\tilde{\vv{X}}^\T_i\tilde{\vv{X}}_j)\circ(\tilde{\vv{X}}^\T_k\tilde{\vv{X}}_l)].
\end{align*}

\begin{figure}
\includegraphics[width = 80 mm, height = 60 mm]{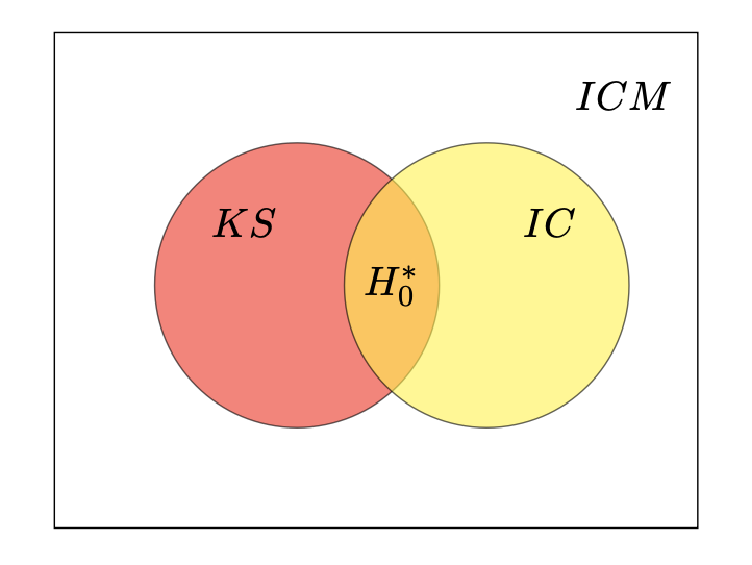}
\caption{The relationship among different types of transposable data. $ICM =$ the independent component model. $KS =$ the data with the Kronecker product dependence structure. $IC =$ the data with independent columns. $H_0^*=$ the data satisfies the hypothesis $H_0^*$ \eqref{eq:hykron}.}
\label{fig:venn}
\end{figure}

Although both our procedure and 
the test in \cite{TA2019} are valid   for testing the hypothesis $H_0^*$,  the two  tests  have different focuses. Our test statistic focuses on detection of non-proportionality of column covariance matrices while the test statistic $W_n$ is more sensible to detect dependence between columns. The Venn diagram in Figure~\ref{fig:venn} helps further illustrate this difference.

Our procedure is more sensible of alternatives lying in the yellow region where the transposable data has independent columns whose covariance matrices are not mutually proportional (this implies that the data does not follow the Kronecker product structure).  
In this region, the procedure in \cite{TA2019}  
might accept  $H_0^*$ since $W_n$ is likely to have low values for independent columns.
On the other hand, the procedure with $W_n$ is more powerful for the detection of alternatives in the red region where the transposable data does have the Kronecker product structure with however dependent columns. Indeed in this situation, our test statistic based on non-proportionality between column covariance matrices will show a low value and is more likely to accept $H_0^*$ because the column covariance matrices are still proportional under the Kronecker product structure. 
To summarize, for testing the hypothesis $H_0^*$,    
if one believes that the covariance matrices of columns are not proportional, then our procedure is preferable. Otherwise,   the proportionality of column covariance matrices is accepted and the procedure in \cite{TA2019} is preferable.  

To support our assertions, we conduct a simulation study under the following setting: $p = 100$, $q = 50$, and the common sample size $n = 150$. The following two scenarios are considered.
\begin{itemize}
\item \textbf{Case I}: The $p\times q$ matrix-valued population $\vv{X} = (\vv{x}_1,\ldots,\vv{x}_q)$ consists of $q$ mutually independent columns with $\vv{x}_i=\vv{\Sigma}_i^{1/2}\vv{z}_i$ as in Assumption \ref{assm:2.1}. The population covariance matrix $\vv{\Sigma}_i = w_i[(1-\beta)\vv{\Lambda}_0 + \beta\vv{\Lambda}_i]$, in which, for $j=0,\ldots,q$, the basis $\vv{\Lambda}_j = \vv{U}^\top_j \vv{D}_j \vv{U}_j$ with $\vv{U}_j$ a random orthogonal matrix and  $\vv{D}_j$ a diagonal matrix with entries randomly picked from $(e^{-3},e^3)$ such that $\vv{\Lambda}_0,\ldots,\vv{\Lambda}_q$ are non-proportional, and $\{w_j\}$ are randomly selected positive weights in $(0.5, 1.5)$. When $\beta = 0$, the $q$ populations are proportional to $\vv{\Lambda}_0$ so $H_0$ holds. When $\beta$ deviates from $0$, the $q$ populations become non-proportional.

\item \textbf{Case II}: The $p\times q$ matrix-valued population $\vv{X} = (\vv{x}_1, \ldots, \vv{x}_q)$ satisfies the Kronecker product dependence structure $\vv{X} = \vv{\Sigma}_R^{1/2}\vv{Z}\vv{\Sigma}_C^{1/2}$, where $\vv{\Sigma}_C^{1/2} = \mathrm{diag}(\vv{w}) + \beta(\vv{e}_1\vv{e}_2^\top + \vv{e}_2\vv{e}_1^\top)$ with $\vv{w} = \{w_j\}$ positive weights in $(0.5,1.5)$ and $\vv{e}_j$ a $q$-dimensional vector having its $j$th position one and rest positions zero, $j=1,2$. When $\beta = 0$, $\vv{\Sigma}_C^{1/2}$ is diagonal, implying that all these $q$ columns are mutually independent with their covariance matrices proportional to each other. When $\beta > 0$, columns become dependent but keep the column covariance matrices mutually proportional.
\end{itemize}
In summary, Case I maintains column independence but alters the proportionality of column covariance matrices, whereas Case II changes column dependence but keeps the proportionality of column covariance matrices unchanged.
Besides, two different types of noises are also considered: the normal noise $z_{ij}\sim\mathcal{N}(0,1)$ and the gamma noise $z_{ij}\sim {\sf Gamma}(4,2)-2$. 

\begin{figure}
\centering
\includegraphics[width = 43 mm, height = 43 mm]{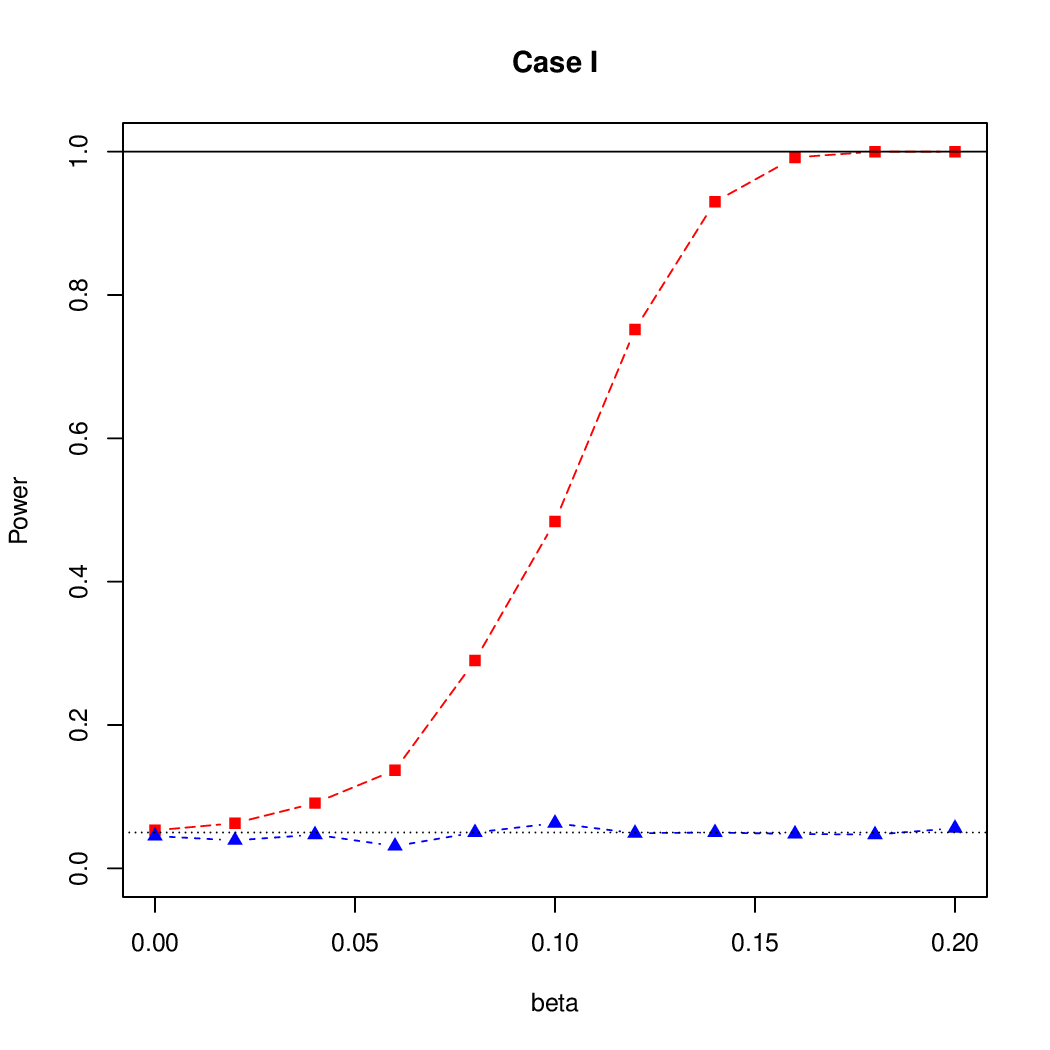}
\includegraphics[width = 43 mm, height = 43 mm]{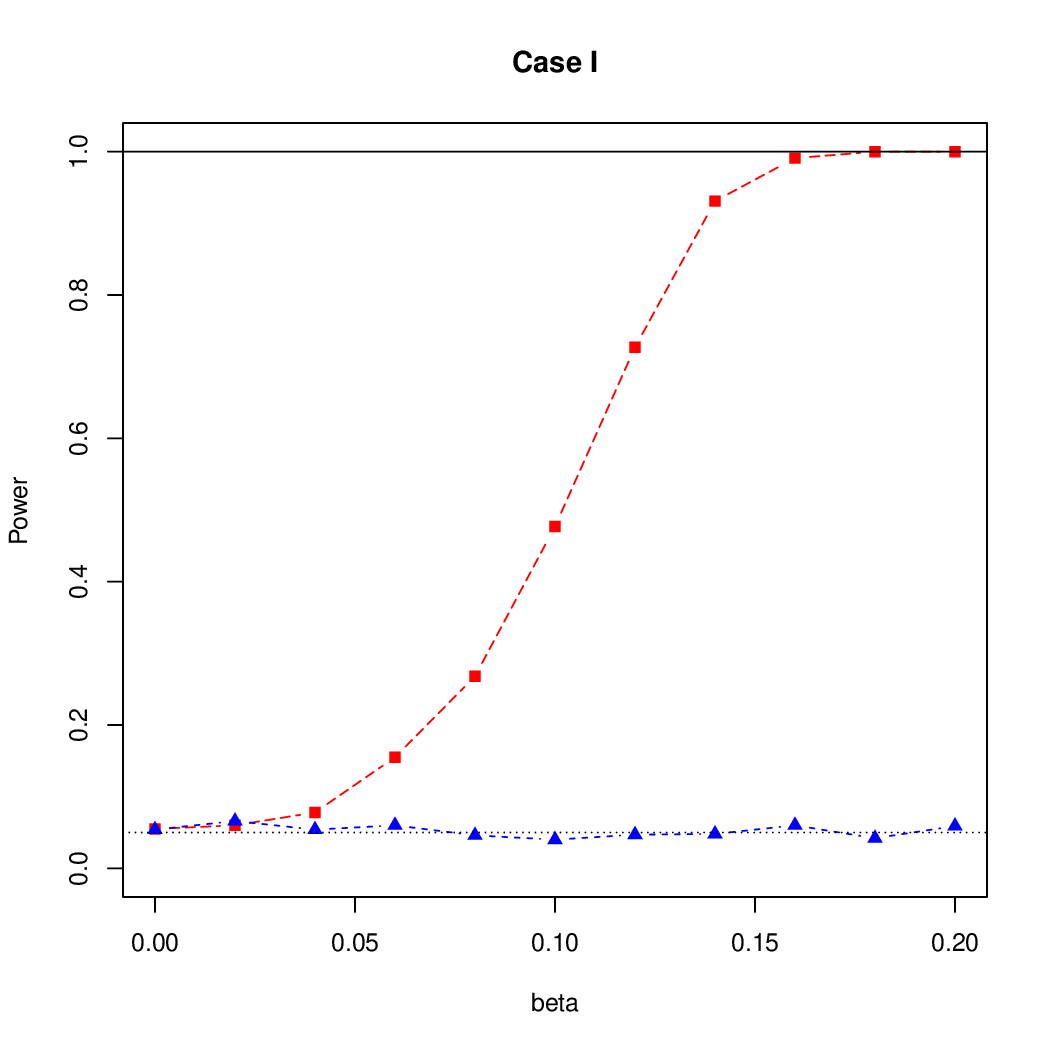}
\includegraphics[width = 43 mm, height = 43 mm]{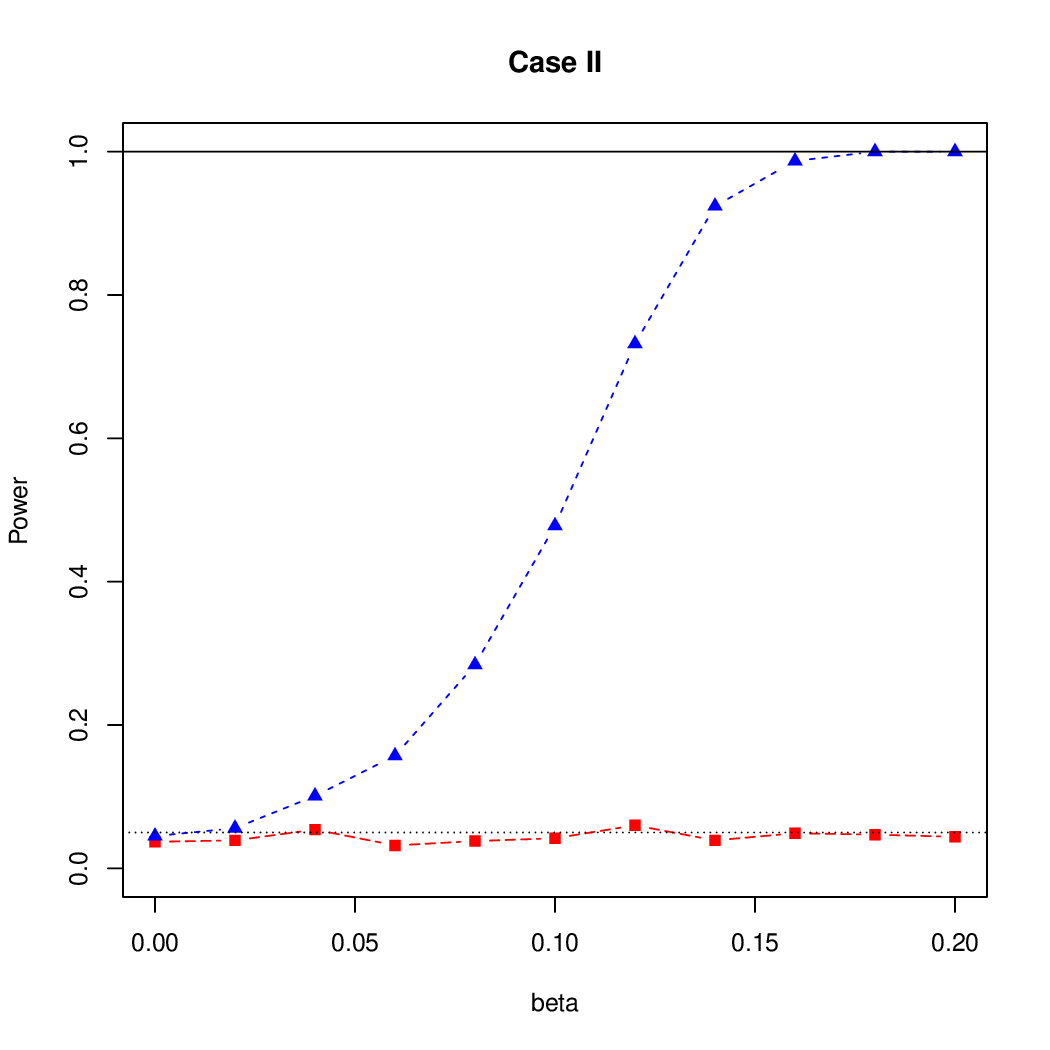}
\includegraphics[width = 43 mm, height = 43 mm]{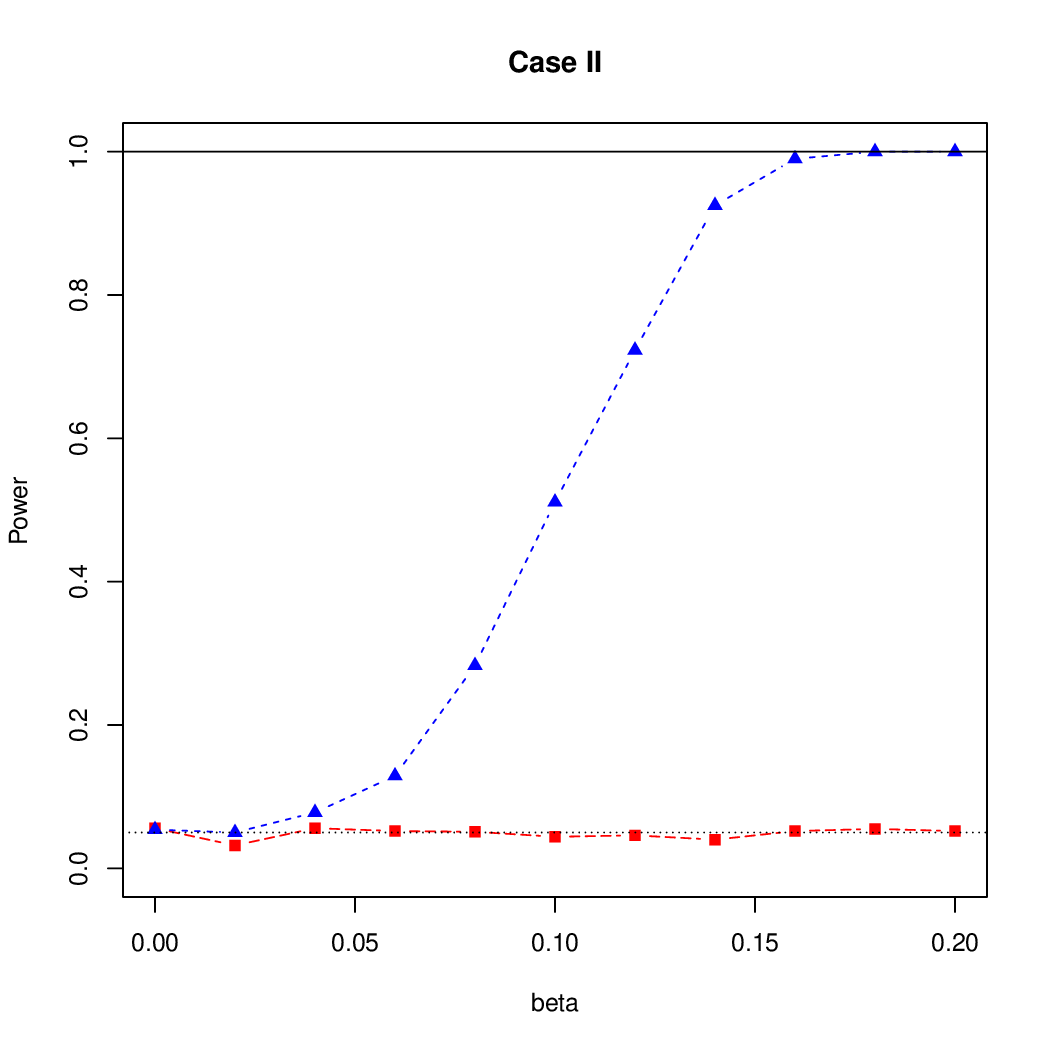}
\caption{Empirical power plot of the comparison of the many-sample proportionality test (red square) and 
the diagonality test (blue triangle) in \cite{TA2019} 
under two cases, I and II, 
with two different noises, normal noise (left panel) and Gamma noise (right panel).}
\label{fig:add}
\end{figure}

Under these two scenarios with two different types of noise, we computed both our test statistic and the one in \cite{TA2019} over 1000 iterations and the empirical power curves plotted in Figure \ref{fig:add}. In Case I, our test shows significant power, while the power of the test in \cite{TA2019} fluctuates around 0.05 even when the alternative does not hold. This indicates the sensitivity of the many-sample proportionality test to covariance proportionality.
In Case II, our test is not sensitive to the variations in column dependence, keeping its power around 0.05. However, the test in \cite{TA2019} efficiently detects the change in column dependence and demonstrates significant power. These results support our previously mentioned assertions.

\subsection{Analysis of the mouse aging data}\label{sec:5.1}
Following \cite{TA2019}, we analyze the gene dataset of Mouse Aging Project collected in \cite{mice07}, which measures gene expression levels from $p=46$ genes in up to $16$ tissues for $40$ mice ($n=40$). \cite{TA2019} examined the dependence structure of genes expression levels among $9$ tissues ($q=9$) by testing the diagonality of $\Sigma_C$, assuming the Kronecker product dependence
structure (\ref{eq:kron}). Their test accepted the diagonality hypothesis $H_0^*$ (\ref{eq:hykron}). 
However, applying our specification test  presented in Section~\ref{ssec:2.6} to this hypothesis  $H_0^*$  leads to the statistics value
$\sqrt{q}U_p/\hat{{\sigma}}_p= 13.592$, which is much larger than the $95\%$ standard normal quantile; so the hypothesis $H_0^*$ is strongly rejected. Therefore, the two tests lead to significantly contradictory conclusions!

To elucidate this, according to the discussion in Section \ref{subsec:2.4}, we see that the procedure in \cite{TA2019} essentially tests the hypothesis that the $q=9$ columns from the tissues are independent; so their conclusion is indeed to accept the independence of the columns. 
According to Lemma \ref{lem:kron}, if, in addition, the Kronecker product structure  (\ref{eq:kron}) is also true,
then the $q$ column covariance matrices must be proportional, which is a fact strongly rejected by our test. 
Therefore,  the only plausible explanation is that 
the $q=9$ tissue populations are independent while the whole data set {
does not follow} the Kronecker product structure  (\ref{eq:kron}).  
 To recap, we conclude that gene expression levels among tissues are independent but their covariance matrices are not proportional, and the Kronecker product covariance model   (\ref{eq:kron}) is not appropriate for this mouse aging data set although the model is widely used in many studies about the dataset.

\begin{remark}
One may wonder that the number of mouse groups is $q=9$ (number of samples) which is not very large and that our asymptotic results for the specification test might not be accurate. 
Following the recent literature in high-dimensional statistics  \citep{Bai09,yao15}, high-dimensional effects take place if the ratio of dimension-to-sample-size is away from zero even though both the dimension and the sample size are not very large.  In our situation,  the effects of many samples can be visible if the ratio $q/n$ of sample number to sample size is away from zero: here we have $q/n = 9/40 = 0.225$ indeed.  It can be expected that the asymptotic results assuming fixed $q$ might be less accurate than our many-sample asymptotic results that consider $q$ and $n$ growing to infinity in comparable magnitude.  

To further verify that our specification test is already accurate enough even for relatively small population number and sample sizes, we conduct a small simulation experiment with dimensions and sample sizes matching those of the mouse aging data, namely,  $p=46$, $q=9$ and $n_1=\cdots =n_{9}=40$, where other design parameters follow those of the experiments (a) and (b)  in Section \ref{ssec:2.5}. 
The empirical sizes of the test under different settings are reported in Table~\ref{tab:2}.
It can be observed that our proposed test controls the size well under the null.  This thus confirms the accuracy of the asymptotic null distribution of the test statistic we derived under the many-sample framework even though the sample number and the sample sizes are not very large.

\begin{table}
\caption{Empirical sizes of the many-sample proportionality test, with $p,q,n$ matching those in \cite{TA2019}.} 
\label{tab:2}
\begin{tabular}{@{}cccc@{}}
\hline
\multicolumn{2}{c}{Case (a)}&
\multicolumn{2}{c}{Case (b)}
 \\ [3pt]
\hline
Normal & Gamma & 
Normal & Gamma
 \\
\hline
{$0.051$} & {$0.046$} &
{$0.054$} & {$0.053$}
 \\
\hline
\end{tabular}
\end{table}

\end{remark}

\begin{remark}

For the case of a small number of populations ($q = 9$), it is of special interest to compare our proposed method with existing methods for a fixed number of populations. Here, we make the comparison with a recent study \cite{RAHMAD2022104865} that addresses multi-sample tests for the proportionality of large-dimensional covariance matrices. In \cite{RAHMAD2022104865}, the author proposed a $(q - 1)$-dimensional random vector $\vv{T}_0$ to characterize the mutual proportionality among groups and proved that the null distribution of $\vv{T}_0$ asymptotically follows a $(q-1)$-dimensional standard normal distribution. Thus, the inner product $t_0 := \vv{T}_0^\top \vv{T}_0$ converges to a chi-square distribution $\chi^2_{q-1}$ with $q-1$ degrees of freedom under the null.

In the simulation, we evaluate both size control and power by comparing our proportionality test with the method in \cite{RAHMAD2022104865} under the same settings: $p = 46, q = 9,$ and $n = 40$ as in the mouse aging data. We keep $q$ populations $\vv{x}_1,\ldots,\vv{x}_q$ to be mutually independent and the $i$th population satisfies the structure $\vv{x}_i=\vv{\Sigma}_i^{1/2}\vv{z}_i$ as in Assumption \ref{assm:2.1}, with two different types of noises $z_{ij}$: $\mathcal{N}(0,1)$ and ${\sf Gamma}(4,2)-2$. 
For $i=1,\ldots,q$, the population covariance matrix $\vv{\Sigma}_i = w_i[(1-\beta)\vv{\Lambda}_0 + \beta\vv{\Lambda}_i]$, where for $j=0,\ldots,q, \vv{\Lambda}_j = \vv{U}^\top_j \vv{D}_j \vv{U}_j$ with $\vv{U}_j$ a random orthogonal matrix and  $\vv{D}_j$ a diagonal matrix with entries randomly picked from $(e^{-3},e^3)$ such that $\vv{\Lambda}_0,\ldots,\vv{\Lambda}_q$ are not proportional to each other, and $\{w_j\}$ are randomly selected positive weights in $(0.5, 1.5)$.  
It can be easily seen that the covariance matrices of these $q$ populations are proportional to $\vv{\Lambda}_0$ when $\beta = 0$. As $\beta$ increases, the $q$ covariance matrices become non-proportional.

Under both the normal and the gamma noise settings, we compute both our test statistic and the statistic $t_0 =\vv{T}_0^\top \vv{T}_0$ in \cite{RAHMAD2022104865} repeatedly for 1000 iterations and derive the empirical sizes table in Table \ref{tab:compare} and the empirical power curves plotted in Figure \ref{fig:compare}.
The results demonstrate that, although the number of populations is small, our test effectively controls the size, whereas the method in \cite{RAHMAD2022104865} fails to do so. Additionally, our test procedure shows a faster increase in power towards one as it deviates from the null hypothesis compared to the method in \cite{RAHMAD2022104865}. These findings validate the efficiency of our many-sample test procedure.
\begin{table}
\caption{Empirical sizes of the many-sample proportionality test and the \cite{RAHMAD2022104865}'s test, with $p,q,n$ matching those in \cite{TA2019}.} 
\label{tab:compare}
\begin{tabular}{@{}cc|cc@{}}
\hline
\hline
\multicolumn{2}{c|}{many-sample test}&
\multicolumn{2}{c}{proportionality test in \cite{RAHMAD2022104865}}
 \\ [3pt]
\hline
Normal & Gamma & 
Normal & Gamma
 \\
\hline
{$0.047$} & {$0.049$} &
{$0.112$} & {$0.121$}
 \\
\hline
\end{tabular}
\end{table}

\begin{figure}
\centering
\includegraphics[width = 43 mm, height = 43 mm]{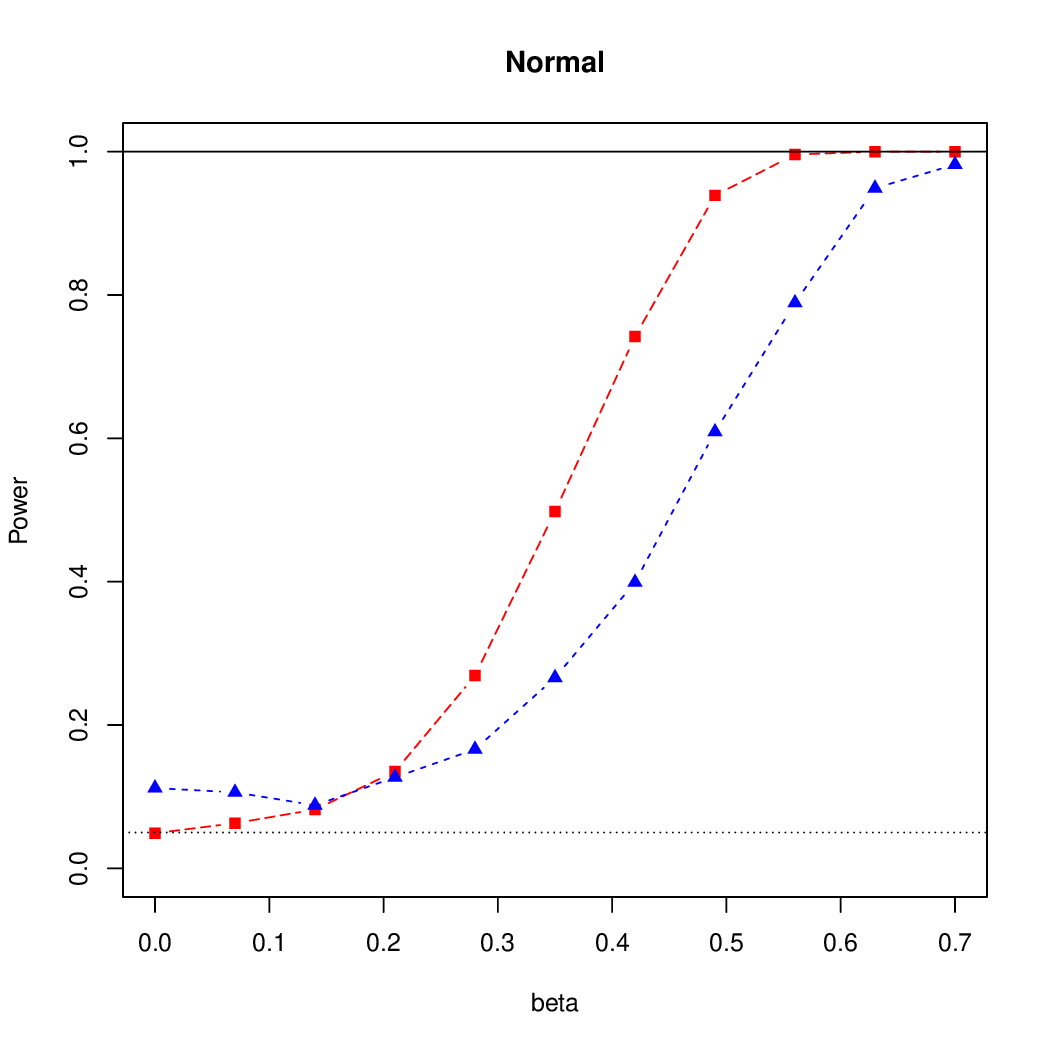}
\includegraphics[width = 43 mm, height = 43 mm]{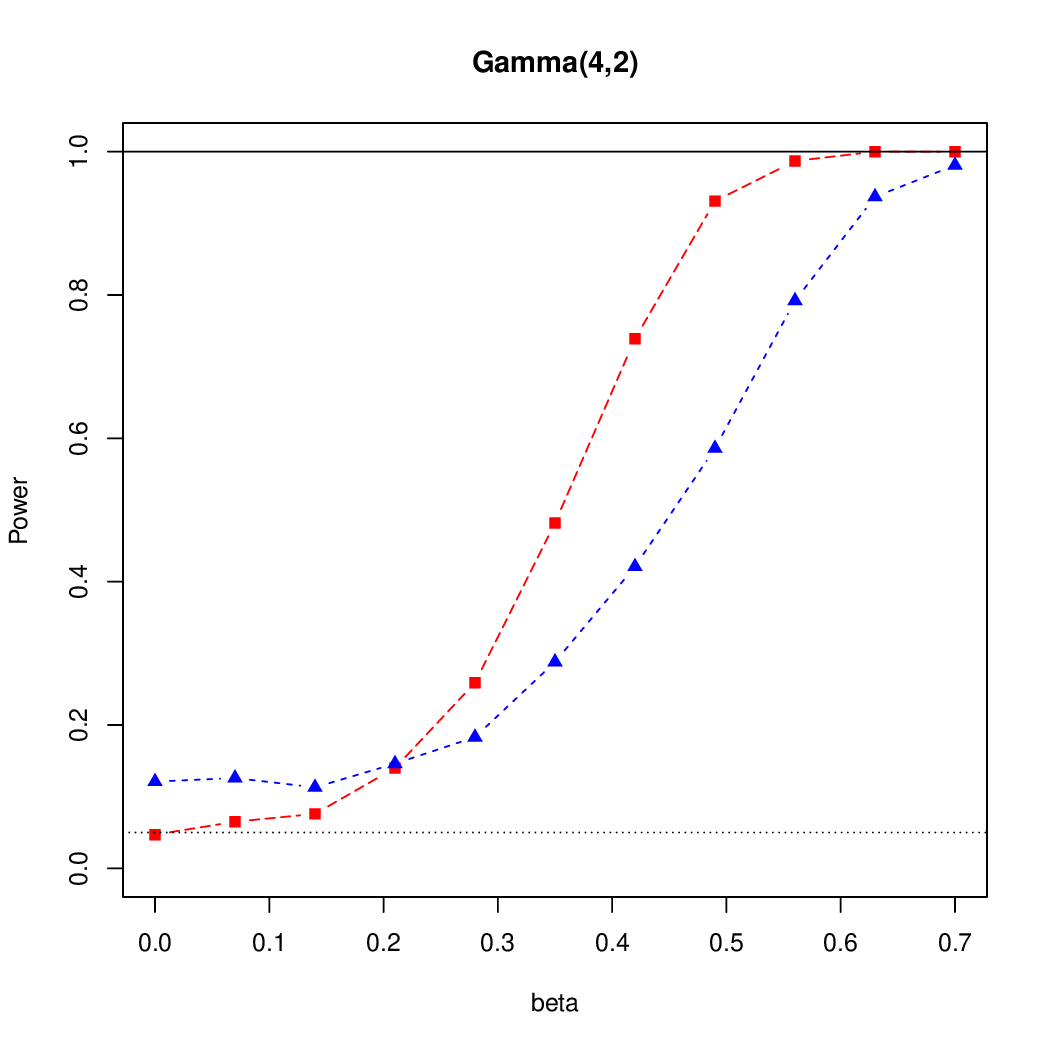}
\caption{Empirical power plot of the comparison of the many-sample proportionality test (red square) and 
the proportionality test in \cite{RAHMAD2022104865} (blue triangle) 
with two different noises, normal noise (left) and Gamma noise (right).}
\label{fig:compare}
\end{figure}

\end{remark}

\section{Testing of equality}\label{sec:3}

\subsection{Basic settings}
Let $\vv{x}_1,\ldots,\vv{x}_q$ be $q$  $p$-dimensional random vectors with covariance matrix $\vv{\Sigma}_i$, $i=1,\ldots,q$. Consider the following equality test:
\begin{equation*}\label{eq:3.1}
H_0:~~~\vv{\Sigma}_i=\vv{\Sigma},~~~i=1,\ldots,q, 
\end{equation*}
for some unknown non-negative definite
matrix $\vv{\Sigma}$. For two $p\times p$ matrices $\vv{A}$ and $\vv{B}$, define
$d_0(\vv{A},\vv{B})={\rm tr}[(\vv{A}-\vv{B})^2]$.
It is easy to see that $d_0$ characterizes the equality of two matrices. 
{
Similar to Section \ref{sec:2}, we define
\begin{equation}\label{eq:m_p}
    m_p = \binom{q}{2}^{-1}\sum_{i<j} d_0(\vv{\Sigma}_i,\vv{\Sigma}_j).
\end{equation}
Then, the testing problem can be equivalently transferred to
\begin{equation}\label{eq:alt_H0}
    H_0:~~~m_p = 0~~~~\hbox{versus}~~~~H_1:~~~m_p>0.
\end{equation}
We construct the test statistic $V_p$ in (\ref{eq:3.20}) based on an unbiased estimation of the distances $d_0(\vv{\Sigma}_i,\vv{\Sigma}_j)$'s and the quantity $m_p$.
}

Testing the equality hypothesis (\ref{eq:alt_H0}) is, in fact, simpler than testing the proportionality hypothesis analyzed in the previous sections.  By repeating a similar analysis, we can construct a test statistic for the many-sample testing of the equality hypothesis. The same Assumptions \ref{assm:2.4}---\ref{assm:2.1} given in Section \ref{sec:2} will be used. 
For each $i=1,\ldots,q$, assume that there are $n_i$ independent observations $\vv{x}_{i,1},\ldots,\vv{x}_{i,n_i}$ from population $\vv{x}_i$. Let $\vv{X}_{ip}=(\vv{x}_{i,1},\ldots,\vv{x}_{i,n_i})$ be the data matrix and the sample covariance matrix $\vv{S}_{i,p}$ defined as in (\ref{eq:2.4}).

\subsection{Unbiased estimation of the distance}
\label{subsec:3.1}
The subsection parallels Section \ref{subsec:2.1}. Hence, we omit the details and directly state the final result.

{
Recall that, for any $i\neq j$,
\begin{equation*}\label{eq:d0_expan}
    d_0(\vv{\Sigma}_i,\vv{\Sigma}_j) = p\left\{\mu_{i,2} + \mu_{j,2} - 2{\rm tr}(\vv{\Sigma}_i\vv{\Sigma}_j)\right\},
\end{equation*}
where $\mu_{i,2} = p^{-1}{\rm tr}(\vv{\Sigma}_i^2)$ for $i=1,\ldots,q$. It should be noted that the second-order spectral moment of the sample covariance matrix is biased to that of the corresponding population covariance matrix due to the high-dimensional effects $p/n_i\to c_i\in(0,\infty)$. Hence,
we apply the unbiased estimator $\hat{\mu}_{i,2}$ in (\ref{eq:ub_est_mu}) and define
\begin{equation}\label{eq:ub_mp}
    g(\vv{X}_{ip},\vv{X}_{jp}) = p\left\{\hat{\mu}_{i,2} + \hat{\mu}_{j,2} - 2\frac{1}{p}{\rm tr}(\vv{S}_{ip}\vv{S}_{jp})\right\}
\end{equation}
It is easy to see from the discussion in Lemma \ref{lem:ub_est} and the fact $\E[{\rm tr}(\vv{S}_{ip}\vv{S}_{jp})] = {\rm tr}(\vv{\Sigma}_i,\vv{\Sigma}_j)$ that $g(\vv{X}_{ip},\vv{X}_{jp})$ is indeed an unbiased estimator for the distance $d_{0}(\vv{\Sigma}_i,\vv{\Sigma}_j)$.
}

\subsection{The test statistic and its asymptotic normality}\label{subsec:3.2}
The proposed test statistic is given as
\begin{equation*}\label{eq:3.20}
V_p=\binom{q}{2}^{-1}\sum_{i<j} g(\vv{X}_{ip},\vv{X}_{jp}).
\end{equation*}
{An immediate consequence of the discussion in Section \ref{subsec:3.1} is that $V_p$ is an unbiased estimator of the quantity $m_p$, i.e., $\E[V_p] = m_p$. In particular, when $H_0$ holds, all distances $d_{0}(\vv{\Sigma}_i,\vv{\Sigma}_j) =0$ and this $\E[V_p] = m_p =0$.

To optimize the computational complexities, an explicit representation of $V_p$ is also given as follows:
\begin{equation}\label{eq:alt_mp}
    \begin{split}
        V_p & = \frac{1}{q(q-1)}\sum_{i\neq j} p\{\hat{\mu}_{i,2} + \hat{\mu}_{j,2} - 2\frac{1}{p}{\rm tr}(\vv{S}_{ip}\vv{S}_{jp})\} \\
        & = \frac{2}{q}\sum_{i=1}^q p\hat{\mu}_{i,2} - \frac{2q}{(q-1)}{\rm tr}\left(\frac{1}{q}\sum_{i=1}^q \vv{S}_{ip}\right)^2 + \frac{2}{q-1}\left(\frac{1}{q}\sum_{i=1}^q {\rm tr}(\vv{S}_{ip}^2)\right).
    \end{split}
\end{equation}
It can be seen from the expression that the summations over two distinct indices are reduced, and the final result only involves mean values of $q$ numbers or matrices. So, the computational time is significantly shortened.

The following theorem reveals the asymptotic normality of our proposed statistic $V_p$ under both the null and the alternative.

\begin{theorem}\label{thm:3.3}
Let 
\begin{equation}\label{eq:tLambdai}
    \tilde{\vv{\Gamma}}_{i,p} = \vv{\Sigma}_i - \frac{1}{q-1}\sum_{j:j\neq i} \vv{\Sigma}_j,~~~~i=1,\ldots,q.
\end{equation}
Then,
under  Assumptions \ref{assm:2.1}---\ref{assm:2.4}, 
it holds that
\[
\frac{\sqrt{q}(V_p - m_p)}{\lambda_p} \overset{d.}{\to}\mathcal{N}(0,1).
\]
The asymptotic variance $\lambda_p^2 = \frac{4}{q}\sum_{i=1}^q \lambda_{i,p}^2$ with

\begin{equation*}
    \lambda_{i,p}^2 = 4c_{ip}^2\mu_{i,2}^2 + 4c_{ip}\langle \tilde{\vv{\Gamma}}_{i,p},\tilde{\vv{\Gamma}}_{i,p}\rangle_{\vv{\Sigma}_i},
\end{equation*} 
where the inner product $\langle\cdot,\cdot\rangle_{\vv{\Sigma}_i}$ is defined as in (\ref{eq:inner_prod}).
\end{theorem}
The proof of the theorem is provided in Appendix \ref{app:2}. From the above theorem,
we see that the asymptotic variance can be decomposed into two parts $\lambda_p^2 = \lambda_{0,p}^2 + \lambda_{r,p}^2$, in which
\begin{equation}\label{eq:var_eq_decomp}
    \lambda_{0,p}^2 = \frac{16}{q}\sum_{i=1}^q c_{ip}^2\mu_{i,2}^2,~~~~~\lambda_{r,p}^2 = \frac{16}{q}\sum_{i=1}^q c_{ip}\langle \tilde{\vv{\Gamma}}_{i,p},\tilde{\vv{\Gamma}}_{i,p}\rangle_{\vv{\Sigma}_i}.
\end{equation}
Under Assumptions \ref{assm:2.1}---\ref{assm:2.3}, $\lambda_{0,p}^2$ is always positive and has both upper and lower bounds. The remainder $\lambda_{r,p}^2$ remains non-negative. In particular, when $H_0$ holds, we see directly from (\ref{eq:tLambdai}) that $\tilde{\vv{\Gamma}}_{i,p} = \vv{O}$ for all $i=1,\ldots,q$ and thus $\lambda_{r,p}^2$ vanishes. Thus, as a direct consequence of Theorem \ref{thm:3.3}, we have
\begin{corollary}\label{cor:eq_test}
    Suppose that $H_0$ in (\ref{eq:alt_H0}) holds. Then, under Assumptions \ref{assm:2.1}---\ref{assm:2.4}, 
    \[
    \frac{\sqrt{q}V_p}{\lambda_{0,p}}\overset{d.}{\to}\mathcal{N}(0,1).
    \]
\end{corollary}
}

To estimate the asymptotic variance $\lambda_{0,p}^2$, we define
\begin{equation*}\label{eq:3.9}
\hat{\lambda}_{p}^2=\frac{16}{q}\sum_{i=1}^q c_{ip}^2\hat{\mu}_{i,2}^2.
\end{equation*}

\begin{theorem}\label{thm:3.4}
Under Assumptions \ref{assm:2.1}---\ref{assm:2.4}, $\hat{\lambda}_{p}^2$ is always consistent to $\lambda_{0,p}^2$ . Consequently, under $H_0$, it holds that 
\[
\frac{\sqrt{q}V_p}{\hat{\lambda}_{p}}\overset{d.}{\to}\mathcal{N}(0,1).
\]
\end{theorem}
The detailed proof of Theorem \ref{thm:3.4} is given in 
Appendix 
\ref{sec:proof_thm3.5}.

In summary, the procedure for testing  $H_0$ in (\ref{eq:alt_H0}) at significance level $\alpha$ is to reject $H_0$ if $q^{1/2}V_p/\hat{\lambda}_p>z_\alpha$,
where $z_{\alpha}$ is the $\alpha$th upper-quantile of the standard normal. 
\subsection{Power under alternative hypothesis}\label{subsec:3.3}

{

This section is devoted to analyzing the power of our proposed test procedure under the alternative $H_1$ in (\ref{eq:alt_H0}).  The next theorem reveals the asymptotic representation of the power function and a sufficient and necessary condition for a strong rejection.
\begin{theorem}\label{thm:generalpower_eq}
    Assume that Assumptions \ref{assm:2.1}---\ref{assm:2.4} holds. For any significant level $\alpha\in(0,1)$, 
    the power function satisfies
    \begin{equation}\label{eq:pf_eq}
        \P_{H_1}\left(\frac{\sqrt{q}V_p}{\hat{\lambda}_{p}} > z_\alpha\right)
        =\Phi\left(\frac{\sqrt{q}m_p}{\lambda_p} - \frac{\lambda_{0,p}}{\lambda_p}z_\alpha\right) + o(1).
    \end{equation}
    Consequently,
    the power function tends to $1$ as $p,q,\{n_j\}$ increase if and only if  $\sqrt{q}m_p \to \infty$. 
\end{theorem}
The proof of the theorem is mainly based on the previous Theorems \ref{thm:3.3} and \ref{thm:3.4} and is provided in Appendix \ref{sec:proof_generalpower_eq}. 
\begin{remark}
Similar to Remark \ref{rem:power_nece}, we define $d'_{\max,q} = \max_{i<j} d_0(\vv{\Sigma}_i,\vv{\Sigma}_j)$. Then, the testing power vanishes if $d'_{\max,q} = o(q^{-1/2})$. In other words, a necessary condition to ensure the power is that the convergence rate of $d'_{\max,q}$ to $0$ cannot be faster than $q^{-1/2}$. 
\end{remark}
}

{
Theorem \ref{thm:generalpower_eq} that the divergence of the mean drift $\sqrt{q}m_p$ guarantees the power. 
Thus, we investigate examples parallel to those in Section \ref{subsec:2.3} to get respective sufficient conditions. 
\begin{example}[$q$ matrices without equal pair]\label{example:1'}
     Suppose that  for any $1\leq i\neq j \leq q$, $d_{\rm prop}(\vv{\Sigma}_i,\vv{\Sigma}_j)>0$. Let
    $d'_{1,q} := \min_{i<j} d_0(\vv{\Sigma}_i,\vv{\Sigma}_j)$.
    Then, the converge rate of $d'_{1,q}$ influences the power. Since $\sqrt{q}M_p \geq \sqrt{q}d_{1,q}$,
    according to Theorem \ref{thm:generalpower_eq}, the power function tends to $1$ when $\sqrt{q}d_{1,q}\to \infty$. 
    \hfill$\Box$
\end{example}

The next example focuses on the case where the $q$ populations can be divided into $m$ subgroups, where populations within the same group are equal to each other while populations from different groups are not. 
\begin{example}[Unequal subgroups without dominant group]\label{example:2'}
    Let $m = m(q)$ be a positive integer smaller than $q$ and may vary as $q$ increases.
    Suppose that $\{1,\ldots,q\}$ can be split into $m$ non-empty disjoint subsets $I_1,\ldots,I_m$.  
    Furthermore, we assume that there are $m$ distinct $p\times p$ non-negative definite matrices $\vv{\Lambda}'_1,\ldots,\vv{\Lambda}'_m$ such that $\vv{\Sigma}_r= \vv{\Lambda}'_i$ if $r\in I_i$. We denote $d'_{2,q} := \min_{1\leq i<j\leq m} d_0(\vv{\Lambda}'_i,\vv{\Lambda}'_j)$.
    Still, $d'_{2,q}$  is positive for fixed $q$, but it can decrease to zero as $q$ varies. Same as Example \ref{example:2}, we use $q_i$ stands for the number of elements in $I_i$, and assume that condition (\ref{eq:theta_0}) holds, which suggests
    It implies that none of the $m$ groups will be dominant. 

    We show that the divergence of the mean drift is determined by $d'_{2,q}$ and $\theta_0$.
    Observe that
    \begin{align*}
        m_p & = \frac{1}{q(q-1)} \sum_{r\neq s} d_{0}(\vv{\Sigma}_r,\vv{\Sigma}_s) \\
        & = \frac{1}{q(q-1)}\sum_{i\neq j} \sum_{r\in I_i,s\in I_j} d_0(\vv{\Lambda}'_i,\vv{\Lambda}'_j)\\
        & = \frac{1}{q(q-1)}\sum_{i\neq j} d_0(\vv{\Lambda}'_i,\vv{\Lambda}'_j)q_iq_j\\
        & \geq  d'_{2,q} \sum_{i=1}^m \frac{q_i}{q}\left(1-\frac{q_i}{q}\right) \\
        & \geq d'_{2,q} \left(1- \theta_0\right).
    \end{align*}
    Thus, if $\sqrt{q}d'_{2,q} \to \infty$, the value of our power function will tend to $1$. 
    \hfill$\Box$
\end{example}

The next example considers the scenario with $\theta_0 =1$, in which the total $q$ populations consist of a dominant subgroup with a small number of outliers.

\begin{example}[One dominant group with a few outliers]\label{example:3'}
    Suppose that there exists a non-empty subset $I_0 \subset\{1,\ldots,q\}$ and an un-specified basis matrix $\vv{\Sigma}$ such that populations outside $I_0$ all equal to $\vv{\Sigma}$, but within the subset, populations are not. We denote $
        d'_{3,q}= \min_{j\in I_0} d_0(\vv{\Sigma}_j,\vv{\Sigma})$.
    Here, $d'_{3,q}$ is positive for any $q$ fixed but may decrease to zero as $q$ grows.
    
     We still focus on the small-$K$ setting: $K = K(q) = o(q)$, where $K = |I_0|$ is the number of elements in $I_0$. This assumption causes that the term $\lambda_{r,p}^2 = o(1)$ in (\ref{eq:var_eq_decomp}) and thus $\lambda_p^2 = \lambda_{0,p}^2 +o(1)$. It then follows that (\ref{eq:pf_eq}) can be further simplified as
    \begin{equation}\label{eq:alterpf_eq}
        \P_{H_1}\left(\frac{\sqrt{q}V_p}{\hat{\lambda}_{p}} > z_\alpha\right)
        =\Phi\left(\frac{\sqrt{q}m_p}{\lambda_{0,p}} - z_\alpha\right) + o(1).
    \end{equation}
    
   We denote $\mu_2 = p^{-1}{\rm tr}(\vv{\Sigma}^2)$ and $
    \bar{c}_{q,12}:=(q-K)^{-1}\sum_{i\in I_0^c}c_{ip}^2$.
    Under Assumptions \ref{assm:2.2}, \ref{assm:2.3} and the small-$K$ condition, it is easy to show that both $\mu_2$ and $\bar{c}_{q,12}$ have positive upper and lower bounds.
    The asymptotic variance can be simplified as
    $
    \lambda_{0,p}^2 = 16 \mu_2^2\bar{c}_{q,12}$.

Observe that
\begin{equation*}
\begin{split}
    m_p &=\frac{2}{q(q-1)} 
    \left\{
    \sum_{\{i,j\}\subset I_0} d_0(\vv{\Sigma}_i,\vv{\Sigma}_j)
    +
    \sum_{i\in I_0^c}\sum_{j\in I_0} 
    d_0(\vv{\Sigma}_i,\vv{\Sigma}_j)
    \right\} \\
    & \geq \frac{2}{q(q-1)}\sum_{i\in I_0^c}\sum_{j\in I_0} 
    d_0(\vv{\Sigma},\vv{\Sigma}_j)\\
    & \geq
    \frac{2K(q-K)}{q(q-1)} d'_{3,q}. 
\end{split}
\end{equation*} 
Therefore, in (\ref{eq:alterpf_eq}), we have
\begin{equation*}\label{eq:alt_meandrift'}
    \frac{\sqrt{q}m_p}{\lambda_{0,p}} \geq \frac{(q-K)Kd'_{3,q}}{2(q-1)\mu_2\sqrt{q\bar{c}_{q,12}}} 
\end{equation*}
Hence, if $Kd'_{3,q}/\sqrt{q} \to \infty$,
then $\sqrt{q}m_p\to \infty$ and thus the power will eventually tend to $1$.
    \hfill$\Box$
\end{example}
}

{
Finally, 
we consider the problem of ``finding a needle in a haystack'', where there is a single outlier $(K=1)$ among these $q$ populations. The relevant results will later be used in  Section \ref{subsec:3.4}.
\begin{example}[Finding a needle in a haystack]\label{example:4'}
Consider the following alternative:
\begin{equation}\label{eq:simu0_eq}
\begin{split}
    H^{(\beta)}_1:~&\vv{\Sigma}_i =
    \begin{cases} 
        \vv{\Sigma}, & \quad i=1,\ldots,q-1; \\
        \vv{\Sigma}+\sqrt{\beta}\vv{\Lambda}, & \quad i = q, 
    \end{cases}  
\end{split}
\end{equation}
where $\vv{\Lambda}\neq \vv{\Sigma}$. Here,  $\beta$, a non-negative variable, is introduced to control the degree of deviation of $H_1^{(\beta)}$ from $H_0$ in a way that $H_1^{(0)} = H_0$ when $\beta =0$, and $H_1^{(\beta)}$ gradually deviates from $H_0$ and as $\beta$ grows. 

For $\beta>0$, since there is only a single outlier, we have
we have 
\[
    d'_{3,q} = d_0(\vv{\Sigma},\vv{\Sigma}+\sqrt{\beta}\vv{\Lambda}) =\beta {\rm tr}(\vv{\Lambda}^2).
\]
In addition, the mean drift
\begin{equation*}
\sqrt{q}m_p = \frac{2\sqrt{q}}{q(q-1)}\sum_{i=1}^{q-1} d_0(\vv{\Sigma}_q,\vv{\Sigma}) = \frac{2}{\sqrt{q}}\beta {\rm tr}(\vv{\Lambda}^2),
\end{equation*}
Meanwhile, we also have $\lambda_{0,p}^2 = 16\mu_2^2 \bar{c}_{q,12}$ and 
\begin{align*}
    \frac{\sqrt{q}m_p}{\lambda_{0,p}} &= \frac{\beta {\rm tr}(\vv{\Lambda}^2)}{2\mu_2\sqrt{q\bar{c}_{q,12}}}. 
\end{align*}
Combining with (\ref{eq:alterpf_eq}), the power function now becomes
\begin{equation}\label{eq:sim1'}
    \P_{H_1}\left(\frac{\sqrt{q}V_p}{\hat{\lambda}_{p}} > z_\alpha\right)
        =\Phi\left(\beta \frac{{\rm tr}(\vv{\Lambda}^2)}{2\mu_2\sqrt{q\bar{c}_{12,q}}} - z_\alpha\right) + o(1).
\end{equation}
Thus, if $\liminf_{p,q\to\infty} {\rm tr}(\vv{\Lambda}^2)/\sqrt{q}>0$,
value of the power function will gradually tend to one.
    \hfill$\Box$
\end{example}  
}

\subsection{Simulation results}\label{subsec:3.4}
The simulation design follows the settings {in Example \ref{example:4'}}. Similar to Section \ref{ssec:2.5}, we set $p=100$, {the number of populations $q = 50, 100, 150$, and
the sample sizes $\{n_j\}$ randomly picked from $\{50,\ldots, 150\}$}.  We consider the standard normal noise $z_{ij}\sim\mathcal{N}(\vv{0},\vv{I}_p)$ and the Gamma noise $z_{ij}\sim {\sf Gamma}(4,2)-2$. 

{
In Example \ref{example:4'}, we set the last population $\vv{
\Sigma}_q
$ as the unique outlier ($K = 1$),
of which  $\vv{\Sigma}_q=\vv{\Sigma}_0+\sqrt{\beta}\vv{\Lambda}_0$,
} and the other populations are all identical to the matrix $\vv{\Sigma}_0$. Here, {
the parameter $\beta$ varies from $0$ to a selected positive number $\beta_{\max}$
according to different scenarios. In addition, for two  given distinct non-negative definite matrices $\vv{\Lambda}$ and $\vv{\Sigma}$, we define 
\[
\vv{\Sigma}_0 = \vv{\Sigma}/\sqrt{p^{-1}{\rm tr}(\vv{\Sigma}^2)}~~~\hbox{and}~~~
\vv{\Lambda}_0=\vv{\Lambda}/\sqrt{p^{-1}{\rm tr}(\vv{\Lambda}^2)}.
\]
The normalization procedure helps reduce  the theoretic curve in (\ref{eq:sim1'}) to be
\begin{equation}\label{eq:sim1'''}
\P_{H_1}\left(\frac{\sqrt{q}V_p}{\hat{\lambda}_{p}} > z_\alpha\right)
        =\Phi\left(\frac{\beta p}{2\sqrt{q\bar{c}_{12,q}}} - z_\alpha\right) + o(1),
\end{equation}
where $\bar{c}_{12,q}= (q-1)^{-1}\sum_{i=1}^{q-1} c_{ip}^2$.}

The following three scenarios of $\vv{\Sigma}$ and $\vv{\Lambda}$ are considered: 

(a) $\vv{\Sigma}=\vv{I}_p$, and $\vv{\Lambda}=\vv{U}{\rm diag}\{\sqrt{\beta},\cdots,\sqrt{\beta},0,\ldots,0\}\vv{U}^\T$, where the first $p/2$ entries equal to ${\beta}^{1/2}$ and $U$ is a randomly chosen orthogonal matrix; 

(b) $\vv{\Sigma}=\vv{U}_1 \vv{D}_1 \vv{U}_1^\T$ and $\vv{\Lambda}=\vv{U}_2 \vv{D}_2 \vv{U}_2^\T$, where $\vv{U}_1$ and $\vv{U}_2$ are two independently and randomly chosen orthogonal matrices, $\vv{D}_1$ and $\vv{D}_2$ are diagonal matrices with their entries randomly picked from the interval $(0.1,10.1)$.

\begin{figure}
\includegraphics[width = 43 mm, height = 43 mm]{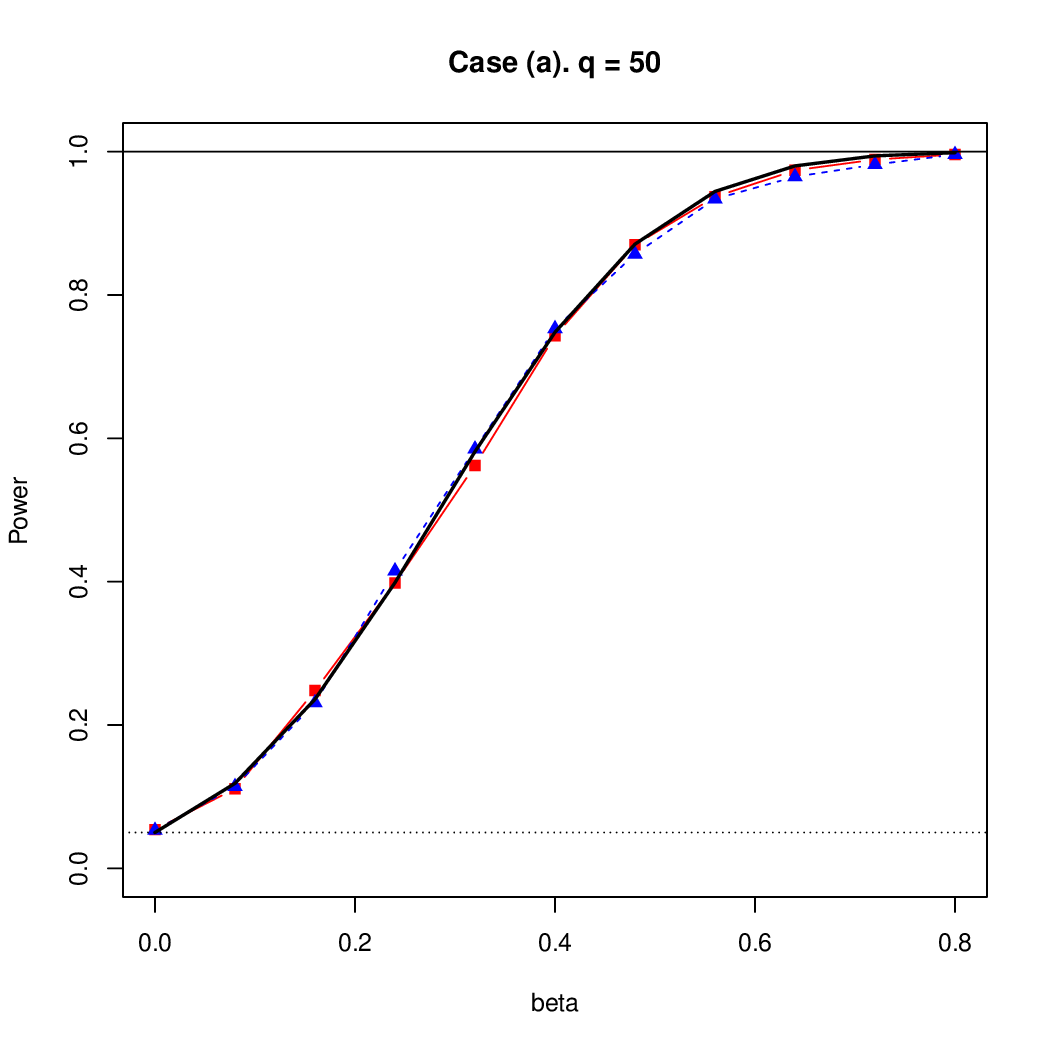}
\includegraphics[width = 43 mm, height = 43 mm]{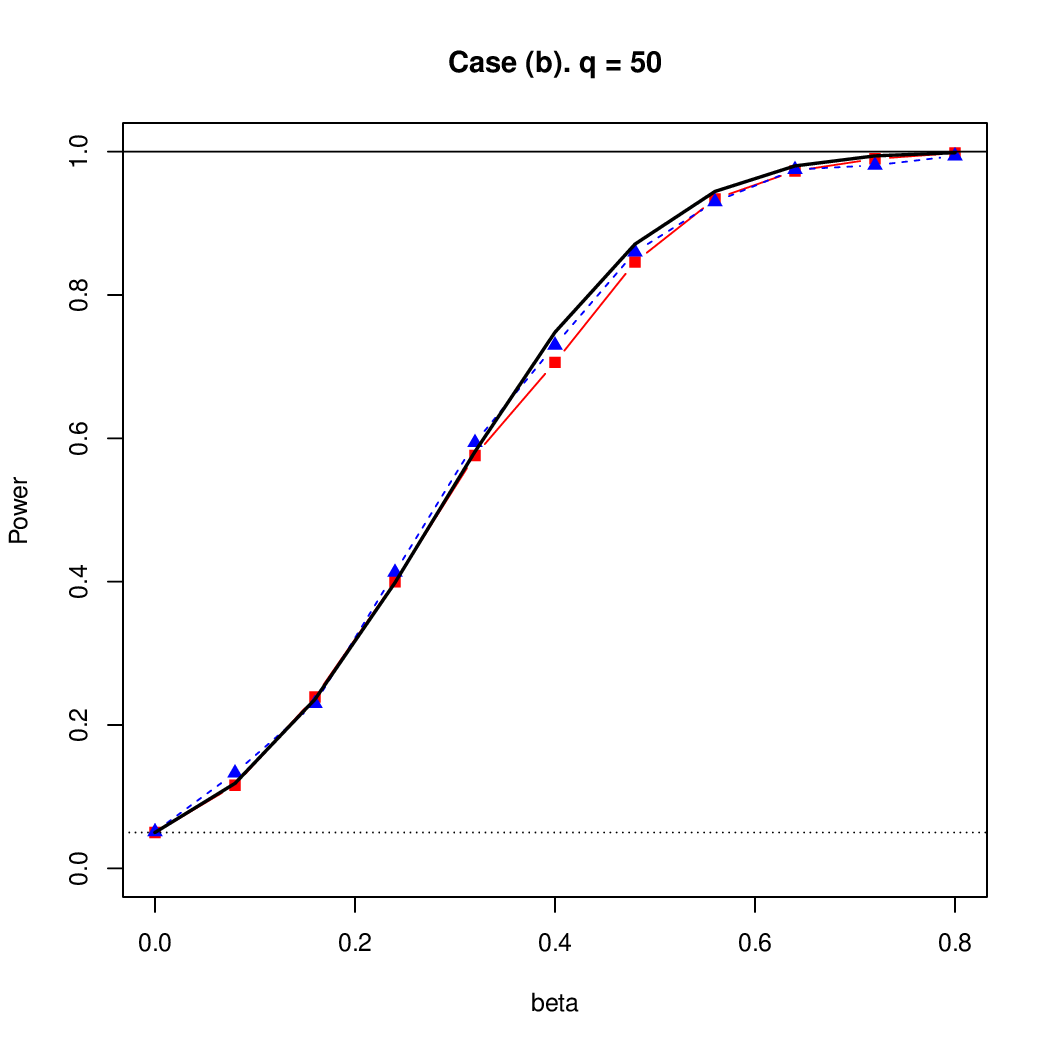}
\includegraphics[width = 43 mm, height = 43 mm]{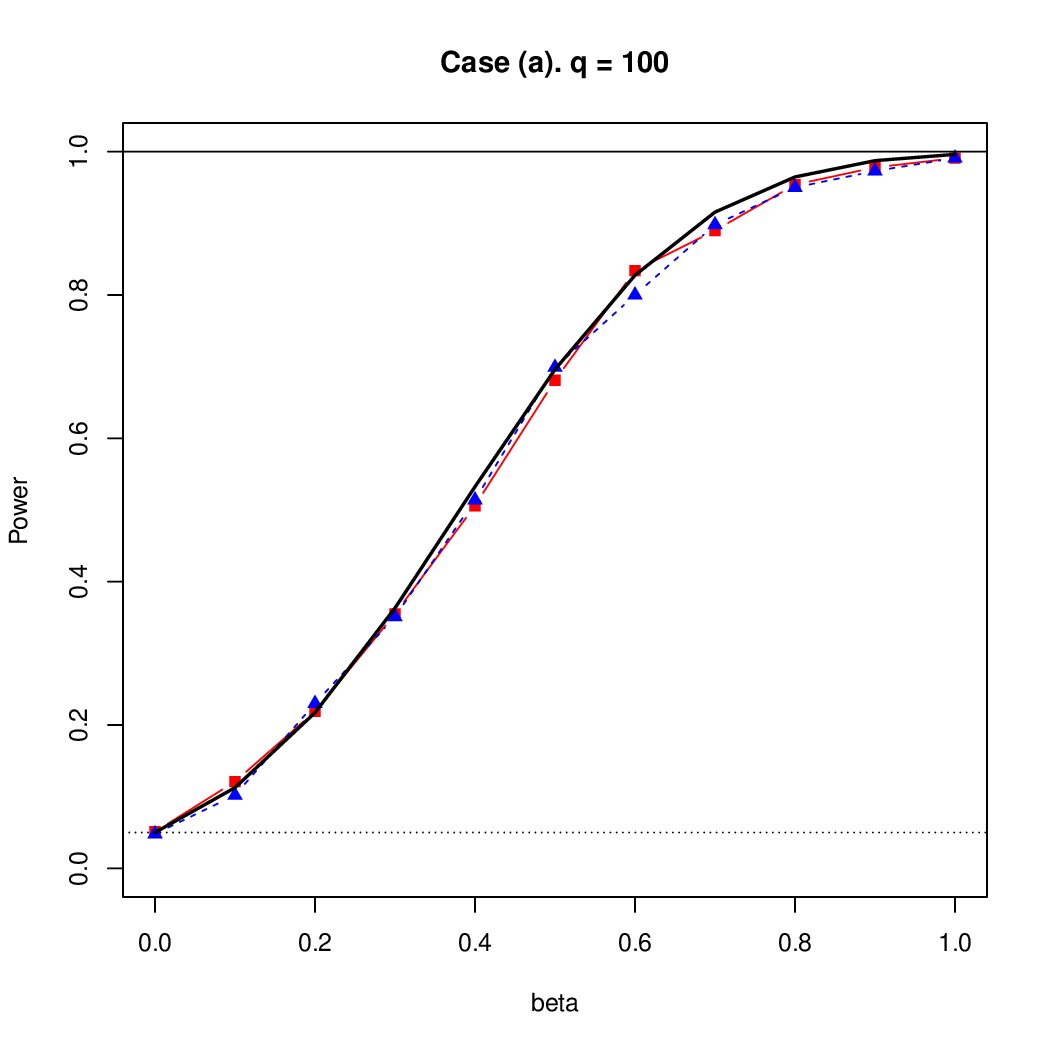}
\includegraphics[width = 43 mm, height = 43 mm]{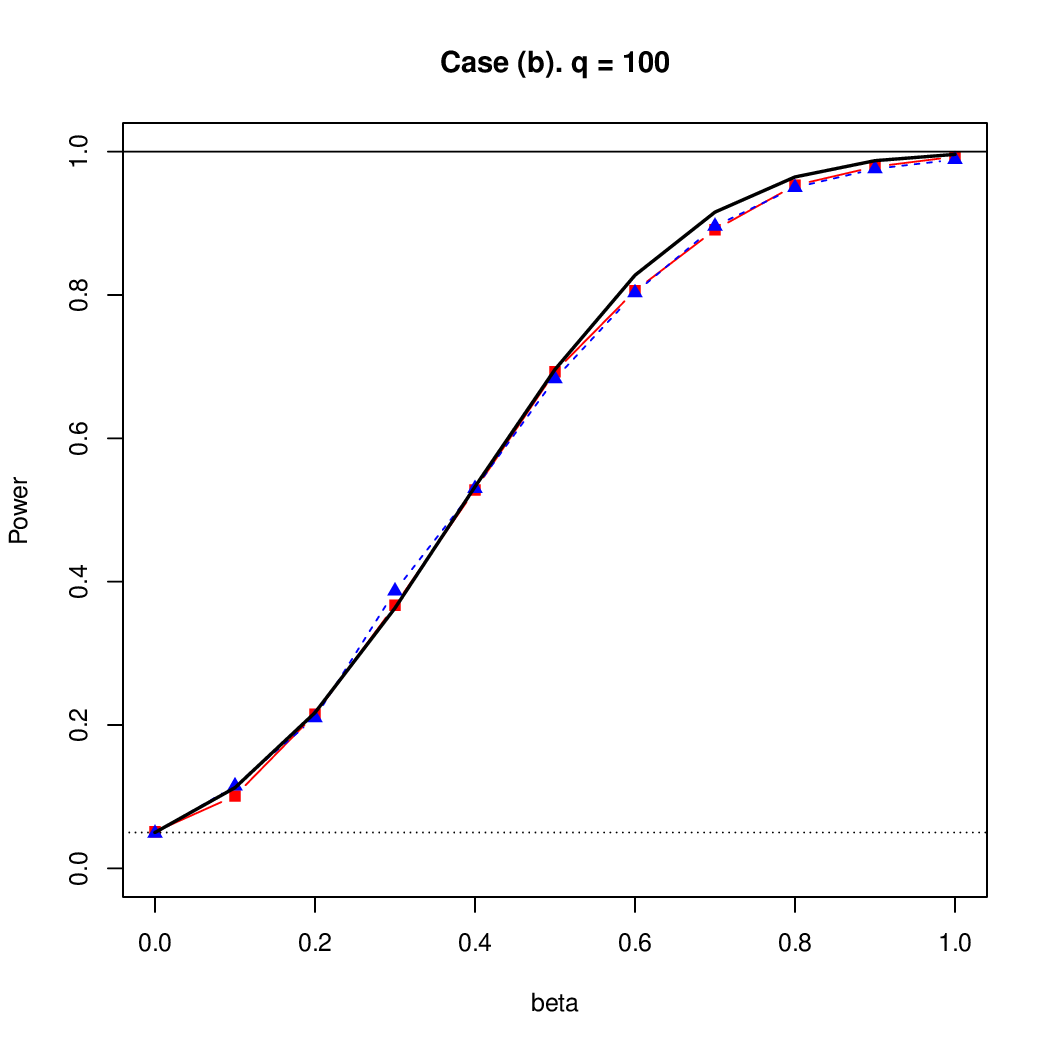}
\includegraphics[width = 43 mm, height = 43 mm]{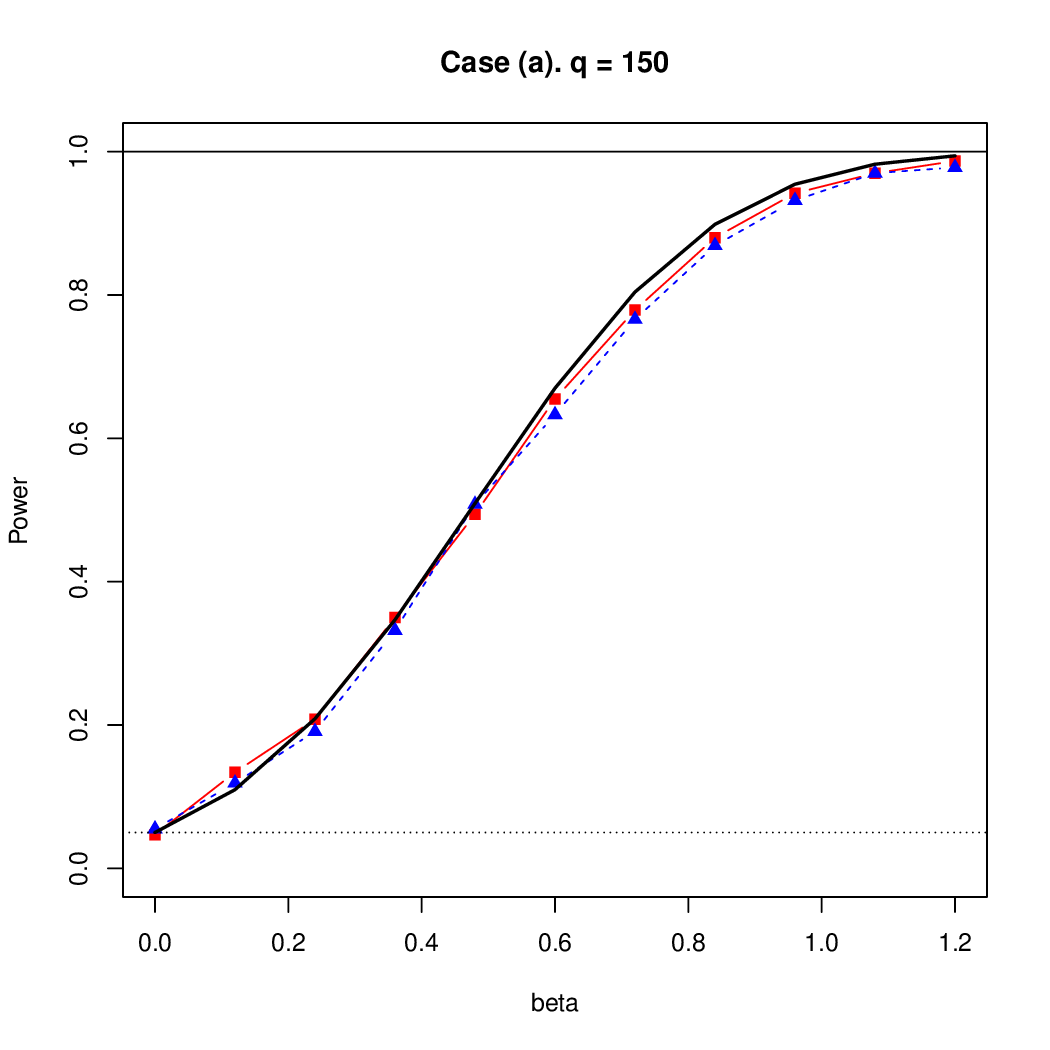}
\includegraphics[width = 43 mm, height = 43 mm]{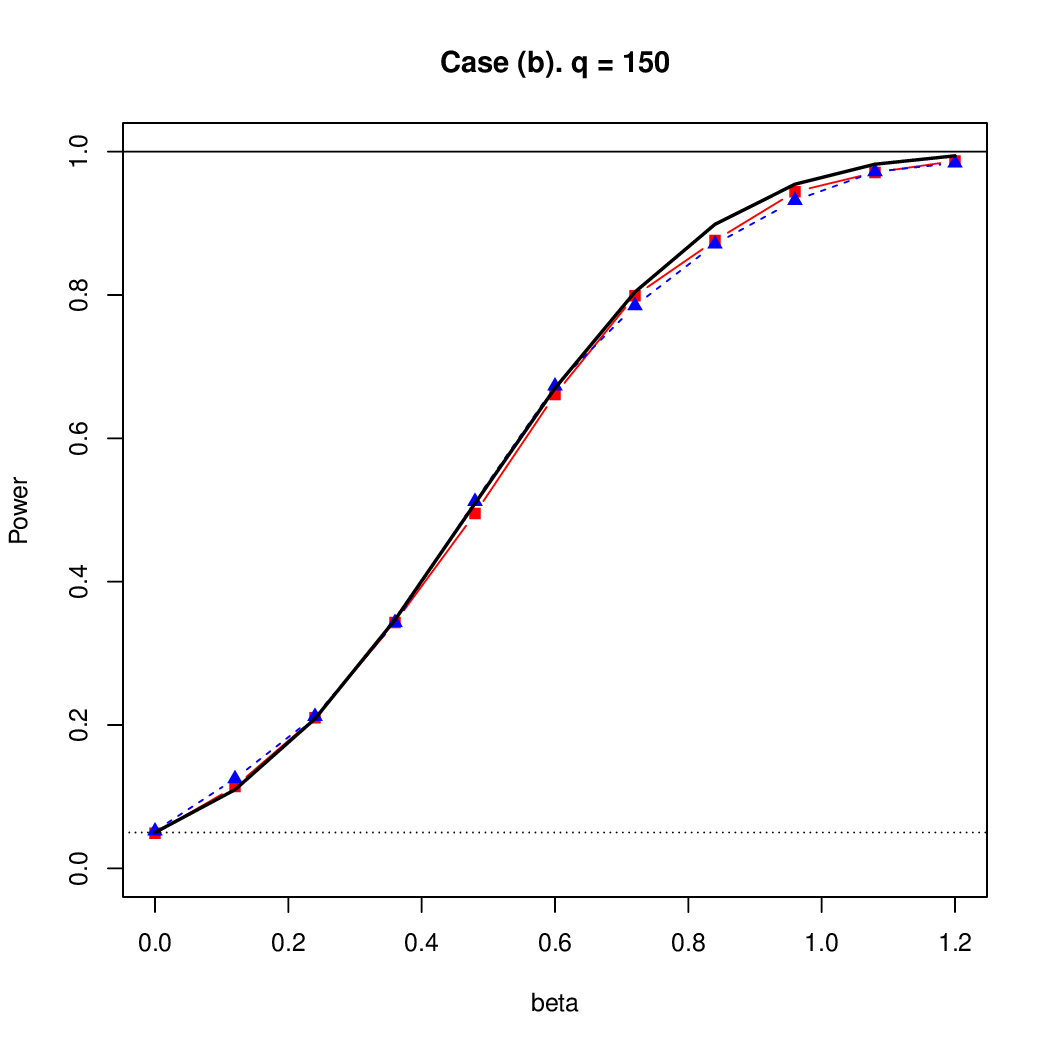}
\caption{Empirical power plots for the many-sample equality test. Two design cases with two different noises, normal noise (red square) and Gamma noise (blue triangle), under three different settings: $q=50, 100$, and $150$. The black dash curves stand for the theoretical power functions given in (\ref{eq:sim1'''}).}
\label{fig:2}
\end{figure}

{For each $\beta$ from $0$ to $\beta_{\max}$ with step size $0.1$}, we repeat simulations $1000$ times and compute the empirical size and power.  The empirical sizes are shown in Table \ref{tab:3}, and the curves of empirical power are displayed in Figure \ref{fig:2}. Our proposed statistic well controls the empirical size under $H_0$, and effectively rejects $H_0$ when the outlier covariance matrix deviates from the majority as $\beta$ increases. {Meanwhile, the empirical power curves are very close to the theoretical ones, which agrees with our previous theoretical results.}

\begin{table}
\caption{Empirical sizes of the many-sample equality test.} 
\label{tab:3}
\begin{tabular}{@{}ccccc@{}}
\hline
&\multicolumn{2}{c}{Case (a)}&
\multicolumn{2}{c}{Case (b)}
 \\ [3pt]
\hline
&Normal & Gamma & 
Normal & Gamma
 \\
\hline
$q = 50$& $0.053$ & $0.045$ & 
$0.054$ & $0.060$
 \\
$q=100$ &  $0.049$& $0.054$ & $0.053$ & 
$0.050$
\\
$q=150$ &  $0.047$& $0.050$ & $0.053$ & 
$0.049$
\\
\hline
\end{tabular}
\end{table}

\subsection{The 1000 Genomes Project (phase 3)}\label{sec:5.2}

The equality assumption on covariance matrices is often assumed in performing analysis of variance or covariance (ANOVA) in genetic data analysis; see for example, \cite{PP2006} and \cite{Price2006}. We here aim to check the validity of such an assumption by investigating a gene dataset from the 1000 Genomes Project phase 3 in \cite{1000gene}. The dataset contains phased genotypes for $112515$ genes in the $22$th chromosome of $2504$ individuals from $q=26$ populations. The list of abbreviations, names, and sample sizes of the  $26$ populations is provided below in Table \ref{tab:4}.

\begin{table}
\caption{Information about $26$ populations in The 1000 Genomes Project phase 3.}
\label{tab:4}
\begin{tabular}{clr}
 \\
{\bf Abbreviation} &  {\bf Name} & {\bf Size} \\
\hline
ACB & (African Caribbean in Barbados) & 96\\
ASW & (African ancestry in Southwest USA) & 61\\
BEB & (Bengali in Bangladesh) & 86\\
CDX & (Chinese Dai in Xishuangbanna,China) & 93\\
\multirow{2}{*}{CEU}& (Utah residents with Northern and &\multirow{2}{*}{99}\\
& ~Western European ancestry) & \\
CHB & (Han Chinese in Beijing, China) & 103\\
CHS & (Southern Han Chinese, China) & 105\\
CLM & (Colombian in Medellin, Colombia) & 94\\
ESN & (Esan in Nigeria) & 99\\
FIN & (Finnish in Finland) & 99\\
GBR & (British in England and Scotland) & 91\\
GIH & (Gujarati Indian in Houston, Texas) & 103\\
GWD & (Gambian in Western Division,The Gambia) & 113\\
IBS & (Iberian populations in Spain) & 107\\
ITU & (Indian Telugu in the UK) & 102\\
JPT & (Japanese in Tokyo,Japan) & 104\\
KHV & (Kinh in Ho Chi Minh City,Vietnam) & 99\\
LWK & (Luhya in Webuye,Kenya) & 99\\
MSL & (Mende in Sierra Leone) & 85\\
MXL & (Mexican ancestry in Los Angeles,California) & 64\\
PEL & (Peruvian in Lima,Peru) & 85\\
PJL & (Punjabi in Lahore,Pakistan) & 96\\
PUR & (Puerto Rican in Puerto Rico) & 104\\
STU & (Sri Lankan Tamil in the UK) & 102\\
TSI & (Toscani in Italy) & 107\\
YRI & (Yoruba in Ibadan,Nigeria) & 108\\
\hline
total& &2504
\end{tabular}
\end{table}

Assume that $26$ populations are mutually independent and for each population, the data are i.i.d. samples from the following model:
\[
\vv{x}_{i,j}=\vv{\mu}_i+\vv{\Sigma}_i^{1/2}\vv{z}_{i,j},~~~~~j=1,\ldots,n_i;~i=1,\ldots,q,
\]   
where $\vv{x}_{i,j}$ stands for the $j$th sample for the $i$th population, $\vv{\mu}_i$ and $\vv{\Sigma}_i$ stand for the population-specific mean and covariance matrix, and $\vv{z}_{i,1},\ldots,\vv{z}_{i,j}$ are i.i.d. random vectors having i.i.d. entries with 
{zero mean, unit variance} and finite eighth moments. Our focus is on checking the validity of  the equality hypothesis:  ``$\vv{\Sigma}_1=\cdots=\vv{\Sigma}_q$''.

As the ultra-high dimension of the observations $112515$ is much larger than the sample sizes, a sub-sampling method is used to investigate the equality of reduced population covariance matrices with size $p=100$.  In other words,  we repeat the experiment $10000$ times. In each experiment, we randomly draw $p=100$ genes from the total of $112515$ ones, centralize the data by the sample mean for each population, perform our testing of equality for sub-covariance matrices, and compute the values of our testing statistic. This reduction to sub-matrices indeed makes sense since the global equality hypothesis  ``$\vv{\Sigma}_1=\cdots=\vv{\Sigma}_q$''implies the equality of reduced covariance matrices from sub-sampled genes, and that all the test statistics is asymptotically standard normal.

It turns out that our testing statistic ${q}^{1/2}V_p/\hat{\lambda}_p$ in these $10000$ experiments ranges from $58.679$ to $90.899$ with mean $74.730$, which are far greater than the $95\%$ standard normal quantile and thus lead to a full and strong rejection of the global equality hypothesis.  Consequently, we conclude that it is indeed not reasonable to assume the equality of the $26$ population covariance matrices in this 1000 Genomes Project phase 3 dataset.

Furthermore, we investigate the contributions of different pairs of populations to the extremely large values of our statistics. Let us take a closer look at our testing statistic 
\[
\frac{{q}^{1/2}V_p}{\hat{\lambda}_p}=\frac{1}{q(q-1)}\sum_{i\neq j} \frac{{q}^{1/2}g(\vv{X}_{ip},\vv{X}_{jp})}{\hat{\lambda}_p}=:\frac{1}{q(q-1)}\sum_{i\neq j}G_{ij},
\]
where $g(\vv{X}_{ip},\vv{X}_{jp})$, $V_p$ and $\hat{\lambda}_p$ are given in (\ref{eq:ub_mp}), (\ref{eq:3.20}) and (\ref{eq:3.9}), respectively. Here, $G_{ij}$ represents the contribution of the {$i$th and $j$th} populations.  Indeed, according to discussions in Section \ref{subsec:3.1} and \ref{subsec:3.2}, $g(\vv{X}_{ip},\vv{X}_{jp})$ is an unbiased and consistent estimator for 
 $d_0(\vv{\Sigma}_i,\vv{\Sigma}_j)$ and $\hat{\lambda}_p$ is a consistent estimator for
\[
4\bigg[\frac{1}{q}\sum_{j=1}^q c_{jp}^2\left(\frac{1}{p}{\rm tr}(\vv{\Sigma}_j^2)\right)^2\bigg]^{1/2},
\]
in which the latter can be viewed as a global normalizer. Denote $\bar{G}_{ij}$ be the averaging value of $G_{ij}$ over the above-mentioned $10000$ repeated experiments: this statistic is, in fact, the two-sample statistic if we were, testing the single hypothesis ``$\vv{\Sigma}_i=\vv{\Sigma}_j $'' and  
a large value means that the two covariance matrices are likely distinct.  

\begin{figure}
\includegraphics[width = 80 mm, height = 70 mm]{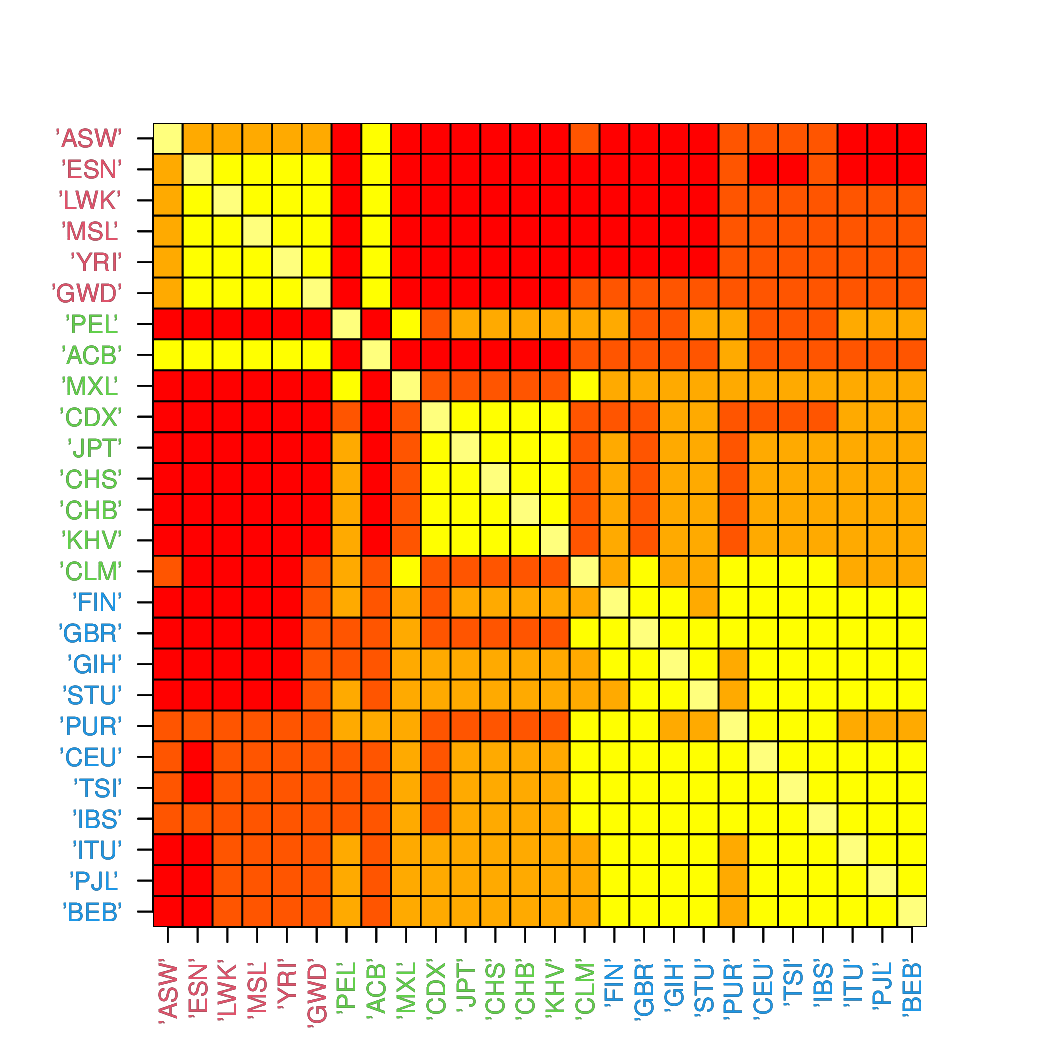}
\caption{A heat-map for pairwise test statistics  $(\bar{G}_{ij})$  with 4 classes. 
  Rows are sorted in descending order according to their row averages. 
  Diagonal values are set to 0.}
\label{fig:4}
\end{figure}

We find that the values of $\bar{G}_{ij}$'s vary from {$24.415$ to $51.521$ with a mean of $36.278$}. Moreover, we also classify the values of $\bar{G}_{ij}$'s into four classes $\{1,2,3,4\}$ by comparing them with the sample quartiles {$Q_{75\%}=40.306,Q_{50\%}=34.475$ and $Q_{25\%}=28.506$}. Figure \ref{fig:4} displays a heat map with these classes for all the population pairs $(i,j)$ 
(the diagonal entries are set to zero) where the populations are sorted according to their row averages  ${q^{-1}}\sum_{j=1}^q \bar{G}_{ij}$ in a decreasing order. 

According to Figure \ref{fig:4}, we can divide $26$ populations into three sub-groups based on their contributions to the extremely large value of our testing statistic. The main contributions come from {the following six populations: ``ASW'', ``ESN'', ``LWK'',``MSL'',}{``YRI'' and ``GWD''}. It can be seen that these six deviate significantly from others since the corresponding statistic values $\bar{G}_{ij}$  almost always lie beyond the third quartile  $Q_{75\%}$. Besides, {there are also nine populations with medium contributions} since their corresponding statistic values mostly lie between the second and third quartiles.  {The other eleven populations tend to take smaller values below the first quartile}, though they are still not small enough to allow us to treat their covariance matrices as equal.

\begin{acks}[Acknowledgments]
The authors wish to express their appreciation to the anonymous reviewers, the Associate Editor, and the Editor for their insightful feedback and constructive suggestions, which have greatly contributed to the improvement of this paper.
\end{acks}

\begin{funding}
Chen Wang was partially supported by Hong Kong RGC
General Research Fund 17301021 and National Natural Science Foundation of China Grant
72033002. Jianfeng Yao was partially supported by NSFC RFIS Grant No. 12350710179.
\end{funding}

\begin{supplement}
\stitle{Supplement to ``Many-sample tests for the equality and the
proportionality hypotheses between large covariance matrices"}
\sdescription{In this supplement, we provide all technical proofs.}
\end{supplement}

\bibliographystyle{imsart-nameyear} 
\bibliography{mybibfile,cov-test}

\newpage
\setcounter{page}{1}
\setcounter{section}{0}

\begin{supplement}
\stitle{Supplement to ``Many-sample tests for the equality and the
proportionality hypotheses between large covariance matrices"}
\sdescription{
\begin{center}
Tianxing Mei, Chen Wang and Jianfeng Yao\\[2ex]
\end{center}
      }
\end{supplement}

\appendix
\section{Auxiliary Lemmas}\label{supp:aux}
{
Let $\vv{x}_1,\ldots,\vv{x}_n$ be an i.i.d. sample from a $p$-dimensional population $\vv{x} =\vv{\Sigma}^{1/2} \vv{z}$,  where $\vv{\Sigma}$ is a semi-definite matrix and $\vv{z}=(z_{j})$ has i.i.d. entries with zero mean, unit variance, finite fourth moment $\nu_4$ and finite eighth moment.

For any semi-definite matrix $\vv{\Lambda}$, and $p\times p$ symmetric matrices $\vv{A}_1$ and $\vv{A}_2$, we denote
\begin{equation}\label{eq:A.0}
\begin{split}
    \langle\vv{A}_1,\vv{A}_2\rangle_{\vv{\Lambda}} &= \frac{2}{p}{\rm tr}(\vv{\Lambda}\vv{A}_1\vv{\Lambda}\vv{A}_2)\\
    &+\frac{(\nu_4-3)}{p}{\rm tr}(\mathcal{D}(\vv{\Lambda}^{1/2}\vv{A}_1\vv{\Lambda}^{1/2})\mathcal{D}(\vv{\Lambda}^{1/2}\vv{A}_2\vv{\Lambda}^{1/2}))
\end{split}
\end{equation}
where $\mathcal{D}(\vv{A})$ is a diagonal matrix made by diagonal entries of $\vv{A}$.
We directly see that when $\vv{\Lambda}$ is fixed, $\langle\cdot,\cdot\rangle_{\vv{\Lambda}}$ is a non-negative definite bi-linear function on the product space of symmetric matrices. 
Moreover, when $\vv{\Lambda}$, $\vv{A}_1$ and $\vv{A}_2$ are bounded uniformly in operator norm $\|\cdot\|$ with respect to $p$, so is $\langle\vv{A}_1,\vv{A}_2\rangle_{\vv{\Lambda}}$.

Let $\vv{S}_n$ be the 
sample covariance matrix, that is, $\vv{S}_n=\frac{1}{n}\sum_{j=1}^n \vv{x}_i\vv{x}_i^\T$. 
Here, we focus on moments of the following statistics related to $\vv{S}_n$:
\[
\frac{1}{p}{\rm tr}(\vv{S}_n^2)~~~\hbox{and}~~~~\left(\frac{1}{p}{\rm tr}(\vv{S}_n)\right)\left(\frac{1}{p}{\rm tr}(\vv{S}_n\vv{A})\right),
\]
where $\vv{A}$ is a $p\times p$ symmetric deterministic matrix.
These two quantities play a crucial role in deriving the asymptotic normality of our proposed $U$-statistic in Sections \ref{sec:2} and \ref{sec:3}.

The following lemmas provide the expectations, covariances, and bounds of the fourth moment about these statistics. The proofs of the results are based on the formulas of quadratic forms of random vectors discussed in \cite{ullah2004finite} and \cite{Bao20101193} with some lengthy but tedious calculations. We provide detailed proofs for some of the crucial formulas in the last Appendix \ref{proof:aux}. 

The first result displays the expectation formulas for the two statistics.
\begin{lemma}\label{lem:expectation}
Suppose that $\{\vv{x}_j\}$ satisfies the structure given at the beginning of this section. Then, it holds that
\begin{equation}\label{eq:A.1}
\begin{split}
    \E\left[\frac{1}{p}{\rm tr}(\vv{S}_n^2)\right] & = \frac{1}{p}{\rm tr}(\vv{\Sigma}^2) + \frac{p}{n}\left(\frac{1}{p}{\rm tr}(\vv{\Sigma})\right)^2 \\
    & + \frac{1}{n}\left\{\frac{1}{p}{\rm tr}(\vv{\Sigma}^2)+\frac{(\nu_4-3)}{p}{\rm tr}(\mathcal{D}(\vv{\Sigma})^2)\right\}.
\end{split}
\end{equation}
Meanwhile, for any symmetric matrix $\vv{A}$ , we also have
\begin{equation}\label{eq:A.2}
\begin{split}
\mathbb{E}\left[\left(\frac{1}{p}{\rm tr}(\vv{S}_{n})\right)\left(\frac{1}{p}{\rm tr}(\vv{S}_{n}\vv{A})\right)\right]
&=\left(\frac{1}{p}{\rm tr}(\vv{\Sigma})\right)\left(\frac{1}{p}{\rm tr}(\vv{\Sigma}\vv{A})\right)\\
&+\frac{1}{pn}\langle\vv{I}_p,\vv{A}\rangle_{\vv{\Sigma}}.
\end{split}
\end{equation}

\end{lemma}
The next lemma provides the covariance formulas.
\begin{lemma}\label{lem:jointcov}
Suppose that $\{\vv{x}_j\}$ satisfies the pressumed structure. Then, we have
\begin{equation}\label{eq:A.3}
\begin{split}
{\rm Var}\left({\rm tr}(\vv{S}_{n}^2)\right)&=\frac{1}{n^4}\bigg\{ 4n{\rm tr}(\vv{\Sigma})^2\mu(\vv{I}_p,\vv{I}_p|\vv{\Sigma}) +8n^2 {\rm tr}(\vv{\Sigma})\langle\vv{\Sigma},\vv{I}_p\rangle_{\vv{\Sigma}}\\
&~~~~~~~~~+4n^2{\rm tr}(\vv{\Sigma}^2)^2+4n^3\langle\vv{\Sigma},\vv{\Sigma}\rangle_{\vv{\Sigma}}\bigg\}(1+ o(1)).
\end{split}
\end{equation}
In addition, for $p\times p$ symmetric matrices $\vv{A}$, $\vv{A}_1$ and $\vv{A}_2$, we also have 
\begin{equation}\label{eq:A.4}
\begin{split}
&{\rm Cov}\left(\frac{1}{p}{\rm tr}(\vv{S}_{n}){\rm tr}(\vv{S}_n\vv{A}_1),\frac{1}{p}{\rm tr}(\vv{S}_{n}){\rm tr}(\vv{S}_n\vv{A}_2)\right)\\
=&\frac{1}{pn}
\left\{\begin{aligned}
&{\rm tr}(\vv{\Sigma}){\rm tr}(\vv{\Sigma})\langle\vv{A}_1,\vv{A}_2\rangle_{\vv{\Sigma}}\\
+&{\rm tr}(\vv{\Sigma}){\rm tr}(\vv{\Sigma}\vv{A}_1)\langle\vv{I}_p,\vv{A}_2\rangle_{\vv{\Sigma}}\\
+&{\rm tr}(\vv{\Sigma}){\rm tr}(\vv{\Sigma}\vv{A}_2)\langle\vv{I}_p,\vv{A}_1\rangle_{\vv{\Sigma}}\\
+&{\rm tr}(\vv{\Sigma}\vv{A}_1){\rm tr}(\vv{\Sigma}\vv{A}_2)\langle\vv{I}_p,\vv{I}_p\rangle_{\vv{\Sigma}}
\end{aligned}\right\}(1+o(1))
\end{split}
\end{equation}
and
\begin{equation}\label{eq:A.4.1}
\begin{split}
&{\rm Cov}\left({\rm tr}(\vv{S}_{n}^2),\frac{1}{p}{\rm tr}(\vv{S}_{n}){\rm tr}(\vv{S}_{n}\vv{A})\right)\\
=&\frac{2}{n^4}\left\{\begin{aligned}
&n^2\cdot{\rm tr}(\vv{\Sigma}){\rm tr}(\vv{\Sigma})\langle\vv{A},\vv{I}_p\rangle_{\vv{\Sigma}}\\
+&n^2\cdot{\rm tr}(\vv{\Sigma}\vv{A}){\rm tr}(\vv{\Sigma})\langle\vv{I}_p,\vv{I}_p\rangle_{\vv{\Sigma}}\\
+&n^3\cdot {\rm tr}(\vv{\Sigma}\vv{A})\langle\vv{I}_p,\vv{\Sigma}\rangle_{\vv{\Sigma}}\\
+&n^3\cdot {\rm tr}(\vv{\Sigma})\langle\vv{A},\vv{\Sigma}\rangle_{\vv{\Sigma}}
\end{aligned}\right\}(1+o(1)). 
\end{split}
\end{equation}
Consequently, it also holds that
\begin{equation}\label{eq:A.5}
\begin{split}
{\rm Var}\left({\rm tr}(\vv{S}_{n}^2)-\frac{1}{n}{\rm tr}(\vv{S}_{n})^2\right)=4c_{p}^2\left(\frac{1}{p}{\rm tr}(\vv{\Sigma}^2)\right)^2+4c_{p}\langle\vv{\Sigma},\vv{\Sigma}\rangle_{\vv{\Sigma}}
\end{split}
\end{equation}
\begin{equation}\label{eq:A.6}
\begin{split}
&{\rm Cov}\left({\rm tr}(\vv{S}_{n}^2)-c_{p}\frac{1}{p}{\rm tr}(\vv{S}_{n})^2,\frac{1}{p}{\rm tr}(\vv{S}_{n}){\rm tr}(\vv{S}_{n}\vv{A})\right)\\
&=\frac{2}{n}\bigg\{{\rm tr}(\vv{\Sigma}\vv{A})\langle\vv{\Sigma},\vv{I}_p\rangle_{\vv{\Sigma}}+{\rm tr}(\vv{\Sigma})\langle\vv{\Sigma},\vv{A}\rangle_{\vv{\Sigma}}\bigg\}(1+o(1))\\
\end{split}
\end{equation}
where $c_p = p/n$.
\end{lemma}
The next lemma shows the covariances of two statistics considered in this section with a statistic in the form of ${\rm tr}(\vv{S}_n\vv{A})$.
\begin{lemma}\label{lem:single_cov}
Let $\{\vv{x}_j\}$ be the i.i.d. sample from the population satisfying the structure given before. Then, for any $p\times p$ symmetric matrix $\vv{A}$, we have
    \begin{equation}\label{eq:lemA3.1}
    \begin{split}
    &{\rm Cov}\left({\rm tr}(\vv{S}_n^2), {\rm tr}(\vv{S}_n\vv{A}) \right) \\
    =& 2c_p^2\left(\frac{1}{p}{\rm tr}(\vv{\Sigma})\right)\langle\vv{I}_p,\vv{A}\rangle_{\vv{\Sigma}} + 2c_p \langle\vv{\Sigma},\vv{A}\rangle_{\vv{\Sigma}}+o(1)
\end{split}
\end{equation}
and
\begin{equation}\label{eq:lemA3.2}
\begin{split}
    {\rm Cov}\left(\frac{1}{n}{\rm tr}(\vv{S}_n)^2, {\rm tr}(\vv{S}_n\vv{A})\right) 
    =2c_p^2\left(\frac{1}{p}{\rm tr}(\vv{\Sigma})\right)\langle\vv{A},\vv{I}_p\rangle_{\vv{\Sigma}} + o(1).
\end{split}
\end{equation}
Consequently, we have
\begin{equation}
    \begin{split}
        {\rm Cov}\left({\rm tr}(\vv{S}_n^2) - \frac{1}{n}{\rm tr}(\vv{S}_n)^2, {\rm tr}(\vv{S}_n\vv{A})\right) = 2c_p\langle\vv{\Sigma},\vv{A}\rangle_{\vv{\Sigma}} + o(1).
    \end{split}
\end{equation}
\end{lemma}

The final results conclude the finiteness of their fourth moments.
\begin{lemma}\label{lem:fourthmom}
    For $\{\vv{x}_j\}$ given at the beginning of the section, there exists a positive constant $C$ dependent  only on $\|\vv{\Sigma}\|$ and the ratio $c_p=p/n$ such that
    \[
    \E\left[{\rm tr}(\vv{S}_{n}^2)  - \E[{\rm tr}(\vv{S}_{n}^2)]\right]^4 \leq C.
    \] 
    Meanwhile, for any symmetric matrices $\vv{A}_1$ and $\vv{A}_2$ that are bounded in operator norm, there is another positive constant $C'$ relying only on $c_p$ and the operator norms of $\vv{\Sigma}$, $\vv{A}_1$ and $\vv{A}_2$ such that
    \[
    \E\left[\frac{1}{p}{\rm tr}(\vv{S}_{n}){\rm tr}(\vv{S}_n\vv{A})  - \E\left(\frac{1}{p}{\rm tr}(\vv{S}_{n}){\rm tr}(\vv{S}_n\vv{A})\right)\right]^4 \leq C'.
    \]
\end{lemma}
}

\section{Proofs  in Section \ref{subsec:2.1}}
{
\subsection{Proof of Lemma \ref{lem:mom}}\label{sec:lem_mom}

Recall the definitions of $\hat{\nu}_{i,12}$, $\hat{\nu}_{i,2}$ and $\hat{\nu}_{i,4}$ in Page \pageref{eq:sample_nu}. It is easy see from Lemma \ref{lem:expectation} that
\begin{align*}
    \mathbb{E}[\hat{\nu}_{i,2}(\vv{A})] = &\mathbb{E}\left(\frac{1}{p}{\rm tr}(\boldsymbol{S}_{ip}^2\boldsymbol{A})\right) \\
     = &\frac{1}{p}{\rm tr}(\boldsymbol{\Sigma}_i^2\boldsymbol{A}) + c_n \left(\frac{1}{p} {\rm tr}(\boldsymbol{\Sigma}_i)\right) \left(\frac{1}{p} {\rm tr}(\boldsymbol{\Sigma}_i\boldsymbol{A})\right) \\
    + &\frac{1}{n_i}\bigg\{\frac{1}{p}{\rm tr}(\boldsymbol{\Sigma}_i^2\boldsymbol{A})
    + \frac{(\nu_{4,i}-3)}{p}{\rm tr}(\mathcal{D}(\boldsymbol{\Sigma}_i)\mathcal{D}(\boldsymbol{\Sigma}_i^{1/2}\boldsymbol{A}_2\boldsymbol{\Sigma}_i^{1/2}))\bigg\}\\
     = & \mu_{i,2}(\vv{A}) + c_{ip}\mu_{i,12}(\vv{A}) + \frac{1}{n_i} (\mu_{i,2}(\vv{A}) + \mu_{i,4}(\vv{A})),
\end{align*}
and
\begin{align*}
    \mathbb{E}[\hat{\nu}_{i,12}(\vv{A})]=&\E\left\{\left(\frac{1}{p} {\rm tr}(\boldsymbol{S}_{ip})\right) \left(\frac{1}{p} {\rm tr}(\boldsymbol{S}_{ip}\boldsymbol{A})\right)\right\} \\
    =& \left(\frac{1}{p} {\rm tr}(\boldsymbol{\Sigma}_i)\right) \left(\frac{1}{p} {\rm tr}(\boldsymbol{\Sigma}_i\boldsymbol{A})\right) \\
    +& \frac{1}{pn_i}\bigg\{\frac{2}{p}{\rm tr}(\boldsymbol{\Sigma}_i^2\boldsymbol{A})
    + \frac{(\nu_{4,i}-3)}{p}{\rm tr}(\mathcal{D}(\boldsymbol{\Sigma}_i)\mathcal{D}(\boldsymbol{\Sigma}_i^{1/2}\boldsymbol{A}\boldsymbol{\Sigma}_i^{1/2}))\bigg\} \\
    = & \mu_{i,12}(\vv{A}) + \frac{1}{pn_i}(2\mu_{i,2}(\vv{A}) + \mu_{i,4}(\vv{A})).
\end{align*}

Observe that for any $1\leq r\neq s\leq n_i$,
\begin{align*}
    &2 {\rm Cov}(\boldsymbol{x}_{i,r}^\T \boldsymbol{x}_{i,r}^\T, \boldsymbol{x}_{i,r}^\T\vv{A} \boldsymbol{x}_{i,r}^\T)\\
    =&\mathbb{E}\left[(\boldsymbol{x}^\T_{i,r}\boldsymbol{x}^\T_{i,r} - \boldsymbol{x}^\T_{i,s}\boldsymbol{x}^\T_{i,s})(\boldsymbol{x}^\T_{i,r}\vv{A}\boldsymbol{x}^\T_{i,r} - \boldsymbol{x}^\T_{i,s}\vv{A}\boldsymbol{x}^\T_{i,s})\right] \\
    =& 2p\left\{\frac{2}{p}{\rm tr}(\boldsymbol{\Sigma}_i^2\boldsymbol{A})
    + \frac{(\nu_{4,i}-3)}{p}{\rm tr}(\mathcal{D}(\boldsymbol{\Sigma}_i)\mathcal{D}(\boldsymbol{\Sigma}^{1/2}_i\boldsymbol{A}\boldsymbol{\Sigma}^{1/2}_i))\right\}\\
    =& 2p(2\mu_{i,2}(\vv{A}) + \mu_{i,4}(\vv{A})).
\end{align*}
Therefore, the statistic $\hat{\nu}_{i,4}(\vv{A})$ is actually a the $U$-statistic, i.e.,
\begin{equation}\label{eq:U-nu4}
\begin{split}
    \hat{\nu}_{i,4}(\vv{A}) & = \frac{1}{pn_i(n_i-1)}\sum_{r<s} (\boldsymbol{x}^\T_{i,r}\boldsymbol{x}_{i,r} - \boldsymbol{x}^\T_{i,s}\boldsymbol{x}_{i,s})(\boldsymbol{x}^\T_{i,r}\vv{A}\boldsymbol{x}_{i,r} - \boldsymbol{x}^\T_{i,s}\vv{A}\boldsymbol{x}_{i,s}) \\
    & = \frac{1}{pn_i(n_i-1)}\sum_{r<s} \left\{\begin{aligned}
        &(\boldsymbol{x}^\T_{i,r}\boldsymbol{x}_{i,r})(\boldsymbol{x}^\T_{i,r}\vv{A}\boldsymbol{x}_{i,r}) + (\boldsymbol{x}^\T_{i,s}\boldsymbol{x}_{i,s})(\boldsymbol{x}^\T_{i,s}\vv{A}\boldsymbol{x}_{i,s}) \\
        -& (\boldsymbol{x}^\T_{i,r}\boldsymbol{x}_{i,r})(\boldsymbol{x}^\T_{i,s}\vv{A}\boldsymbol{x}_{i,s}) - (\boldsymbol{x}^\T_{i,s}\boldsymbol{x}_{i,s})(\boldsymbol{x}^\T_{i,r}\vv{A}\boldsymbol{x}_{i,r})
    \end{aligned}\right\} \\
        & = \frac{1}{p(n_i-1)}\sum_{r=1}^{n_i} \bigg\{
        (\boldsymbol{x}^\T_{i,r}\boldsymbol{x}_{i,r})(\boldsymbol{x}^\T_{i,r}\vv{A}\boldsymbol{x}_{i,r}) 
         -  {\rm tr}(\boldsymbol{S}_{ip}) {\rm tr}(\boldsymbol{S}_{ip}\boldsymbol{A})\bigg\}
\end{split}
\end{equation}
and thus
\[
\mathbb{E}[\hat{\nu}_{i,4}(\vv{A})] = 2\mu_{i,2}(\vv{A}) + \mu_{i,4}(\vv{A}).
\]
The conclusion in Lemma \ref{lem:mom} then follows.

\subsection{Proof of Lemma \ref{lem:ub_est}}\label{sec:lem_ub_est}

Let
\[
\vv{F}_{p,n_i} = \begin{pmatrix}
        1+\frac{1}{n_i} & c_{ip} & \frac{1}{n_i} \\
        \frac{2}{pn_i} & 1 & \frac{1}{pn_i} \\
        2 & 0 & 1
    \end{pmatrix}.
\]
We conclude that the matrix is invertible and next will calculate the explicit form of its inverse. Indeed, The determinant of $\vv{F}_{p,n_i}$ is 
\begin{align*}
    {\rm det}(\vv{F}_{p,n_i}) &= \left(1+\frac{1}{n_i}\right) - \frac{2}{n_i} = 1- \frac{1}{n_i}>0.
\end{align*}
With this in mind, we find the inverse of $\vv{F}_{p,n_i}$:
\begin{equation}\label{eq:inverse}
\begin{split}
    \vv{F}_{p,n_i}^{-1} & = \frac{\vv{F}_{p,n_i}^*}{{\rm det}(\vv{F}_{p,n_i})}  \\
    & = \frac{n_i}{n_i-1} \begin{pmatrix}
       1 & -c_{ip} & -\frac{n_i-1}{n_i^2}\\
       0 & \frac{n_i-1}{n_i} & -\frac{n_i-1}{pn_i^2}\\
       -2 & 2c_{ip} & \frac{(n_i-1)(n_i+2)}{n_i^2}\\
    \end{pmatrix},
    \end{split}
\end{equation}
where $\vv{F}_{p,n_i}^*$ stands for the  adjugate of $\vv{F}_{p,n_i}$.

Hence, by solving the system of equations, we find that
\begin{align*}
    \E\left\{\frac{n_i}{n_i-1} \begin{pmatrix}
       1 & -c_{ip} & -\frac{n_i-1}{n_i^2}\\
       0 & \frac{n_i-1}{n_i} & -\frac{n_i-1}{pn_i^2}\\
       -2 & 2c_{ip} & \frac{(n_i-1)(n_i+2)}{n_i^2}\\
    \end{pmatrix} \begin{pmatrix}
        \hat{\nu}_{i,2}(\vv{A})\\
        \hat{\nu}_{i,12}(\vv{A})\\
        \hat{\nu}_{i,4}(\vv{A})
    \end{pmatrix} \right\}= \begin{pmatrix}
        \mu_{i,2}(\vv{A})\\
        \mu_{i,12}(\vv{A})\\
        \mu_{i,4}(\vv{A})
    \end{pmatrix}.
\end{align*}

We first simplify the expression of $\hat{\nu}_{i,4}(\vv{A})$. Recall the discussion in the last section that $\hat{\nu}_{i,4}(\vv{A})$ is a $U$-statistic and satisfies (\ref{eq:U-nu4}). Observe that 
\begin{align*}
    &\frac{1}{n_i-1}\sum_{r=1}^{n_i} (\boldsymbol{x}^\T_{i,r}\boldsymbol{x}_{i,r})(\boldsymbol{x}^\T_{i,r}\vv{A}\boldsymbol{x}_{i,r})\\
    & = {\rm tr}\left(\vv{A} \left\{\frac{1}{n_i-1} \sum_{r=1}^{n_i}\vv{x}_{i,r}\vv{x}_{i,r}^\T(\vv{x}_{i,r}^T\vv{x}_{i,r})\right\}\right) \\
    & = {\rm tr}\left(\vv{A} \left\{\frac{1}{n_i-1} \vv{X}_{ip}{\rm diag}(\vv{X}^\T_{ip}\vv{X}_{ip}) \vv{X}^\T_{ip}\right\}\right).
\end{align*}
It then follows that
\begin{align*}
    \hat{\nu}_{i,4}(\vv{A}) & = \frac{1}{p}{\rm tr}\left(\vv{A}\left\{\frac{1}{n_i-1}\vv{X}_{ip}{\rm diag}(\vv{X}_{ip}^\T\vv{X}_{ip})\vv{X}_{ip}^\T-\frac{n_i}{n_i-1}{\rm tr}(\vv{S}_{ip})\vv{S}_{ip}\right\}\right).
\end{align*}
For simplicity, we denote
\begin{equation}
    \vv{R}_{i} = \frac{1}{n_i(n_i-1)} \vv{X}_{ip}{\rm diag}(\vv{X}_{ip}^\T\vv{X}_{ip})\vv{X}_{ip}^\T
\end{equation}

The unbiased estimator of $\mu_{i,12}(\vv{A})$
now becomes
\begin{align*}
    \hat{\mu}_{i,12}(\vv{A}):=&\hat{\nu}_{i,12}(\vv{A}) - \frac{1}{pn_i}\hat{\nu}_{i,4}(\vv{A}) \\
    =& \frac{1}{p}{\rm tr}\left(\vv{A}\left\{\left(\frac{1}{p}{\rm tr}(\vv{S}_{ip})\right)\vv{S}_{ip} - \frac{1}{p}\vv{R}_i + \frac{1}{n_i-1}\left(\frac{1}{p}{\rm tr}(\vv{S}_{ip})\right)\vv{S}_{ip}\right\}\right)\\
    =&\frac{1}{p}{\rm tr}\left(\vv{A}\left\{\frac{n_i}{n_i-1}\left(\frac{1}{p}{\rm tr}(\vv{S}_{ip})\right)\vv{S}_{ip} - \frac{1}{p}\vv{R}_{i}\right\}\right)\\
    =& \frac{1}{p}{\rm tr}(\vv{A}\vv{R}_{i,12}),
\end{align*}
in which $\vv{R}_{i,12}$ corresponds to (\ref{eq:R12}).

The unbiased estimator of $\mu_{i,2}(\vv{A})$
is
\begin{align*}
    \hat{\mu}_{i,2}(\vv{A}):=&\frac{n_i}{n_i-1}\hat{\nu}_{i,2}(\vv{A}) - \frac{p}{n_i-1}\hat{\nu}_{i,12}(\vv{A}) - \frac{1}{n_i}\hat{\nu}_{i,4}(\vv{A}) \\
    =& \frac{1}{p}{\rm tr}\left(\vv{A}\left\{\frac{n_i}{n_i-1}\vv{S}_{ip}^2 -\frac{1}{n_i-1}{\rm tr}(\vv{S}_{ip})\vv{S}_{ip} - \vv{R}_i + \frac{1}{n_i-1}{\rm tr}(\vv{S}_{ip})\vv{S}_{ip}\right\}\right)\\
    =& \frac{1}{p}{\rm tr}\left(\vv{A}\left\{\frac{n_i}{n_i-1}\vv{S}_{ip}^2  - \vv{R}_i \right\}\right)\\
    =& \frac{1}{p}{\rm tr}(\vv{A}\vv{R}_{i,2}),
\end{align*}
in which $\vv{R}_{i,2}$ corresponds to (\ref{eq:R2}).

Finally, the unbiased estimator of $\mu_{i,4}(\vv{A})$ is
\begin{align*}
    \hat{\mu}_{i,2}(\vv{A}):=&-\frac{2n_i}{n_i-1}\hat{\nu}_{i,2}(\vv{A}) + \frac{2p}{n_i-1}\hat{\nu}_{i,12}(\vv{A}) + \frac{n_i+2}{n_i}\hat{\nu}_{i,4}(\vv{A}) \\
    =& \frac{1}{p}{\rm tr}\left(\vv{A}\left\{\frac{2n_i}{n_i-1}\vv{S}_{ip}^2 - (n_i+2)\vv{R}_i + \frac{n_i}{n_i-1}{\rm tr}(\vv{S}_{ip})\vv{S}_{ip}\right\}\right)\\
    =& \frac{1}{p}{\rm tr}\left(\vv{A}\left\{\frac{n_i}{n_i-1}\vv{S}_{ip}^2  - \vv{R}_i \right\}\right)\\
    =& \frac{1}{p}{\rm tr}(\vv{A}\vv{R}_{i,4}),
\end{align*}
in which $\vv{R}_{i,4}$ corresponds to (\ref{eq:R4}). Hence, we complete the proof of Lemma \ref{lem:ub_est}. 
}

\section{Proofs in Section \ref{subsec:2.2}}
\subsection{Proofs of Theorem \ref{thm:2.4}}\label{app:thm2.4}

The main idea to derive the limiting distribution of a statistic resembling a $U$-statistic is to approximate it by its  Hajeck projection. (See \citesupp[Chapters 11 and 12]{vander2000} for details.) 

Let $\mathcal{S}_n$ be the set of all variables of the form $\sum_{i=1}^q f_{i}(\vv{X}_{ip})$  
for arbitrary measurable functions $f_i(\cdot)$ of the observations related to the $i$-th population with $\E[f_i(\vv{X}_{ip})^2]<\infty$. 
Then, the Haj\'{e}ck projection of $U_p$ satisfies
\begin{equation}\label{eq:hajeck_prop}
\begin{split}
\hat{U}_p&=\E(U_p|\mathcal{S}_n)-M_p=\sum_{i=1}^q \{\E(U_p|\vv{X}_{ip})-M_p\}\\
&=\frac{2}{q}\sum_{i=1}^q \left[\frac{1}{q-1}\sum_{j:j\neq i} \left\{\E(h(\vv{X}_{ip},\vv{X}_{jp})|\vv{X}_{ip})-d_{\rm prop}(\vv{\Sigma}_i,\vv{\Sigma}_j)\right\}\right]\\
&=:\frac{2}{q}\sum_{i=1}^q \left\{h_{1,i}(\vv{X}_{ip}) - \E[h_{1,i}(\vv{X}_{ip})]\right\},
\end{split}
\end{equation}
where
\begin{equation}\label{eq:h1i}
    h_{1,i}(\vv{X}_{ip}) = \frac{1}{q-1}\sum_{j:j\neq i} \E[h(\vv{X}_{ip},\vv{X}_{jp})|\vv{X}_{ip}].
\end{equation}
It can be easily seen that  $\hat{U}_p$ is a sum of independent and centered random variables. 

To show $\hat{U}_p$ is a good approximation to $U_p - M_p$, we observe that,
under Assumptions \ref{assm:2.1}---\ref{assm:2.3}, there exists a positive constant that bounds uniformly ${\rm Var}(h(\vv{X}_{ip},\vv{X}_{jp}))$ for any $1\leq i\neq j\leq q$. (The proof of the uniform boundedness is tedious and totally similar to the discussion for the variance ${\rm Var}(h_{1,i}(\vv{X}_{ip}))$ below, and thus omitted.) Thus, we can apply the similar argument for the $U$-statistic of order two as in Theorem 12.3 of \citesupp{vander2000} to obtain 
\[
\frac{{\rm var}(U_p)}{{\rm var}(\hat{U}_p)}\to 1,~~~\hbox{as}~p\to\infty.
\]
It then follows by Theorem 11.2 in \citesupp{vander2000} that
\[
\frac{U_p-M_p}{\sqrt{{\rm Var}(U_p)}}-\frac{\hat{U}_p}{\sqrt{{\rm Var}(\hat{U}_p)}}\overset{p.}{\to}0.
\]
Hence, it suffices for us to prove the asymptotic normality of $\hat{U}_p$.

Before proceeding to the asymptotic normality of $\hat{U}_p$, we first simplify the representation of $h_{1,i}(\vv{X}_{i,p})$. Observe that
\begin{equation}
    \begin{split}
        \E[h(\vv{X}_{ip},\vv{X}_{jp})|\vv{X}_{ip}] &= p\left\{\hat{\mu}_{i,2}\mu_{j,1}^2 + \hat{\mu}_{i,12}\mu_{j,2} - 2\hat{\mu}_{i,12}(\widetilde{\vv{\Sigma}_j})\right\},
    \end{split}
\end{equation}
where $\widetilde{\vv{\Sigma}_j}$ is given in (\ref{eq:tildeSgm}). 
Recall the definitions of $\alpha_{i,p}$, $\beta_{i,p}$ and $\vv{\Lambda}_{i,p}$ given in (\ref{eq:alpha_beta_Lambda}).
By (\ref{eq:inverse}), we further find 
\begin{equation}\label{eq:h_1i}
\begin{split}
h_{1,i}(\vv{X}_{ip})&=p\left\{\hat{\mu}_{i,2} \alpha_{i,p} + \hat{\mu}_{i,12}\beta_{i,p} - 2\hat{\mu}_{i,12}(\vv{\Lambda}_{i,p})\right\} \\
&= p\left\{(\hat{\nu}_{i,2}(\vv{I}_p) - c_{ip}\hat{\nu}_{i,12}(\vv{I}_p))\alpha_{i,p} + \hat{\nu}_{i,12}(\vv{I}_p)\beta_{i,p}   -2\hat{\nu}_{i,12}(\vv{\Lambda}_{i,p})\right\}+ R_{1,i}\\
& =: H_{1,i} +R_{1,i}, 
\end{split}
\end{equation}
where the remainder satisfies
\begin{equation}\label{eq:r_1i}
\begin{split}
R_{1,i}(\vv{X}_{ip})
&= \frac{p}{n_i-1}\left\{\hat{\nu}_{i,2}(\vv{I}_p) - c_{ip}\hat{\nu}_{i,12}(\vv{I}_p)\right\}\alpha_{i,p} + \frac{p}{n_i}\hat{\nu}_{i,4}(\vv{I}_p)\alpha_{i,p} \\
& + \frac{1}{n_i}\hat{\nu}_{i,4}(\vv{I}_p)\beta_{i,p} - 2\frac{p}{n_i}\hat{\nu}_{i,4}(\vv{\Lambda}_{i,p}).
\end{split}
\end{equation}
The Assumptions \ref{assm:2.1}---\ref{assm:2.3} ensure that variances of $\{p(\hat{\nu}_{i,2}(\vv{I}_p) - c_{ip}\hat{\nu}_{i,12}(\vv{I}_p))\}_{1\leq i\leq q}$, $\{p\hat{\nu}_{i,4}(\vv{I}_p)\}_{1\leq i\leq q}$ and $\{p\hat{\nu}_{i,4}(\vv{\Lambda}_{i,p})\}_{1\leq i\leq q}$ are all uniformly bounded by some positive constant $C$ independent of $p,q$ and $\{n_j\}$. Therefore, variances of the remainders $\{R_{1,i}\}_{1\leq i\leq q}$ are also uniformly bounded and of order $O(p^{-2})$, that is, $\max_{1\leq i\leq q} {\rm Var}(R_{1,i}) \leq C'/p^2$ for another positive number, say, $C'$. Consequently, we see that
\[
\hat{U}_p = \frac{2}{q}\sum_{i=1}^q \{H_{1,i}-\E[H_{1,i}]\} + O_p(\frac{1}{\sqrt{q}p}).
\]
Thus, it only remains to show the asymptotic normality of the sum of $H_{1,i}$'s.
Observe that $(H_{1,i}-\E[H_{1,i}])$'s are mutually independent, we immediately have $\hat{U}_p/\sqrt{{\rm var}(\hat{U}_p)}\overset{d.}{\to}\mathcal{N}(0,1)$ provided that the Lyaponouv condition
\begin{equation}\label{eq:b12}
\sum_{i=1}^q \E|H_{1,i} - \E[H_{1,i}]|^4\bigg/\left(\sum_{i=1}^q {\rm Var}(H_{1,i})^2\right)^2\to 0
\end{equation}
holds. Hence, it remains (i) to derive the explicit representation of ${\rm Var}(\hat{U}_{p})$ in terms of the spectral moments of $\{\vv{\Sigma}_j\}$, and (ii) to verify the Lyaponouv condition (\ref{eq:b12}). 

Let ${\sigma}_{i,p}^2={\rm Var}(H_{1,i})$. Then,
\[
{\rm Var}(\hat{U}_p) = \frac{4}{q^2} \sum_{i=1}^q \sigma_{i,p}^2 + O(\frac{1}{qp^2}).
\]
Hence, to achieve the goal (i), it suffices to find the explicit formula of $\sigma_{i,p}^2$. 

Recall the definition of $H_{1,i}$ in (\ref{eq:h_1i}). Applying Lemma \ref{lem:jointcov}, we have
\begin{align*}
    \sigma_{i,p}^2 & = \alpha_{i,p}^2 {\rm Var}(p\{\hat{\nu}_{i,2}(\vv{I}_p) - c_{ip}\hat{\nu}_{i,12}(\vv{I}_p)\}) \\
    &+ \beta_{i,p}^2 {\rm Var}(p\hat{\nu}_{i,12}(\vv{I}_p)) 
    + 4{\rm Var}(p\hat{\nu}_{i,12}(\vv{\Lambda}_{i,p})) \\
    &+ 2\alpha_{i,p}\beta_{i,p}{\rm Cov}(p\{\hat{\nu}_{i,2}(\vv{I}_p) - c_{ip}\hat{\nu}_{i,12}(\vv{I}_p)\},p\hat{\nu}_{i,12}(\vv{I}_p))\\
    & -4 \alpha_{i,p}{\rm Cov}(p\{\hat{\nu}_{i,2}(\vv{I}_p) - c_{ip}\hat{\nu}_{i,12}(\vv{I}_p)\},p\hat{\nu}_{i,12}(\vv{\Lambda}_{i,p})) \\
    &- 4 \beta_{i,p}{\rm Cov}(p\hat{\nu}_{i,12}(\vv{I}_p),p\hat{\nu}_{i,12}(\vv{\Lambda}_{i,p}))\\
    & = \alpha_{i,p}^2\left\{4c_{ip}^2\mu_{i,2}^2 + 4c_{ip}\langle\vv{\Sigma}_i, \vv{\Sigma}_i\rangle_{\vv{\Sigma}_i}\right\}
     + 4c_{ip}\beta_{i,p}^2 \mu_{i,1}^2 \langle\vv{I}_p, \vv{I}_p\rangle_{\vv{\Sigma}_i}\\
    & + 4 c_{ip}\langle\mu_{i,1}\vv{\Lambda}_{i,p}+\kappa_{1,i}\vv{I}_p,\mu_{i,1}\vv{\Lambda}_{i,p}+\kappa_{1,i}\vv{I}_p\rangle_{\vv{\Sigma}_i}\\
        & + 8c_{ip}\alpha_{i,p}\beta_{i,p} \mu_{i,1}\langle\vv{\Sigma}_i,\vv{I}_p\rangle_{\vv{\Sigma}_i}\\
        & - 8c_{ip}\alpha_{i,p}\left\{\mu_{i,1}\langle\vv{\Sigma}_i,\vv{\Lambda}_{i,p}\rangle_{\vv{\Sigma}_i} + \kappa_{1,i}\langle\vv{\Sigma}_i,\vv{I}_p\rangle_{\vv{\Sigma}_i}\right\}\\
        & - 8c_{ip}\beta_{i,p} \mu_{i,1}\left\{\mu_{i,1}\mu(\vv{\Lambda}_{i,p},\vv{I}_p|\vv{\Sigma}_i) + \kappa_{1,i}\mu(\vv{I}_p,\vv{I}_p|\vv{\Sigma}_i)\right\}+ r_{i,p},
\end{align*}
in which, under Assumptions \ref{assm:2.1}---\ref{assm:2.3}, the remainders $\{r_{i,p}\}_{1\leq i\leq q}$ can be uniformly bounded and of order $O(1/p)$, i.e., $\max_{1\leq i\leq q}|r_{i,p}| \leq C''/p$ for some constant $C''$. Adopting the notation $\vv{\Gamma}_{i,p}$ in (\ref{eq:Gamma_ip}), we can then simplify the variance $\sigma_{i,p}^2$ to be
\[
\sigma_{i,p}^2  = 4c_{ip}^2 \alpha_{i,p}^2\mu_{i,2}^2
         + 4c_{ip}\langle\vv{\Gamma}_{i,p},\vv{\Gamma}_{i,p}\rangle_{\vv{\Sigma}_i} + r_{i,p}.
\]

Since both $\{c_{ip}\}$ and $\{\vv{\Sigma}_j\}$ are uniformly bounded under Assumptions \ref{assm:2.1}---\ref{assm:2.3}, it is easy to see that 
\[
\sigma_p^2 = \sum_{i=1}^q \sigma_{i,p}^2 = O(q).
\]
Also, by Assumption \ref{assm:2.2}, we have
\begin{align*}
\alpha_{i,p} &= \frac{q-1}{q}\left\{\frac{1}{q}\sum_{j=1}^q \mu_{j,1}^2 - \frac{\mu_{j,1}}{q}\right\} \geq  \frac{q-1}{q}\left\{\left(\frac{1}{q}\sum_{j=1}^q \mu_{j,1}\right)^2 - \frac{\mu_{j,1}}{q}\right\}\\
& \geq \frac{q-1}{q}\left(c^2 - \frac{C}{q}\right) >0.
\end{align*}
and thus
\begin{align*}
    \frac{1}{q}\sum_{i=1}^q \sigma_{i,p}^2 & \geq \frac{1}{q}\sum_{i=1}^q 4c_{ip}^2 \alpha_{i,p}^2\mu_{i,2}^2 \geq \frac{4c_0^2(q-1)}{q}\left(c^2 - \frac{C}{q}\right)\left(\frac{1}{q}\sum_{i=1}^q \mu_{i,2}^2\right) \\
    & \geq  \frac{4c_0^2(q-1)}{q}\left(c^2 - \frac{C}{q}\right)\left(\frac{1}{q}\sum_{i=1}^q \mu_{i,1}\right)^4 \\
    & \geq \frac{4c_0^2(q-1)}{q}\left(c^2 - \frac{C}{q}\right) c^4 >0
\end{align*}
It suggests that the denominator in (\ref{eq:b12}) satisfies
\[
c' q^2 \leq \left(\sum_{i=1}^q {\rm Var}(H_{1,i})\right)^2 \leq q^2C'''
\]
for some positive constants $c'$ and $C'''$.

Meanwhile, according to Lemma \ref{lem:fourthmom}, together with Assumptions \ref{assm:2.1}---\ref{assm:2.3}, there exists some constant $A_0$ such that
\begin{equation}\label{eq:B}
\sum_{i=1}^q \E|h_{1,i}(\vv{X}_{ip})|^4\leq q\cdot A_0.
\end{equation}
which means the numerator in (\ref{eq:b12}) has the order of $O(q)$. Hence, the Lyaponouv condition (\ref{eq:b12}) indeed holds. 

Combining the discussions above, we finally conclude that \[
(U_p-M_p)/\sqrt{{\rm Var}(U_p)}\overset{d.}{\to}\mathcal{N}(0,1),
\]
where  
\begin{align*}
{\rm Var}(U_p)&=\frac{4}{q^2}\sum_{i=1}^q \sigma_{i,p}^2+o(1)\\
&=\frac{1}{q}\left\{\frac{4}{q}\sum_{i=1}^q \left\{4c_{ip}^2 \alpha_{i,p}^2\mu_{i,2}^2
         + 4c_{ip}\langle\vv{\Gamma}_{i,p},\vv{\Gamma}_{i,p}\rangle_{\vv{\Sigma}_i}\right\}+ O(\frac{1}{p})\right\} +o(1)\\
&=\frac{1}{q}\sigma_{p}^2+o(1).
\end{align*}
The proof is complete.

\subsection{Proof of Corollary \ref{cor:2.5}}\label{sec:cor_prop}
According to Theorem \ref{thm:2.4}, it suffices for us to verify that $\vv{\Gamma}_{i,p} = \vv{O}$ under $H_0$. In fact, when $H_0$ holds, we define a basis matrix $\vv{\Sigma}$ such that $\vv{\Sigma}_i = \mu_{i,1}\vv{\Sigma}$ and $p^{-1}{\rm tr}(\vv{\Sigma}) = 1$. The matrix $\vv{\Lambda}_{i,p}$ in Theorem \ref{thm:2.4} now becomes
\begin{align*}
    \vv{\Lambda}_{i,p} = \left(\frac{1}{q-1}\sum_{j:j\neq i}\mu_{j,1}^2\right)\vv{\Sigma} = \alpha_{i,p}\vv{\Sigma}.
\end{align*}
It then follows that
\begin{align*}
    \kappa_{1,i} &= \frac{1}{p}{\rm tr}(\vv{\Lambda}_{i,p}\vv{\Sigma}_i) = \alpha_{i,p} \mu_{i,1} \left(\frac{1}{p}{\rm tr}(\vv{\Sigma}^2)\right) \\
    & = \mu_{i,1}\left\{\frac{1}{q-1}\sum_{j:j\neq i}\mu_{j,1}^2 \left(\frac{1}{p}{\rm tr}(\vv{\Sigma}^2)\right)\right\}\\
    & =\mu_{i,1}\left\{\frac{1}{q-1}\sum_{j:j\neq i} \left(\frac{1}{p}{\rm tr}(\vv{\Sigma}_j^2)\right)\right\}\\
    & = \mu_{i,1}\beta_{i,p}
\end{align*}
Hence, by (\ref{eq:Gamma_ip}), we have $\vv{\Gamma}_{i,p} = \vv{O}$ and thus the inner product $\langle \vv{\Gamma}_{i,p},\vv{\Gamma}_{i,p}\rangle_{\vv{\Sigma}_i} = 0$. Observe that $\sigma_p^2 =\sigma_{0,p}^2 + \sigma_{r,p}^2$. We finally see that $\sigma_p^2=\sigma_{0,p}^2$.

\subsection{Proof of Theorem \ref{thm:2.5}}\label{app:thm2.5}

According to Theorem \ref{thm:2.4}, it only remains to show that $\hat{\sigma}_p^2$ is consistent to $\sigma_{0,p}^2$. From the discussion in the last Section \ref{app:thm2.4}, we have already seen that both $\{{\rm Var}(p\hat{\mu}_{i,2})\}$ and $\{{\rm Var}(p\hat{\mu}_{i,12})\}$ can be uniformly bounded under Assumptions \ref{assm:2.1}---\ref{assm:2.3}. In other words, there is a positive constant  $A_1$ independent of $p,q$ and $\{n_j\}$ such that 
\[
\max_{1\leq i\leq q} {\rm Var}(p\hat{\mu}_{i,2}) + {\rm Var}(p\hat{\mu}_{i,12}) \leq A_1.
\]

Meanwhile, we recall that $\hat{\mu}_{i,2}$ and $\hat{\mu}_{i,12}$ are normailzed spectral moments of $\vv{R}_{i,2}$ and $\vv{R}_{i,12}$ as in (\ref{eq:R2}) and (\ref{eq:R12}), which are linear combinations of matrices
\[
\frac{1}{n_i(n_i-1)}\vv{X}_{ip}{\rm diag}(\vv{X}_{ip}^\T\vv{X}_{ip})\vv{X}_{ip}^\T,~~~\left(\frac{1}{p}{\rm tr}(\vv{S}_{ip})\right)\vv{S}_{ip},~~~\hbox{and}~~~\vv{S}_{ip}^2.
\]
Since $\{\vv{\Sigma}_i\}$ are uniformly bounded in Assumption \ref{assm:2.2}, we know from \cite[Theorem 5.8]{BS10} that $\|\vv{X}_{ip}\|^2/n_i = \|\vv{S}_{ip}\|$ are also uniformly bounded with probability one. 
Consequently, there exists another constant $M_0$ such that for any $i=1,\ldots,q$,
\[
\max_{1\leq i\leq q} \{\hat{\mu}_{i,2}^2 + \mu_{i,2}^2+ \hat{\mu}_{i,12}^2 +\mu_{i,1}^2\}\leq M_0.
\]
It then follows that
\begin{align*}
\bigg|\frac{1}{q}\sum_{j=1}^q c_{jp}^2\left\{\hat{\mu}_{i,2}^2-\mu_{i,2}^2\right\}\bigg|
\leq\frac{M_0}{q}\sum_{j=1}^q c_{jp}^2\left|\hat{\mu}_{i,2}-\mu_{i,2}\right|
\end{align*}
and 
\begin{align*}
    \E\bigg|\frac{1}{q}\sum_{j=1}^q c_{jp}^2\left\{\hat{\mu}_{i,2}^2-\mu_{i,2}^2\right\}\bigg|
&\leq\frac{M_0}{q}\sum_{j=1}^q c_{jp}^2\E\left|\hat{\mu}_{i,2}-\mu_{i,2}\right| \\
&\leq \frac{M_0}{qp}\sum_{j=1}^q c_{jp}^2\sqrt{{\rm Var}(p\hat{\mu}_{i,2})} \leq \frac{C_0^2M_0A_1}{p} \to 0.
\end{align*}
Hence, we have
\[
\frac{1}{q}\sum_{j=1}^q c_{jp}^2\left|\hat{\mu}_{i,2}^2 - \mu_{i,2}^2\right|\overset{p.}{\to}0.
\]

Similarly,
\begin{align*}
    \E\bigg|\frac{1}{q}\sum_{j=1}^q \left\{\hat{\mu}_{i,12}-\mu_{i,12}\right\}\bigg|
&\leq\frac{1}{q}\sum_{j=1}^q \E\left|\hat{\mu}_{i,12}-\mu_{i,12}\right| \\
&\leq \frac{1}{qp}\sum_{j=1}^q \sqrt{{\rm Var}(p\hat{\mu}_{i,12})} \leq \frac{M_0A_1}{p} \to 0.
\end{align*}
Thus, we also have
\[
\frac{1}{q}\sum_{j=1}^q |\hat{\mu}_{i,12}-\mu_{i,12}|\overset{p.}{\to}0.
\]

Combining the above results with (\ref{eq:sigma_0_app}), we see that $\hat{\sigma}_p^2$ is indeed a consistent estimator for $\sigma_{0,p}^2$. Thus, the proof is complete.

\section{Proof of Theorem \ref{thm:generalpower}}\label{sec:proof_generalpower}
{
    According to Theorem \ref{thm:2.4}, $U_p$ is still asymptotically normal under alternative, that is,
    \[
        \frac{\sqrt{q}(U_p - M_p)}{\sigma_p} \overset{d.}{\to}\mathcal{N}(0,1).
    \]

    Still, according to Theorem \ref{thm:2.5}, $\hat{\sigma}_p^2$ is always consistent to $\sigma_{0,p}^2$ whenever $H_0$ holds or not. Hence, $\hat{\sigma}_{p}/\sigma_{0,p}\overset{d.}{\to} 1$.

    Observe that
    \begin{align*}
        \P_{H_1}\left(\frac{\sqrt{q}U_p}{\hat{\sigma}_p} > z_\alpha\right) &= \P_{H_1}\left(\frac{\sqrt{q}(U_p-M_p)}{\sigma_p} > \frac{\hat{\sigma}_p}{\sigma_p} z_\alpha - \frac{\sqrt{q}M_p}{\sigma_p}\right) \\
        & = \Phi\left(\frac{\sqrt{q}M_p}{\sigma_p} - \frac{\sigma_{0,p}}{\sigma_p} z_\alpha\right) + o(1),
    \end{align*}
    which proves the first conclusion.

    To see the second conclusion, we notice that $\sigma_p^2$ is upper bounded under Assumption \ref{assm:2.2}. Hence, there exists some positive constant $C$ such that $\sigma_p^2\leq C'$. Meanwhile, since $\sigma_p^2 \geq \sigma_{0,p}^2$ and the later has a positive lower bound $c'$ under Assumption \ref{assm:2.1}. In short, we have $0<c'\leq \sigma_p^2 \leq C'$, where $c',C'$ are positive constants.
    With this fact in mind, we immediately have 
    \[
    \frac{\sqrt{q}M_p}{\sqrt{c'}}\geq \frac{\sqrt{q}M_p}{\sigma_p} \geq \frac{\sqrt{q}M_p}{\sqrt{C'}}.
    \]
    Meanwhile, since $\sigma_{0,p}/\sigma_p\leq 1$, we have
    \[
    \left|\frac{\sigma_{0,p}}{\sigma_p} z_\alpha\right| \leq |z_\alpha|,
    \]
    where the right-hand side remains unchanged when $p,q$ and $\{n_j\}$ increase. Consequently, we have
    \begin{align*}
        \Phi\left(\frac{\sqrt{q}M_p}{\sqrt{c'}} + |z_\alpha|\right)& \geq \Phi\left(\frac{\sqrt{q}M_p}{\sigma_p} - \frac{\sigma_{0,p}}{\sigma_p} z_\alpha\right)
        \geq \Phi\left(\frac{\sqrt{q}M_p}{\sqrt{C'}} - |z_\alpha|\right).
    \end{align*}
    From the above discussion, it is not difficult to see that the power function (\ref{eq:pf}) tends to $1$ if and only if $\sqrt{q}M_p \to \infty$.
    The proof is then complete. 

}

\section{Proof of Lemma \ref{lem:kron}}\label{app:4}

``(i) $\Longrightarrow$ (iii)'':
 Assume that $\tilde{\vv{X}}$ satisfies (\ref{eq:kron}) with $\vv{\Sigma}_C$ diagonal. Let $\vv{x}_1,\ldots,\vv{x}_q$ and $\vv{z}_1,\ldots,\vv{z}_q$ be columns of $\tilde{\vv{X}}$ and $\vv{Z}$ respectively. Denote $a_{ij}$ be the $(i,j)$th entry of $\vv{\Sigma}_C^{1/2}$ and write $\vv{a}_{(i)}$ be the $i$th column of $\vv{\Sigma}_C$. Since $\tilde{\vv{X}}$ satisfies (\ref{eq:kron}), we have
\[
\tilde{\vv{X}}=\vv{\Sigma}_R^{1/2}\vv{Z}\vv{\Sigma}_C^{1/2}=(\vv{\Sigma}_R^{1/2}\vv{z}_1,\ldots,\vv{\Sigma}_R^{1/2}\vv{z}_q)(\vv{a}_{(1)},\cdots,\vv{a}_{(q)}).
\]
Hence, columns of $\tilde{\vv{X}}$ can be represented by linear combinations of $\vv{\Sigma}_R^{1/2}\vv{z}_j$'s with the coefficients the corresponding columns of $\vv{\Sigma}_C^{1/2}$, that is,
\[
\vv{x}_i=\sum_{j=1}^q a_{ji}\vv{\Sigma}^{1/2}_R\vv{z}_j,~~~~i=1,\ldots,q.
\]
When $\vv{\Sigma}_C$ is diagonal, so is $\vv{\Sigma}_C^{1/2}$. Then, we have $\vv{x}_i=a_{ii}\vv{\Sigma}_R^{1/2} z_i$. It  follows that $\vv{x}_1,\ldots,\vv{x}_q$ are independent and their covariance matrices differ only by a factor.

``(iii)$\Longrightarrow$(ii)'': Assume that $\vv{x}_1,\ldots,\vv{x}_q$ are independent with their covariance matrices proportional to each other. Divide $\vv{\Sigma}^{1/2}$ into $q^2$ blocks of size $p\times p$ and denote the $(i,j)$th block by $\vv{A}_{ij}$. Correspondingly, denote $\vv{\Sigma}_{ij}$ the $(i,j)$th block of $\vv{\Sigma}$. Then,
\[
\vv{\Sigma}_{ij}=\sum_{l=1}^q \vv{A}_{il}\vv{A}^{\T}_{jl},~~~i,j=1,\ldots,q,
\]
and
\[
\begin{pmatrix}\vv{x}_1\\\vv{x}_2\\\vdots\\\vv{x}_q\end{pmatrix}={\rm vec}(\tilde{\vv{X}})=\vv{\Sigma}^{1/2}\vv{z}=\begin{pmatrix}
\vv{A}_{11}&\vv{A}_{12}&\cdots&\vv{A}_{1q}\\
\vv{A}_{21}&\vv{A}_{22}&\cdots&\vv{A}_{2q}\\
\vdots&\vdots&\ddots&\vdots\\
\vv{A}_{q1}&\vv{A}_{q2}&\cdots&\vv{A}_{qq}\end{pmatrix}\begin{pmatrix}\vv{z}_1\\\vv{z}_2\\\vdots\\\vv{z}_q\end{pmatrix}.
\]
Since $\vv{x}_i$ and $\vv{x}_j$ are independent for any $i\neq j$, we have $\vv{\Sigma}_{ij}=\vv{O}$ for all $i\neq j$. Hence, $\vv{A}_{ij}=\vv{O}$ for $i\neq j$, $\vv{A}_{ii}=\vv{\Sigma}_{ii}^{1/2}$ and $\vv{x}_i=\Sigma_{ii}^{1/2} \vv{z}_i$. Moreover, since covariance matrices of $\vv{x}_1,\ldots,\vv{x}_q$ are proportional to each other, there must exists a sequence of $w_1,\ldots,w_q$ and a $p\times p$ matrix $\vv{\Sigma}_1$ such that $\vv{\Sigma}_{ii}=w_i^2\vv{\Sigma}_1$. Let $\vv{\Sigma}_R=\vv{\Sigma}_1$ and $\vv{\Sigma}_C={\rm diag}(w_1^2,\ldots,w_q^2)$. We have $\vv{\Sigma}=\vv{\Sigma}_C\otimes \vv{\Sigma}_R$ and therefore $\tilde{\vv{X}}=\vv{\Sigma}_R^{1/2}\vv{Z}\vv{\Sigma}_C^{1/2}$.

``(ii) $\Longrightarrow$ (i)'': Assume that $\vv{\tilde{X}}$ satisfies (\ref{eq:kron}) with independent columns. Note that 
\[
{\rm cov}(\tilde{\vv{x}}_{j},\tilde{\vv{x}}_k) = (\vv{\Sigma}_C)_{jk} \vv{\Sigma}_R,~~~~j,k\in\{1,\ldots,q\}
\]
Due to the column independence, ${\rm cov}(\tilde{\vv{x}}_j,\tilde{\vv{x}}_k)=\vv{0}$ for any $j\neq k$. Hence, $(\vv{\Sigma}_C)_{jk}=0$ for any $j\neq k$. In other words, $\vv{\Sigma}_C$ is diagonal.

\section{Proofs in Section \ref{subsec:3.2}}
\subsection{Proof of Theorem \ref{thm:3.3}}\label{app:2}

The proof is roughly parallel to that in Section \ref{app:thm2.4}, by the method of Haj\'{e}ck projection. The Haj\'{e}ck projection of $V_p$ satisfies
\begin{equation}\label{eq:b31}
\begin{split}
\hat{V}_p&=\sum_{i=1}^q (\E(V_p|\vv{X}_{ip})-m_p)\\
&=\frac{2}{q}\sum_{i=1}^q \left[\frac{1}{q-1}\sum_{j:j\neq i} \{\E(g(\vv{X}_{ip},\vv{X}_{jp})|\vv{X}_{ip})- d_0(\vv{\Sigma}_i,\vv{\Sigma}_j)\}\right]\\
&=:\frac{2}{q}\sum_{i=1}^q \{g_{1,i}(\vv{X}_{ip}) -\E[g_{1,i}(\vv{X}_{ip})]\},
\end{split}
\end{equation}
where 
\begin{equation}
   g_{1,i}(\vv{X}_{ip}) = \frac{1}{q-1}\sum_{j:j\neq i} \E(g(\vv{X}_{ip},\vv{X}_{jp})|\vv{X}_{ip}).
\end{equation}
It can be easily seen that $\hat{V}_p$ is a sum of independent variables with zero mean. 

Under Assumptions \ref{assm:2.1}---\ref{assm:2.4}, ${\rm var}(g(X_{ip},X_{jp}))$ are bounded uniformly for all $1\leq i\neq j\leq q$. Thus, we apply  a  similar argument as in Theorem 12.3 and Theorem 11.2 in \citesupp{vander2000} to have 
\[
\frac{{\rm Var}(V_p)}{{\rm Var}(\hat{V}_p)}\to 1,~~~\hbox{and}~~~
\frac{V_p-m_p}{\sqrt{{\rm Var}(V_p)}}-\frac{\hat{V}_p}{\sqrt{{\rm Var}(\hat{V}_p)}}\overset{p}{\to}0, ~~~\hbox{as}~p\to\infty.
\]

To simplify the representation of $g_{1,i}(\vv{X}_{i,p})$, we observe that
\begin{equation}\label{eq:g1i}
    \E[g(\vv{X}_{ip},\vv{X}_{jp})|\vv{X}_{ip}] = p\hat{\mu}_{i,2} + p\mu_{j,2} - 2{\rm tr}(\vv{S}_{ip}\vv{\Sigma}_j).
\end{equation}
Let 
\[
\tilde{\vv{\Lambda}}_{i,p} = \frac{1}{q-1}\sum_{j:j\neq i} \vv{\Sigma}_j.
\]
Then, $\tilde{\vv{\Gamma}}_{i,p} = \vv{\Sigma}_i - \tilde{\vv{\Lambda}}_{i,p}$. 
With this notation, we immediately find that
\begin{equation}
    \begin{split}
        g_{1,i}(\vv{X}_{ip}) &= p\hat{\mu}_{i,2} + p\left(\frac{1}{q-1}\sum_{j:j\neq i}\mu_{j,2}\right) - 2{\rm tr}(\vv{S}_{ip}\tilde{\vv{\Lambda}}_{i,p})\\
        & = p\hat{\nu}_{i,2}(\vv{I}_p) + p\left(\frac{1}{q-1}\sum_{j:j\neq i}\mu_{j,2}\right)- 2{\rm tr}(\vv{S}_{ip}\tilde{\vv{\Lambda}}_{i,p})+\tilde{R}_{1,i}\\
        & =: G_{1,i}+\tilde{R}_{1,i},
    \end{split}
\end{equation}
where the remainder satisfies
\begin{equation}
    \begin{split}
        \tilde{R}_{1,i} &=\frac{p}{n_i-1}\{\hat{\nu}_{i,2}(\vv{I}_p) - c_{ip}\hat{\nu}_{i,12}(\vv{I}_p)\}+ \frac{p}{n_i}\hat{\nu}_{4,i}(\vv{I}_p).
    \end{split}
\end{equation}
Under Assumptions \ref{assm:2.2}---\ref{assm:2.3}, the variances of $\tilde{R}_{1,i}$ can be uniformly bounded in the form: $\max_{1\leq i\leq q}{\rm Var}(\tilde{R}_{1,i})\leq C'/p^2$. Hence, we see that
\[
\hat{V}_p = \frac{2}{q}\sum_{i=1}^q \{G_{1,i} - \E[G_{1,i}]\} + O_p(\frac{1}{\sqrt{q}p}).
\]
Consequently,
to obtain $\hat{V}_p/\sqrt{{\rm var}(\hat{V}_p)}\overset{d.}{\to}\mathcal{N}(0,1)$, it suffices to verify the Lyaponouv condition
\begin{equation}\label{eq:c12}
\sum_{i=1}^q \E|G_{1,i}-\E[G_{1,i}]|^4\bigg/\left(\sum_{i=1}^q {\rm var}(G_{1,i})^2\right)^2\to 0.
\end{equation}

Let $\lambda_{i,p}^2={\rm Var}(G_{1,i})$. Then, we have
\[
{\rm Var}(\hat{V}_p) = \frac{4}{q^2}\sum_{i=1}^q \lambda_{i,p}^2 + O(\frac{1}{qp^2}).
\]
Moreover, by Lemmas \ref{lem:jointcov} and \ref{lem:single_cov}, we have
\begin{align*}
    \lambda_{i,p}^2 &= {\rm Var}(p\hat{\nu}_{i,2}(\vv{I}_p))+4{\rm Var}({\rm tr}(\vv{S}_{ip}\tilde{\vv{\Lambda}}_{i,p})) 
     - 4 {\rm Cov}(p\hat{\nu}_{i,2}(\vv{I}_p), {\rm tr}(\vv{S}_{ip}\tilde{\vv{\Gamma}}_{i,p})) \\
    & = 4c_{ip}^2\mu_{i,2}^2 + 4c_{ip}\langle\vv{\Sigma}_i,\vv{\Sigma}_i\rangle_{\vv{\Sigma}_i} + 4c_{ip}\langle\tilde{\vv{\Lambda}}_{i,p},\tilde{\vv{\Lambda}}_{i,p}\rangle_{\vv{\Sigma}_i} 
    - 8c_{ip}\langle\vv{\Sigma}_i,\tilde{\vv{\Lambda}}_{i,p}\rangle_{\vv{\Sigma}_i} + \tilde{r}_{i,p}\\
    & = 4c_{ip}^2 \mu_{i,2}^2 + 4c_{ip}\langle\vv{\Sigma}_i -\tilde{\vv{\Lambda}}_{i,p},\vv{\Sigma}_i -\tilde{\vv{\Lambda}}_{i,p}\rangle_{\vv{\Sigma}_i}+\tilde{r}_{i,p}\\
    & =4c_{ip}^2 \mu_{i,2}^2 + 4c_{ip}\langle\tilde{\vv{\Gamma}}_{i,p},\tilde{\vv{\Gamma}}_{i,p}\rangle_{\vv{\Sigma}_i}+\tilde{r}_{i,p},
\end{align*}
in which $\max_{1\leq i\leq q}\tilde{r}_{i,p} \leq C''/p$.
Therefore, 
\begin{align*}
\frac{1}{q}\sum_{i=1}^q \lambda_{i,p}^2 & \geq \left(\frac{1}{q}\sum_{i=1}^{q} 4c_{i,p}^2 \mu_{i,2}^2\right)(1+o(1))\\
&\geq 4c_0^2 \left(\frac{1}{q}\sum_{i=1}^{q} \mu_{i,1}\right)^2(1+o(1))> 4c_0^2 c^2 (1 + o(1))>0.
\end{align*}
It suggests that the denominator in (\ref{eq:c12}) has an order $q^2$.

Meanwhile, similar to the discussion of (\ref{eq:B}), it also holds
\begin{equation}\label{eq:D}
\sum_{i=1}^q \E|G_{1,i} - \E[G_{1,i}]|^4=O(q).
\end{equation}
Hence, the Lyaponouv condition (\ref{eq:c12}) indeed holds and we have
\[
(V_p-m_p)/\sqrt{{\rm Var}(V_p)}\overset{d.}{\to}\mathcal{N}(0,1)
\]
with  
\begin{align*}
{\rm Var}(V_p)&=\frac{4}{q^2}\sum_{i=1}^q {\lambda}_{i,p}^2+o(1)\\
&=\frac{16}{q}\left\{\frac{4}{q}\sum_{i=1}^q \{4c_{ip}^2 \mu_{i,2}^2 + 4c_{ip}\langle\tilde{\vv{\Gamma}}_{i,p},\tilde{\vv{\Gamma}}_{i,p}\rangle_{\vv{\Sigma}_i}\}+O(\frac{1}{p})\right
\}+o(1)\\
&=\frac{1}{q}{\lambda}_{p}^2+o(1).
\end{align*}
The proof is then complete.

\subsection{Proof of Theorem \ref{thm:3.4}}\label{sec:proof_thm3.5}

The proof is, in general, parallel to that in Section \ref{app:thm2.5}. Here, we show  the consistency of $\hat{\lambda}_p^2$ to $\lambda_{0,p}^2$, i.e., $\hat{\lambda}_p^2-\lambda_{0,p}^2\overset{p}{\to}0$.  
Recall the discussion in Section \ref{app:thm2.5}, $\{\hat{\mu}_{i,2}\}_{1\leq i\leq q}$ are uniformly bounded with probability one, that is, there exists a uniform constant $M_1$ such that for any $i=1,\ldots,q$,
\begin{align*}
  \max_{1\leq i\leq q}\{\hat{\mu}_{i,2}+\mu_{i,2}\}\leq M_1.
\end{align*}
It then follows that
\begin{align*}
\left|\frac{1}{q}\sum_{j=1}^q c_{jp}^2(\hat{\mu}_{j,2}^2 -\mu_{j,2}^2)\right|
\leq \frac{M_1}{q}\sum_{j=1}^q c_{jp}^2\left|\hat{\mu}_{j,2}-\mu_{j,2}\right|,
\end{align*}
so that
\begin{align*}
\E\left|\frac{1}{q}\sum_{j=1}^q c_{jp}^2(\hat{\mu}_{j,2}^2 -\mu_{j,2}^2)\right|&\leq \frac{C^2_0M_1}{pq}\sum_{j=1}^q \sqrt{{\rm Var}(p\hat{\mu}_{i,2})} \leq \frac{C^2_0A_0M_1}{p} \to 0.
\end{align*}
In other words,
\[
\frac{1}{q}\sum_{j=1}^q c_{jp}^2(\hat{\mu}_{j,2}^2 -\mu_{j,2}^2) \overset{p.}{\to} 0.
\]
The proof is then complete.

\section{Proof of Theorem \ref{thm:generalpower_eq}}\label{sec:proof_generalpower_eq}
{
    According to Theorem \ref{thm:3.3}, $U_p$ is still asymptotically normal under alternative, that is,
    \[
        \frac{\sqrt{q}(V_p - M_p)}{\sigma_p} \overset{d.}{\to}\mathcal{N}(0,1).
    \]
    Meanwhile, by Theorem \ref{thm:3.4}, $\hat{\lambda}_p^2$ is always consistent to $\lambda_{0,p}^2$ in general. In other words, $\hat{\lambda}_{p}/\lambda_{0,p}\overset{d.}{\to} 1$.

    Observe that
    \begin{align*}
        \P_{H_1}\left(\frac{\sqrt{q}V_p}{\hat{\lambda}_p} > z_\alpha\right) &= \P_{H_1}\left(\frac{\sqrt{q}(U_p-m_p)}{\lambda_p} > \frac{\hat{\lambda}_p}{\lambda_p} z_\alpha - \frac{\sqrt{q}m_p}{\lambda_p}\right) \\
        & = \Phi\left(\frac{\sqrt{q}m_p}{\lambda_p} - \frac{\lambda_{0,p}}{\lambda_p} z_\alpha\right) + o(1).
    \end{align*}
    So, the first conclusion is proved.

    As for the second conclusion, recall that $\lambda_p^2$ is upper bounded under Assumption \ref{assm:2.2}. In other words, there must be a positive constant $C'$ such that $\lambda_p^2\leq C'$. Meanwhile, we also have $\lambda_p^2 \geq \lambda_{0,p}^2 \geq c'>0$ under Assumption \ref{assm:2.1}. In total, we now see that $0<c'\leq \lambda_p^2 \leq C'$ and thus 
    \[
    \frac{\sqrt{q}m_p}{\sqrt{c'}}\geq \frac{\sqrt{q}m_p}{\lambda_p} \geq \frac{\sqrt{q}m_p}{\sqrt{C'}}.
    \]
    Meanwhile, since $\lambda_{0,p}/\lambda_p\leq 1$, we have
    \[
    \left|\frac{\lambda_{0,p}}{\lambda_p} z_\alpha\right| \leq |z_\alpha|,
    \]
    where the right-hand side is independent of $p,q$ and $\{n_j\}$. Consequently, we have
    \begin{align*}
        \Phi\left(\frac{\sqrt{q}m_p}{\sqrt{c'}} + |z_\alpha|\right)& \geq \Phi\left(\frac{\sqrt{q}m_p}{\lambda_p} - \frac{\lambda_{0,p}}{\lambda_p} z_\alpha\right)
        \geq \Phi\left(\frac{\sqrt{q}m_p}{\sqrt{C'}} - |z_\alpha|\right).
    \end{align*}
    From the above inequality, it is easy to see that the power function (\ref{eq:pf_eq}) tends to $1$ if and only if $\sqrt{q}m_p \to \infty$.
    The proof is then complete. 
}

\section{Proofs of Auxiliary Lemmas}\label{proof:aux}
\subsection{Auxiliary Covariance Formulas}\label{app:a1}

The proofs of results in Appendix \ref{supp:aux} rely heavily on computing the covariances of products of quadratic forms of random vectors. Suppose that $\vv{z}$ is a $p$-dimensional random vector with i.i.d. entries having zero mean, unit variance, and finite fourth moment $\nu_4$, which is independent of coordinates. 

Suppose that $\vv{A}_1,\ldots,\vv{A}_4$ be arbitrary $p\times p$ deterministic symmetric matrices such that traces of all possible $\vv{A}_i$'s, $\vv{A}_i\vv{A}_j$'s, $\vv{A}_i\vv{A}_j\vv{A}_l$'s and $\vv{A}_i\vv{A}_j\vv{A}_l\vv{A}_r$'s are of $O(p)$. This is always the case when they are bounded in the operator norm, and their singular value distributions converge weakly to certain deterministic limits.
Under this assumption, we now focus on the terms in the highest orders of the following two covariances:
\[
{\rm Cov}\left((\vv{z}^{\T}\vv{A}_1\vv{z})(\vv{z}^{\T}\vv{A}_2\vv{z}),(\vv{z}^{\T}\vv{A}_3\vv{z})(\vv{z}^{\T}\vv{A}_4\vv{z})\right)
\]
and
\[{\rm Cov}\left((\vv{z}^{\T}\vv{A}_1\vv{z})(\vv{z}^{\T}\vv{A}_2\vv{z}),\vv{z}^{\T}\vv{A}_3\vv{z}\right).
\]
These leading terms determine the asymptotic variances occurring in Theorems \ref{thm:2.4} and \ref{thm:3.3}. The following lemma provides the corresponding answers of our interest.
\begin{lemma}\label{lem:appapp0}
Suppose that $\vv{z}$ and $\vv{A}_i$'s satisfy the above-mentioned conditions. Then, we have
\begin{equation}\label{eq:appapp1}
\begin{split}
&{\rm Cov}\left((\vv{z^{\T}}\vv{A}_1\vv{z})(\vv{z^{\T}}\vv{A}_2\vv{z}),(\vv{z^{\T}}\vv{A}_3\vv{z})(\vv{z^{\T}}\vv{A}_4\vv{z})\right)\\
=&\left\{\begin{aligned}
    &{\rm tr}(\vv{A}_2){\rm tr}(\vv{A}_4)\{2{\rm tr}(\vv{A}_1\vv{A}_3) + (\nu_4-3){\rm tr}(\mathcal{D}(\vv{A}_1)\mathcal{D}(\vv{A}_3))\}\\
    +&{\rm tr}(\vv{A}_2){\rm tr}(\vv{A}_4)\{2{\rm tr}(\vv{A}_1\vv{A}_4) + (\nu_4-3){\rm tr}(\mathcal{D}(\vv{A}_1)\mathcal{D}(\vv{A}_4))\}\\
    +&{\rm tr}(\vv{A}_1){\rm tr}(\vv{A}_3)\{2{\rm tr}(\vv{A}_2\vv{A}_4) + (\nu_4-3){\rm tr}(\mathcal{D}(\vv{A}_2)\mathcal{D}(\vv{A}_4))\}\\
    +&{\rm tr}(\vv{A}_1){\rm tr}(\vv{A}_4)\{2{\rm tr}(\vv{A}_2\vv{A}_3) + (\nu_4-3){\rm tr}(\mathcal{D}(\vv{A}_2)\mathcal{D}(\vv{A}_3))\}\\
\end{aligned}\right\}+O(p^2),
\end{split}
\end{equation}
and
\begin{equation}\label{eq:appapp2}
\begin{split}
&{\rm Cov}\left((\vv{z^{\T}}\vv{A}_1\vv{z})(\vv{z^{\T}}\vv{A}_2\vv{z}),\vv{z^{\T}}\vv{A}_3\vv{z}\right)\\
=&\left\{\begin{aligned}
    &{\rm tr}(\vv{A}_2)\{2{\rm tr}(\vv{A}_1\vv{A}_3) + (\nu_4-3){\rm tr}(\mathcal{D}(\vv{A}_1)\mathcal{D}(\vv{A}_3))\}\\
    +&{\rm tr}(\vv{A}_1)\{2{\rm tr}(\vv{A}_2\vv{A}_3) + (\nu_4-3){\rm tr}(\mathcal{D}(\vv{A}_2)\mathcal{D}(\vv{A}_3))\}
\end{aligned}\right\}+O(p).
\end{split}
\end{equation}
\end{lemma}
The proof of the lemma is mainly based on a direct expansion of the expectations of products of $(\vv{z}^\T\vv{A}_i\vv{z})$'s with the help of results in \cite{ullah2004finite} and \cite{Bao20101193}. To save space, we omit the tedious and lengthy calculations.  



\subsection{Proof of Equation (\ref{eq:A.3})}\label{sp:2.1}

Observe that
\begin{align*}
{\rm tr}(\vv{S}_{n}^2)=\frac{1}{n_i^2}\sum_{r,s=1}^n (\vv{z}^{\T}_r\vv{\Sigma}\vv{z}_s)^2.
\end{align*}
We have
\begin{align*}
n^{4}{\rm Var}({\rm tr}(\vv{S}_n^2))&=\sum_{r} {\rm var}\left((\vv{z^{\T}}_r\vv{\Sigma}\vv{z}_r)^2\right)+4\sum_{r\neq s}{\rm cov}\left((\vv{z^{\T}}_r\vv{\Sigma}\vv{z}_r)^2,(\vv{z^{\T}}_r\vv{\Sigma}\vv{z}_s)^2\right)\\
&+\sum_{r\neq s,k\neq l}{\rm cov}\left((\vv{z^{\T}}_r\vv{\Sigma}\vv{z}_s)^2,(\vv{z^{\T}}_k\vv{\Sigma}\vv{z}_l)^2\right)\\
&=\sum_{r} {\rm Var}\left((\vv{z^{\T}}_r\vv{\Sigma}\vv{z}_r)^2\right)+4\sum_{r\neq s}{\rm Cov}\left((\vv{z^{\T}}_r\vv{\Sigma}\vv{z}_r)^2,(\vv{z^{\T}}_r\vv{\Sigma}\vv{z}_s)^2\right)\\
&+2\sum_{r\neq s}{\rm Var}\left((\vv{z^{\T}}_r\vv{\Sigma}\vv{z}_s)^2\right)+4\sum_{r\neq s\neq k}{\rm Cov}\left((\vv{z^{\T}}_r\vv{\Sigma}\vv{z}_s)^2,(\vv{z^{\T}}_r\vv{\Sigma}\vv{z}_k)^2\right).
\end{align*}

Let us consider these four types of covariances step by step.

For the first term ${\rm Var}\left((\vv{z^{\T}}_r\vv{\Sigma}\vv{z}_r)^2\right)$, by letting $\vv{A}_1=\vv{A}_2=\vv{A}_3=\vv{A}_4=\vv{\Sigma}$ in (\ref{eq:appapp1}), we have
\begin{align*}
 {\rm Var}\left((\vv{z^{\T}}_r\vv{\Sigma}\vv{z}_r)^2\right)=8{\rm tr}(\vv{\Sigma}^2){\rm tr}(\vv{\Sigma})^2+4(\nu_{4}-3){\rm tr}(\mathcal{D}(\vv{\Sigma})^2){\rm tr}(\vv{\Sigma})^2+o(p^3).
\end{align*}

For the second term ${\rm Cov}\left((\vv{z^{\T}}_r\vv{\Sigma}\vv{z}_r)^2,(\vv{z^{\T}}_r\vv{\Sigma}\vv{z}_s)^2\right)$, we observe
\[
\E\left[(\vv{z^{\T}}_r\vv{\Sigma}\vv{z}_s)^2-{\rm tr}(\vv{\Sigma}^2)|\vv{z}_r\right]=\vv{z^{\T}}_r\vv{\Sigma}^2\vv{z}_r-{\rm tr}(\vv{\Sigma}^2).
\]
By letting $\vv{A}_1=\vv{A}_2=\vv{\Sigma}$ and $\vv{A}_3=\vv{\Sigma}^2$ in (\ref{eq:appapp2}), we then have
\begin{align*}
&{\rm Cov}\left((\vv{z^{\T}}_r\vv{\Sigma}\vv{z}_r)^2,(\vv{z^{\T}}_r\vv{\Sigma}\vv{z}_s)^2\right)\\
=&{\rm Cov}\left((\vv{z^{\T}}_r\vv{\Sigma}\vv{z}_r)^2,\vv{z^{\T}}_r\vv{\Sigma}^2\vv{z}_r\right)\\
=&4{\rm tr}(\vv{\Sigma}^3){\rm tr}(\vv{\Sigma})+2(\nu_{4}-3){\rm tr}(\mathcal{D}(\vv{\Sigma})\vv{\Sigma}^2){\rm tr}(\vv{\Sigma})+o(p^2),
\end{align*}
in which ${\rm tr}(\mathcal{D}(\vv{\Sigma})\vv{\Sigma}^2)={\rm tr}(\mathcal{D}(\vv{\Sigma})\mathcal{D}(\vv{\Sigma}^2))$.

Next, for the third term ${\rm Var}\left((\vv{z^{\T}}_r\vv{\Sigma}\vv{z}_s)^2\right)$, we have
\begin{align*}
{\rm Var}\left((\vv{z^{\T}}_r\vv{\Sigma}\vv{z}_s)^2\right)&=\E(\vv{z^{\T}}_r\vv{\Sigma}\vv{z}_s)^4-{\rm tr}(\vv{\Sigma}^2)^2.
\end{align*}
Since
\[
(\vv{z^{\T}}_r\vv{\Sigma}\vv{z}_s)^4=\left(\vv{z^{\T}}_s\left(\vv{\Sigma}\vv{z}_r\vv{z^{\T}}_r\vv{\Sigma}\right)\vv{z}_s\right)^2,
\]
we have
\begin{align*}
\E((\vv{z^{\T}}_r\vv{\Sigma}_i\vv{z}_s)^4|\vv{z}_r)&={\rm tr}(\vv{\Sigma}_i\vv{z}_r\vv{z^{\T}}_r\vv{\Sigma}_i)^2+2{\rm tr}\left[(\vv{\Sigma}_i\vv{z}_r\vv{z^{\T}}_r\vv{\Sigma}_i)^2\right]\\
&+(\nu_{4}-3){\rm tr}\left[\mathcal{D}(\vv{\Sigma}_i\vv{z}_r\vv{z^{\T}}_r\vv{\Sigma}_i)^2\right]\\
&=3(\vv{z^{\T}}_r\vv{\Sigma}_i\vv{z}_r)^2+(\nu_{4}-3){\rm tr}\left[\mathcal{D}(\vv{\Sigma}_i\vv{z}_r\vv{z^{\T}}_r\vv{\Sigma}_i)^2\right].
\end{align*}
Therefore,
\begin{align*}
{\rm Var}\left((\vv{z^{\T}}_r\vv{\Sigma}_i\vv{z}_s)^2\right)&=\E(\vv{z^{\T}}_r\vv{\Sigma}_i\vv{z}_s)^4-{\rm tr}(\vv{\Sigma}_i^2)^2\\
&=3\left[{\rm tr}(\vv{\Sigma}_i^2)^2+2{\rm tr}(\vv{\Sigma}_i^4)+(\nu_{4}-3){\rm tr}(\mathcal{D}(\vv{\Sigma}_i^2)^2)\right]\\
&+(\nu_{4}-3)\E{\rm tr}\left[\mathcal{D}(\vv{\Sigma}_i\vv{z}_r\vv{z^{\T}}_r\vv{\Sigma}_i)^2\right]-{\rm tr}(\vv{\Sigma}_i^2)^2\\
&= 2{\rm tr}(\vv{\Sigma}_i^2)^2+o(p^2)\\
&+(\nu_{4}-3)\E{\rm tr}\left[\mathcal{D}(\vv{\Sigma}_i\vv{z}_r\vv{z^{\T}}_r\vv{\Sigma}_i)^2\right]\\
&=2{\rm tr}(\vv{\Sigma}_i^2)^2+o(p^2).
\end{align*}
It remains to compute the last term ${\rm Cov}\left((\vv{z^{\T}}_r\vv{\Sigma}\vv{z}_s)^2,(\vv{z^{\T}}_r\vv{\Sigma}\vv{z}_k)^2\right)$.To see this, we let $e_l$ be the $p$-dimensional vector with its $l$th entry $1$ and the others $0$. Then,
\[
\mathcal{D}(\vv{A})=\sum_{l=1}^p \vv{e}_l\vv{e^{\T}}_l\cdot(\vv{e^{\T}}_l\vv{A}\vv{e}_l).
\]
It follows that
\[
\mathcal{D}(\vv{A})^2=\sum_{l=1}^p \vv{e}_l\vv{e^{\T}}_l\cdot(\vv{e^{\T}}_l\vv{A}\vv{e}_l)^2,~~~~~
{\rm tr}(\mathcal{D}(\vv{A})^2)=\sum_{l=1}^p (\vv{e^{\T}}_l\vv{A}\vv{e}_l)^2.
\]
Let $\vv{A}=\vv{\Sigma}\vv{z}_r\vv{z^{\T}}_r\vv{\Sigma}$. Then, we have
$$
\vv{e^{\T}}_l\vv{A}\vv{e}_l=(\vv{e^{\T}}_l\vv{\Sigma}_i\vv{z}_r)^2~~~~\hbox{so that}~~~~ (\vv{e^{\T}}_l\vv{A}\vv{e}_l)^2=(\vv{e^{\T}}_l\vv{\Sigma}_i\vv{z}_r)^4.
$$
Consequently, it holds that
\begin{align*}
\E{\rm tr}\left[\mathcal{D}(\vv{\Sigma}\vv{z}_r\vv{z^{\T}}_r\vv{\Sigma})^2\right]&=\sum_{l=1}^p \E\left[\vv{z^{\T}}_r(\vv{\Sigma}\vv{e}_l\vv{e^{\T}}_l\vv{\Sigma})\vv{z}_r\right]^2\\
&=\sum_{l=1}^p 3(\vv{e^{\T}}_l\vv{\Sigma}^2\vv{e}_l)^2+(\nu_{4}-3){\rm tr}\left[\mathcal{D}(\vv{\Sigma}\vv{e}_l\vv{e^{\T}}_l\vv{\Sigma})^2\right]\\
&=3{\rm tr}(\mathcal{D}(\vv{\Sigma}^2))+(\nu_{4}-3)\sum_{l,k=1}^p (\vv{e^{\T}}_k\vv{\Sigma}\vv{e}_l\vv{e^{\T}}_l\vv{\Sigma}_i\vv{e}_k)^2\\
&=3{\rm tr}(\mathcal{D}(\vv{\Sigma}^2))+(\nu_{4}-3)\sum_{l,k=1}^p \vv{\Sigma}_{kl}^4\\
&=3{\rm tr}(\mathcal{D}(\vv{\Sigma}^2))+(\nu_{4}-3){\rm tr}\left[(\vv{\Sigma}\circ\vv{\Sigma})^2\right],
\end{align*}
where $\circ$ stands for the Hadamard product of matrices. 
According to Corollary A.21 in \citesupp{BS10}, we have
$$
\|\vv{\Sigma}_i\circ\vv{\Sigma}_i\|_{\rm op}\leq \|\vv{\Sigma}_i\|_{\rm op}^2=O(1),
$$
so that the traces above are of order $O(p)$. 

Back to the last term, we find
\[
\E((\vv{z^{\T}}_r\vv{\Sigma}\vv{z}_s)^2|\vv{z^{\T}}_r)=\vv{z^{\T}}_r\vv{\Sigma}^2\vv{z}_r.
\]
Therefore,
\begin{align*}
{\rm Cov}\left((\vv{z^{\T}}_r\vv{\Sigma}\vv{z}_s)^2,(\vv{z^{\T}}_r\vv{\Sigma}\vv{z}_k)^2\right)
&=\E\left(\vv{z^{\T}}_r\vv{\Sigma}^2\vv{z}_r-{\rm tr}(\vv{\Sigma}^2)\right)^2\\
&=2{\rm tr}(\vv{\Sigma}^4)+(\nu_{4}-3){\rm tr}(\mathcal{D}(\vv{\Sigma}^2)^2).
\end{align*}

To recap, we have
\begin{align*}
{\rm Var}({\rm tr}(\vv{S}_n^2))&=\frac{1}{n^4}\bigg\{n\cdot \left(8{\rm tr}(\vv{\Sigma}^2){\rm tr}(\vv{\Sigma})^2+4(\nu_{4}-3){\rm tr}(\mathcal{D}(\vv{\Sigma})^2){\rm tr}(\vv{\Sigma})^2\right)\\
&+4n^2\cdot\left(4{\rm tr}(\vv{\Sigma}^3){\rm tr}(\vv{\Sigma})+2(\nu_{4}-3){\rm tr}(\mathcal{D}(\vv{\Sigma})\mathcal{D}(\vv{\Sigma}^2)){\rm tr}(\vv{\Sigma})\right)\\
&+2n^2\cdot 2{\rm tr}(\vv{\Sigma}^2)^2\\
&+4n^3\cdot \left(2{\rm tr}(\vv{\Sigma}^4)+(\nu_{4}-3){\rm tr}(\mathcal{D}(\vv{\Sigma}^2)^2)\right)+o(n^4)\bigg\}\\
&=\frac{1}{n^4}\bigg\{ 4n{\rm tr}(\vv{\Sigma})^2\mu(\vv{I}_p,\vv{I}_p|\vv{\Sigma}) +8n^2 {\rm tr}(\vv{\Sigma})\langle\vv{\Sigma},\vv{I}_p\rangle_{\vv{\Sigma}}\\
&~~~~~~~~~+4n^2{\rm tr}(\vv{\Sigma}^2)^2+4n^3\langle\vv{\Sigma},\vv{\Sigma}\rangle_{\vv{\Sigma}}\bigg\}(1+ o(1)).
\end{align*}

\subsection{Proof of Equation (\ref{eq:A.4})}\label{sp:2.2}
Recall that
\begin{align*}
{\rm tr}(\vv{S}_{n}\vv{A}_1)\left(\frac{1}{p}{\rm tr}(\vv{S}_{n})\right)=\frac{1}{pn^2}\sum_{r,s=1}^n (\vv{z}^{\T}_r\vv{\Sigma}^{1/2}\vv{A}_1\vv{\Sigma}^{1/2}\vv{z}_r)(\vv{z}^{\T}_s\vv{\Sigma}\vv{z}_s),\\
{\rm tr}(\vv{S}_{n}\vv{A}_2)\left(\frac{1}{p}{\rm tr}(\vv{S}_{n})\right)=\frac{1}{pn^2}\sum_{r,s=1}^n (\vv{z}^{\T}_r\vv{\Sigma}^{1/2}\vv{A}_2\vv{\Sigma}^{1/2}\vv{z}_r) (\vv{z}^{\T}_s\vv{\Sigma}\vv{z}_s).
\end{align*}
We have
\begin{align*}
&p^2n^4{\rm Cov}\left({\rm tr}(\vv{S}_{n}\vv{A}_1)\left(\frac{1}{p}{\rm tr}(\vv{S}_{n})\right),{\rm tr}(\vv{S}_{n}\vv{A}_2)\left(\frac{1}{p}{\rm tr}(\vv{S}_{n})\right)\right)\\
=&\sum_r {\rm Cov}\left((\vv{z^{\T}}_r\vv{\Sigma}^{1/2}\vv{A}_1\vv{\Sigma}^{1/2}\vv{z}_r)(\vv{z^{\T}}_r\vv{\Sigma}\vv{z}_r),(\vv{z^{\T}}_r\vv{\Sigma}^{1/2}\vv{A}_2\vv{\Sigma}^{1/2}\vv{z}_r)(\vv{z^{\T}}_r\vv{\Sigma}\vv{z}_r)\right)\\
+&\sum_{r\neq s}{\rm Cov}\left((\vv{z^{\T}}_r\vv{\Sigma}^{1/2}\vv{A}_1\vv{\Sigma}^{1/2}\vv{z}_r)(\vv{z^{\T}}_r\vv{\Sigma}\vv{z}_r),(\vv{z^{\T}}_r\vv{\Sigma}^{1/2}\vv{A}_2\vv{\Sigma}^{1/2}\vv{z}_r)(\vv{z^{\T}}_s\vv{\Sigma}\vv{z}_s)\right)\\
+&\sum_{r\neq s}{\rm Cov}\left((\vv{z^{\T}}_r\vv{\Sigma}^{1/2}\vv{A}_1\vv{\Sigma}^{1/2}\vv{z}_r)(\vv{z^{\T}}_r\vv{\Sigma}\vv{z}_r),(\vv{z^{\T}}_s\vv{\Sigma}^{1/2}\vv{A}_2\vv{\Sigma}^{1/2}\vv{z}_s)(\vv{z^{\T}}_r\vv{\Sigma}\vv{z}_r)\right)\\
+&\sum_{r\neq s}{\rm Cov}\left((\vv{z^{\T}}_s\vv{\Sigma}^{1/2}\vv{A}_1\vv{\Sigma}^{1/2}\vv{z}_s)(\vv{z^{\T}}_r\vv{\Sigma}\vv{z}_r),(\vv{z^{\T}}_r\vv{\Sigma}^{1/2}\vv{A}_2\vv{\Sigma}^{1/2}\vv{z}_r)(\vv{z^{\T}}_r\vv{\Sigma}\vv{z}_r)\right)\\
+&\sum_{r\neq s}{\rm Cov}\left((\vv{z^{\T}}_r\vv{\Sigma}^{1/2}\vv{A}_1\vv{\Sigma}^{1/2}\vv{z}_r)(\vv{z^{\T}}_s\vv{\Sigma}\vv{z}_s),(\vv{z^{\T}}_r\vv{\Sigma}^{1/2}\vv{A}_2\vv{\Sigma}^{1/2}\vv{z}_r)(\vv{z^{\T}}_r\vv{\Sigma}\vv{z}_r)\right)\\
+&\sum_{r\neq s}{\rm Cov}\left((\vv{z^{\T}}_r\vv{\Sigma}^{1/2}\vv{A}_1\vv{\Sigma}^{1/2}\vv{z}_r)(\vv{z^{\T}}_s\vv{\Sigma}\vv{z}_s),(\vv{z^{\T}}_r\vv{\Sigma}^{1/2}\vv{A}_2\vv{\Sigma}^{1/2}\vv{z}_r)(\vv{z^{\T}}_s\vv{\Sigma}\vv{z}_s)\right)\\
+&\sum_{r\neq s}{\rm Cov}\left((\vv{z^{\T}}_r\vv{\Sigma}^{1/2}\vv{A}_1\vv{\Sigma}^{1/2}\vv{z}_r)(\vv{z^{\T}}_s\vv{\Sigma}\vv{z}_s),(\vv{z^{\T}}_s\vv{\Sigma}^{1/2}\vv{A}_2\vv{\Sigma}^{1/2}\vv{z}_s)(\vv{z^{\T}}_r\vv{\Sigma}\vv{z}_r)\right)\\
+&\sum_{r\neq s\neq k}{\rm Cov}\left((\vv{z^{\T}}_r\vv{\Sigma}^{1/2}\vv{A}_1\vv{\Sigma}^{1/2}\vv{z}_r)(\vv{z^{\T}}_s\vv{\Sigma}\vv{z}_s),(\vv{z^{\T}}_r\vv{\Sigma}^{1/2}\vv{A}_2\vv{\Sigma}^{1/2}\vv{z}_r)(\vv{z^{\T}}_k\vv{\Sigma}\vv{z}_k)\right)\\
+&\sum_{r\neq s\neq k}{\rm Cov}\left((\vv{z^{\T}}_r\vv{\Sigma}^{1/2}\vv{A}_1\vv{\Sigma}^{1/2}\vv{z}_r)(\vv{z^{\T}}_s\vv{\Sigma}\vv{z}_s),(\vv{z^{\T}}_k\vv{\Sigma}^{1/2}\vv{A}_2\vv{\Sigma}^{1/2}\vv{z}_k)(\vv{z^{\T}}_r\vv{\Sigma}\vv{z}_r)\right)\\
+&\sum_{r\neq s\neq k}{\rm Cov}\left((\vv{z^{\T}}_r\vv{\Sigma}^{1/2}\vv{A}_1\vv{\Sigma}^{1/2}\vv{z}_r)(\vv{z^{\T}}_s\vv{\Sigma}\vv{z}_s),(\vv{z^{\T}}_s\vv{\Sigma}^{1/2}\vv{A}_2\vv{\Sigma}^{1/2}\vv{z}_s)(\vv{z^{\T}}_k\vv{\Sigma}\vv{z}_k)\right)\\
+&\sum_{r\neq s\neq k}{\rm Cov}\left((\vv{z^{\T}}_r\vv{\Sigma}^{1/2}\vv{A}_1\vv{\Sigma}^{1/2}\vv{z}_r)(\vv{z^{\T}}_s\vv{\Sigma}\vv{z}_s),(\vv{z^{\T}}_k\vv{\Sigma}^{1/2}\vv{A}_2\vv{\Sigma}^{1/2}\vv{z}_k)(\vv{z^{\T}}_s\vv{\Sigma}\vv{z}_s)\right).
\end{align*}

In what follows, we compute these different types of terms step by step.

First, by letting $\tilde{\vv{A}}_1=\vv{\Sigma}^{1/2}\vv{A}_1\vv{\Sigma}^{1/2}$, $\tilde{A}_3 = \vv{\Sigma}^{1/2}\vv{A}_2\vv{\Sigma}^{1/2}$ and $\vv{A}_2=\vv{A}_4=\vv{\Sigma}$ in  (\ref{eq:appapp1}), we have 
\begin{align*}
&{\rm Cov}\left((\vv{z^{\T}}_r\vv{\Sigma}^{1/2}\vv{A}_1\vv{\Sigma}^{1/2}\vv{z}_r)(\vv{z^{\T}}_r\vv{\Sigma}\vv{z}_r),(\vv{z^{\T}}_r\vv{\Sigma}^{1/2}\vv{A}_2\vv{\Sigma}^{1/2}\vv{z}_r)(\vv{z^{\T}}_r\vv{\Sigma}\vv{z}_r)\right)\\
=&\{2{\rm tr}(\vv{\Sigma}\vv{A}_1\vv{\Sigma}\vv{A}_2)+(\nu_4 -3){\rm tr}(\mathcal{D}(\vv{\Sigma}^{1/2}\vv{A}_1\vv{\Sigma}^{1/2})\mathcal{D}(\vv{\Sigma}^{1/2}\vv{A}_2\vv{\Sigma}^{1/2}))\}{\rm tr}(\vv{\Sigma})^2\\
+&\{2{\rm tr}(\vv{\Sigma}^2\vv{A}_1)+(\nu_4-3){\rm tr}(\mathcal{D}(\vv{\Sigma}_i)\mathcal{D}(\vv{\Sigma}^{1/2}\vv{A}_1\vv{\Sigma}^{1/2}))\}{\rm tr}(\vv{\Sigma}\vv{A}_2){\rm tr}(\vv{\Sigma})\\
+&\{2{\rm tr}(\vv{\Sigma}^2\vv{A}_2)+(\nu_4-3){\rm tr}(\mathcal{D}(\vv{\Sigma}_i)\mathcal{D}(\vv{\Sigma}^{1/2}\vv{A}_2\vv{\Sigma}^{1/2}))\}{\rm tr}(\vv{\Sigma}\vv{A}_1){\rm tr}(\vv{\Sigma})\\
+&\{2{\rm tr}(\vv{\Sigma}^2)+(\nu_4-3){\rm tr}(\mathcal{D}(\vv{\Sigma})^2)\}{\rm tr}(\vv{\Sigma}\vv{A}_1){\rm tr}(\vv{\Sigma}\vv{A}_2)+o(p^3).
\end{align*}

Next, we observe that
\begin{align*}
&\E\left[(\vv{z^{\T}}_r\vv{\Sigma}^{1/2}\vv{A}_1\vv{\Sigma}^{1/2}\vv{z}_r)(\vv{z^{\T}}_r\vv{\Sigma}\vv{z}_r)(\vv{z^{\T}}_r\vv{\Sigma}^{1/2}\vv{A}_2\vv{\Sigma}^{1/2}\vv{z}_r)(\vv{z^{\T}}_s\vv{\Sigma}\vv{z}_s)|\vv{z}_r\right]\\
=&(\vv{z^{\T}}_r\vv{\Sigma}^{1/2}\vv{A}_1\vv{\Sigma}^{1/2}\vv{z}_r)(\vv{z^{\T}}_r\vv{\Sigma}\vv{z}_r)(\vv{z^{\T}}_r\vv{\Sigma}^{1/2}\vv{A}_2\vv{\Sigma}^{1/2}\vv{z}_r){\rm tr}(\vv{\Sigma}).
\end{align*}
It then follows by letting $\tilde{\vv{A}}_1=\vv{\Sigma}^{1/2}\vv{A}_1\vv{\Sigma}^{1/2}$, $\tilde{A}_3 = \vv{\Sigma}^{1/2}\vv{A}_2\vv{\Sigma}^{1/2}$ and $\vv{A}_2=\vv{\Sigma}$ in (\ref{eq:appapp2}) that
\begin{align*}
&{\rm Cov}\left((\vv{z^{\T}}_r\vv{\Sigma}^{1/2}\vv{A}_1\vv{\Sigma}^{1/2}\vv{z}_r)(\vv{z^{\T}}_r\vv{\Sigma}\vv{z}_r),(\vv{z^{\T}}_r\vv{\Sigma}^{1/2}\vv{A}_2\vv{\Sigma}^{1/2}\vv{z}_r)(\vv{z^{\T}}_s\vv{\Sigma}\vv{z}_s)\right)\\
=&{\rm tr}(\vv{\Sigma}){\rm Cov}\left((\vv{z^{\T}}_r\vv{\Sigma}^{1/2}\vv{A}_1\vv{\Sigma}^{1/2}\vv{z}_r)(\vv{z^{\T}}_r\vv{\Sigma}\vv{z}_r),(\vv{z^{\T}}_r\vv{\Sigma}^{1/2}\vv{A}_2\vv{\Sigma}^{1/2}\vv{z}_r)\right)\\
=&{\rm tr}(\vv{\Sigma})\bigg\{ [2{\rm tr}(\vv{\Sigma}\vv{A}_1\vv{\Sigma}\vv{A}_2)+(\nu_4-3){\rm tr}(\mathcal{D}(\vv{\Sigma}^{1/2}\vv{A}_1\vv{\Sigma}^{1/2})\mathcal{D}(\vv{\Sigma}^{1/2}\vv{A}_2\vv{\Sigma}^{1/2}))]{\rm tr}(\vv{\Sigma})\\
&~~~~~~~+[2{\rm tr}(\vv{\Sigma}^2\vv{A}_2)+(\nu_4-3){\rm tr}(\mathcal{D}(\vv{\Sigma})\mathcal{D}(\vv{\Sigma}^{1/2}\vv{A}_2\vv{\Sigma}^{1/2}))]{\rm tr}(\vv{\Sigma}\vv{A}_1)\bigg\}+o(p^3).
\end{align*}
Similarly, we also have
\begin{align*}
    &{\rm Cov}\left((\vv{z^{\T}}_r\vv{\Sigma}^{1/2}\vv{A}_1\vv{\Sigma}^{1/2}\vv{z}_r)(\vv{z^{\T}}_r\vv{\Sigma}\vv{z}_r),(\vv{z^{\T}}_s\vv{\Sigma}^{1/2}\vv{A}_2\vv{\Sigma}^{1/2}\vv{z}_s)(\vv{z^{\T}}_r\vv{\Sigma}\vv{z}_r)\right)\\
=&{\rm tr}(\vv{\Sigma}\vv{A}_2)\bigg\{ [2{\rm tr}(\vv{\Sigma}^2\vv{A}_1)+(\nu_4-3){\rm tr}(\mathcal{D}(\vv{\Sigma}^{1/2}\vv{A}_1\vv{\Sigma}^{1/2})\mathcal{D}(\vv{\Sigma}))]{\rm tr}(\vv{\Sigma})\\
&~~~~~~~+[2{\rm tr}(\vv{\Sigma}^2)+(\nu_4-3){\rm tr}(\mathcal{D}(\vv{\Sigma})^2)]{\rm tr}(\vv{\Sigma}\vv{A}_1)\bigg\}+o(p^3),\\
&{\rm Cov}\left((\vv{z^{\T}}_r\vv{\Sigma}^{1/2}\vv{A}_1\vv{\Sigma}^{1/2}\vv{z}_r)(\vv{z^{\T}}_s\vv{\Sigma}\vv{z}_s),(\vv{z^{\T}}_r\vv{\Sigma}^{1/2}\vv{A}_2\vv{\Sigma}^{1/2}\vv{z}_r)(\vv{z^{\T}}_r\vv{\Sigma}\vv{z}_r)\right)\\
=&{\rm tr}(\vv{\Sigma})\bigg\{ [2{\rm tr}(\vv{\Sigma}\vv{A}_1\vv{\Sigma}\vv{A}_2)+(\nu_4-3){\rm tr}(\mathcal{D}(\vv{\Sigma}^{1/2}\vv{A}_1\vv{\Sigma}^{1/2})\mathcal{D}(\vv{\Sigma}^{1/2}\vv{A}_2\vv{\Sigma}^{1/2}))]{\rm tr}(\vv{\Sigma})\\
&~~~~~~~+[2{\rm tr}(\vv{\Sigma}^2\vv{A}_1)+(\nu_4-3){\rm tr}(\mathcal{D}(\vv{\Sigma})\mathcal{D}(\vv{\Sigma}^{1/2}\vv{A}_1\vv{\Sigma}^{1/2}))]{\rm tr}(\vv{\Sigma}\vv{A}_2)\bigg\}+o(p^3),\\
&{\rm Cov}\left((\vv{z^{\T}}_s\vv{\Sigma}^{1/2}\vv{A}_1\vv{\Sigma}^{1/2}\vv{z}_s)(\vv{z^{\T}}_r\vv{\Sigma}\vv{z}_r),(\vv{z^{\T}}_r\vv{\Sigma}^{1/2}\vv{A}_2\vv{\Sigma}^{1/2}\vv{z}_r)(\vv{z^{\T}}_r\vv{\Sigma}\vv{z}_r)\right)\\
=&{\rm tr}(\vv{\Sigma}\vv{A}_1)\bigg\{ [2{\rm tr}(\vv{\Sigma}^2\vv{A}_2)+(\nu_4-3){\rm tr}(\mathcal{D}(\vv{\Sigma}^{1/2}\vv{A}_2\vv{\Sigma}^{1/2})\mathcal{D}(\vv{\Sigma}))]{\rm tr}(\vv{\Sigma})\\
&~~~~~~~+[2{\rm tr}(\vv{\Sigma}^2)+(\nu_4-3){\rm tr}(\mathcal{D}(\vv{\Sigma})^2)]{\rm tr}(\vv{\Sigma}\vv{A}_2)\bigg\}+o(p^3).
\end{align*}

Furthermore, we focus on the cases in which there are two different indices, and each index occurs four times. We observe that
\begin{align*}
&\E\left[(\vv{z^{\T}}_r\vv{\Sigma}^{1/2}\vv{A}_1\vv{\Sigma}^{1/2}\vv{z}_r)(\vv{z^{\T}}_s\vv{\Sigma}\vv{z}_s),(\vv{z^{\T}}_r\vv{\Sigma}^{1/2}\vv{A}_2\vv{\Sigma}^{1/2}\vv{z}_r)(\vv{z^{\T}}_s\vv{\Sigma}\vv{z}_s)|\vv{z}_r\right]\\
=&(\vv{z^{\T}}_r\vv{\Sigma}^{1/2}\vv{A}_1\vv{\Sigma}^{1/2}\vv{z}_r)(\vv{z^{\T}}_r\vv{\Sigma}^{1/2}\vv{A}_2\vv{\Sigma}^{1/2}\vv{z}_r)
\E[(\vv{z^{\T}}_s\vv{\Sigma}\vv{z}_s)^2]\\
=&(\vv{z^{\T}}_r\vv{\Sigma}^{1/2}\vv{A}_1\vv{\Sigma}^{1/2}\vv{z}_r)(\vv{z^{\T}}_r\vv{\Sigma}^{1/2}\vv{A}_2\vv{\Sigma}^{1/2}\vv{z}_r)\left({\rm tr}(\vv{\Sigma})^2+2{\rm tr}(\vv{\Sigma}^2)+(\nu_{4}-3){\rm tr}(\mathcal{D}(\vv{\Sigma})^2)\right),
\end{align*}
so that
\begin{align*}
&\E\left[(\vv{z^{\T}}_r\vv{\Sigma}^{1/2}\vv{A}_1\vv{\Sigma}^{1/2}\vv{z}_r)(\vv{z^{\T}}_s\vv{\Sigma}\vv{z}_s),(\vv{z^{\T}}_r\vv{\Sigma}^{1/2}\vv{A}_2\vv{\Sigma}^{1/2}\vv{z}_r)(\vv{z^{\T}}_s\vv{\Sigma}\vv{z}_s)\right]\\
=&\E[(\vv{z^{\T}}_r\vv{\Sigma}^{1/2}\vv{A}_1\vv{\Sigma}^{1/2}\vv{z}_r)(\vv{z^{\T}}_r\vv{\Sigma}^{1/2}\vv{A}_2\vv{\Sigma}^{1/2}\vv{z}_r)]\left({\rm tr}(\vv{\Sigma})^2+2{\rm tr}(\vv{\Sigma}^2)+(\nu_{4}-3){\rm tr}(\mathcal{D}(\vv{\Sigma})^2)\right)\\
=& \left\{{\rm tr}(\vv{\Sigma}\vv{A}_1){\rm tr}(\vv{\Sigma}\vv{A}_2) + 2{\rm tr}(\vv{\Sigma}\vv{A}_1\vv{\Sigma}\vv{A}_2) + (\nu_4-3){\rm tr}(\mathcal{D}(\vv{\Sigma}^{1/2}\vv{A}_1\vv{\Sigma}^{1/2})\mathcal{D}(\vv{\Sigma}^{1/2}\vv{A}_2\vv{\Sigma}^{1/2}))\right\}\\
\times&\left\{{\rm tr}(\vv{\Sigma})^2+2{\rm tr}(\vv{\Sigma}^2)+(\nu_{4}-3){\rm tr}(\mathcal{D}(\vv{\Sigma})^2)\right\}
\end{align*}
It then follows that
\begin{align*}
&{\rm Cov}\left((\vv{z^{\T}}_r\vv{\Sigma}^{1/2}\vv{A}_1\vv{\Sigma}^{1/2}\vv{z}_r)(\vv{z^{\T}}_s\vv{\Sigma}\vv{z}_s),(\vv{z^{\T}}_r\vv{\Sigma}^{1/2}\vv{A}_2\vv{\Sigma}^{1/2}\vv{z}_r)(\vv{z^{\T}}_s\vv{\Sigma}\vv{z}_s)\right)\\
=&\E\left[(\vv{z^{\T}}_r\vv{\Sigma}^{1/2}\vv{A}_1\vv{\Sigma}^{1/2}\vv{z}_r)(\vv{z^{\T}}_s\vv{\Sigma}\vv{z}_s),(\vv{z^{\T}}_r\vv{\Sigma}^{1/2}\vv{A}_2\vv{\Sigma}^{1/2}\vv{z}_r)(\vv{z^{\T}}_s\vv{\Sigma}\vv{z}_s)\right]-{\rm tr}(\vv{\Sigma}\vv{A}_1){\rm tr}(\vv{\Sigma}\vv{A}_2){\rm tr}(\vv{\Sigma})\\
=&{\rm tr}(\vv{\Sigma})^2\{2{\rm tr}(\vv{\Sigma}\vv{A}_1\vv{\Sigma}\vv{A}_2) + (\nu_4-3){\rm tr}(\mathcal{D}(\vv{\Sigma}^{1/2}\vv{A}_1\vv{\Sigma}^{1/2})\mathcal{D}(\vv{\Sigma}^{1/2}\vv{A}_2\vv{\Sigma}^{1/2}))\}\\
+&{\rm tr}(\vv{\Sigma}\vv{A}_1){\rm tr}(\vv{\Sigma}\vv{A}_2)\{2{\rm tr}(\vv{\Sigma}^2)+(\nu_{4}-3){\rm tr}(\mathcal{D}(\vv{\Sigma})^2)\}+o(p^3).
\end{align*}

Similarly, we also have
\begin{align*}
&{\rm Cov}\left((\vv{z^{\T}}_r\vv{\Sigma}^{1/2}\vv{A}_1\vv{\Sigma}^{1/2}\vv{z}_r)(\vv{z^{\T}}_s\vv{\Sigma}\vv{z}_s),(\vv{z^{\T}}_s\vv{\Sigma}^{1/2}\vv{A}_2\vv{\Sigma}^{1/2}\vv{z}_s)(\vv{z^{\T}}_r\vv{\Sigma}\vv{z}_r)\right)\\
=&{\rm tr}(\vv{\Sigma}\vv{A}_1){\rm tr}(\vv{\Sigma})\{2{\rm tr}(\vv{\Sigma}^2\vv{A}_2) + (\nu_4-3){\rm tr}(\mathcal{D}(\vv{\Sigma}^{1/2}\vv{A}_2\vv{\Sigma}^{1/2})\mathcal{D}(\vv{\Sigma}))\}\\
+&{\rm tr}(\vv{\Sigma}\vv{A}_2){\rm tr}(\vv{\Sigma})\{2{\rm tr}(\vv{\Sigma}^2\vv{A}_1)+(\nu_4-3){\rm tr}(\mathcal{D}(\vv{\Sigma}^{1/2}\vv{A}_1\vv{\Sigma}^{1/2})\mathcal{D}(\vv{\Sigma}))\}+o(p^3).
\end{align*}

Finally, we focus on the cases with three different indices. Observe that
\begin{align*}
    &{\rm Cov}\left((\vv{z^{\T}}_r\vv{\Sigma}^{1/2}\vv{A}_1\vv{\Sigma}^{1/2}\vv{z}_r)(\vv{z^{\T}}_s\vv{\Sigma}\vv{z}_s),(\vv{z^{\T}}_r\vv{\Sigma}^{1/2}\vv{A}_2\vv{\Sigma}^{1/2}\vv{z}_r)(\vv{z^{\T}}_k\vv{\Sigma}\vv{z}_k)\right)\\
    =& {\rm tr}(\vv{\Sigma})^2 {\rm Cov}((\vv{z^{\T}}_r\vv{\Sigma}^{1/2}\vv{A}_1\vv{\Sigma}^{1/2}\vv{z}_r),(\vv{z^{\T}}_r\vv{\Sigma}^{1/2}\vv{A}_2\vv{\Sigma}^{1/2}\vv{z}_r))\\
    =& {\rm tr}(\vv{\Sigma})^2\{2{\rm tr}(\vv{\Sigma}\vv{A}_1\vv{\Sigma}\vv{A}_2) + (\nu_4-3){\rm tr}(\mathcal{D}(\vv{\Sigma}^{1/2}\vv{A}_1\vv{\Sigma}^{1/2})\mathcal{D}(\vv{\Sigma}^{1/2}\vv{A}_2\vv{\Sigma}^{1/2}))\}.
\end{align*}
Similarly, the other three cases satisfy
\begin{align*}
     &{\rm Cov}\left((\vv{z^{\T}}_r\vv{\Sigma}^{1/2}\vv{A}_1\vv{\Sigma}^{1/2}\vv{z}_r)(\vv{z^{\T}}_s\vv{\Sigma}\vv{z}_s),(\vv{z^{\T}}_k\vv{\Sigma}^{1/2}\vv{A}_2\vv{\Sigma}^{1/2}\vv{z}_k)(\vv{z^{\T}}_r\vv{\Sigma}\vv{z}_r)\right)\\
    =& {\rm tr}(\vv{\Sigma}){\rm tr}(\vv{\Sigma}\vv{A}_2)\{2{\rm tr}(\vv{\Sigma}^2\vv{A}_1) + (\nu_4-3){\rm tr}(\mathcal{D}(\vv{\Sigma}^{1/2}\vv{A}_1\vv{\Sigma}^{1/2})\mathcal{D}(\vv{\Sigma}))\},\\
    &{\rm Cov}\left((\vv{z^{\T}}_r\vv{\Sigma}^{1/2}\vv{A}_1\vv{\Sigma}^{1/2}\vv{z}_r)(\vv{z^{\T}}_s\vv{\Sigma}\vv{z}_s),(\vv{z^{\T}}_s\vv{\Sigma}^{1/2}\vv{A}_2\vv{\Sigma}^{1/2}\vv{z}_s)(\vv{z^{\T}}_k\vv{\Sigma}\vv{z}_k)\right)\\
    =& {\rm tr}(\vv{\Sigma}){\rm tr}(\vv{\Sigma}\vv{A}_1)\{2{\rm tr}(\vv{\Sigma}^2\vv{A}_2) + (\nu_4-3){\rm tr}(\mathcal{D}(\vv{\Sigma}^{1/2}\vv{A}_2\vv{\Sigma}^{1/2})\mathcal{D}(\vv{\Sigma}))\},\\
    &{\rm Cov}\left((\vv{z^{\T}}_r\vv{\Sigma}^{1/2}\vv{A}_1\vv{\Sigma}^{1/2}\vv{z}_r)(\vv{z^{\T}}_s\vv{\Sigma}\vv{z}_s),(\vv{z^{\T}}_k\vv{\Sigma}^{1/2}\vv{A}_2\vv{\Sigma}^{1/2}\vv{z}_k)(\vv{z^{\T}}_s\vv{\Sigma}\vv{z}_s)\right)\\
    =& {\rm tr}(\vv{\Sigma}\vv{A}_1){\rm tr}(\vv{\Sigma}\vv{A}_2)\{2{\rm tr}(\vv{\Sigma}^2) + (\nu_4-3){\rm tr}(\mathcal{D}(\vv{\Sigma})^2)\}.
\end{align*}

To recap, we have
\begin{align*}
&{\rm Cov}\left({\rm tr}(\vv{S}_{n}\vv{A}_1)\left(\frac{1}{p}{\rm tr}(\vv{S}_{n})\right),{\rm tr}(\vv{S}_{n}\vv{A}_2)\left(\frac{1}{p}{\rm tr}(\vv{S}_{n})\right)\right)\\
=&\frac{1}{p^2n}\bigg\{{\rm tr}(\vv{\Sigma})^2\{2{\rm tr}(\vv{\Sigma}\vv{A}_1\vv{\Sigma}\vv{A}_2) + (\nu_4-3){\rm tr}(\mathcal{D}(\vv{\Sigma}^{1/2}\vv{A}_1\vv{\Sigma}^{1/2})\mathcal{D}(\vv{\Sigma}^{1/2}\vv{A}_2\vv{\Sigma}^{1/2}))\}\\
&+{\rm tr}(\vv{\Sigma}){\rm tr}(\vv{\Sigma}\vv{A}_2)\{2{\rm tr}(\vv{\Sigma}^2\vv{A}_1) + (\nu_4-3){\rm tr}(\mathcal{D}(\vv{\Sigma}^{1/2}\vv{A}_1\vv{\Sigma}^{1/2})\mathcal{D}(\vv{\Sigma}))\}\\
&+{\rm tr}(\vv{\Sigma}){\rm tr}(\vv{\Sigma}\vv{A}_1)\{2{\rm tr}(\vv{\Sigma}^2\vv{A}_2) + (\nu_4-3){\rm tr}(\mathcal{D}(\vv{\Sigma}^{1/2}\vv{A}_2\vv{\Sigma}^{1/2})\mathcal{D}(\vv{\Sigma}))\}\\
&+{\rm tr}(\vv{\Sigma}\vv{A}_1){\rm tr}(\vv{\Sigma}\vv{A}_2)\{2{\rm tr}(\vv{\Sigma}^2) + (\nu_4-3){\rm tr}(\mathcal{D}(\vv{\Sigma})^2)\}\bigg\}+o(1)\\
&=\frac{1}{p^2n}
\left\{\begin{aligned}
&{\rm tr}(\vv{\Sigma}){\rm tr}(\vv{\Sigma})\langle\vv{A}_1,\vv{A}_2\rangle_{\vv{\Sigma}}\\
+&{\rm tr}(\vv{\Sigma}){\rm tr}(\vv{\Sigma}\vv{A}_1)\langle\vv{I}_p,\vv{A}_2\rangle_{\vv{\Sigma}}\\
+&{\rm tr}(\vv{\Sigma}){\rm tr}(\vv{\Sigma}\vv{A}_2)\langle\vv{I}_p,\vv{A}_1\rangle_{\vv{\Sigma}}\\
+&{\rm tr}(\vv{\Sigma}\vv{A}_1){\rm tr}(\vv{\Sigma}\vv{A}_2)\langle\vv{I}_p,\vv{I}_p\rangle_{\vv{\Sigma}}
\end{aligned}\right\}(1+o(1))
\end{align*}

\subsection{Proof of Equation (\ref{eq:A.4.1})}\label{sp:2.4}

Recall that 
\begin{align*}
    {\rm tr}(\vv{S}_{n}^2) &= \frac{1}{n^2}\sum_{r,s =1}^n (\vv{z}_r^\T\vv{\Sigma}\vv{z}_s)^2\\
    {\rm tr}(\vv{S}_n\vv{A})\left(\frac{1}{p}{\rm tr}(\vv{S}_n)\right)&=  \frac{1}{pn^2} \sum_{r,s=1}^n (\vv{z}_r^\T\vv{\Sigma}^{1/2}\vv{A}\vv{\Sigma}^{1/2}\vv{z}_r)(\vv{z}_s^\T\vv{\Sigma}\vv{z}_s).
\end{align*}
We have
\begin{align*}
&pn^4{\rm Cov}\left({\rm tr}(\vv{S}_{n}^2),{\rm tr}(\vv{S}_n\vv{A})\left(\frac{1}{p}{\rm tr}(\vv{S}_n)\right)\right)\\
=&\sum_{r} {\rm Cov}\left((\vv{z^{\T}}_r\vv{\Sigma}\vv{z}_r)^2,(\vv{z^{\T}}_r\vv{\Sigma}^{1/2}\vv{A}\vv{\Sigma}^{1/2}\vv{z}_r)(\vv{z^{\T}}_r\vv{\Sigma}\vv{z}_r)\right)\\
+&\sum_{r\neq s}{\rm Cov}\left((\vv{z^{\T}}_r\vv{\Sigma}\vv{z}_r)^2,(\vv{z^{\T}}_r\vv{\Sigma}^{1/2}\vv{A}\vv{\Sigma}^{1/2}\vv{z}_r)(\vv{z^{\T}}_s\vv{\Sigma}\vv{z}_s)\right)\\
+&\sum_{r\neq s}{\rm Cov}\left((\vv{z^{\T}}_r\vv{\Sigma}\vv{z}_r)^2,(\vv{z^{\T}}_s\vv{\Sigma}^{1/2}\vv{A}\vv{\Sigma}^{1/2}\vv{z}_s)(\vv{z^{\T}}_r\vv{\Sigma}\vv{z}_r)\right)\\
+&2\sum_{r\neq s}{\rm Cov}\left((\vv{z^{\T}}_r\vv{\Sigma}\vv{z}_s)^2,(\vv{z^{\T}}_r\vv{\Sigma}^{1/2}\vv{A}\vv{\Sigma}_i^{1/2}\vv{z}_r)(\vv{z^{\T}}_r\vv{\Sigma}\vv{z}_r)\right)\\
+&2\sum_{r\neq s}{\rm Cov}\left((\vv{z^{\T}}_r\vv{\Sigma}\vv{z}_s)^2,(\vv{z^{\T}}_r\vv{\Sigma}^{1/2}\vv{A}\vv{\Sigma}^{1/2}\vv{z}_r)(\vv{z^{\T}}_s\vv{\Sigma}\vv{z}_s)\right)\\
+&2\sum_{r\neq s\neq k}{\rm Cov}\left((\vv{z^{\T}}_r\vv{\Sigma}\vv{z}_s)^2,(\vv{z^{\T}}_r\vv{\Sigma}^{1/2}\vv{A}\vv{\Sigma}^{1/2}\vv{z}_r)(\vv{z^{\T}}_k\vv{\Sigma}\vv{z}_k)\right)\\
+&2\sum_{r\neq s\neq k}{\rm Cov}\left((\vv{z^{\T}}_r\vv{\Sigma}\vv{z}_s)^2,(\vv{z^{\T}}_k\vv{\Sigma}^{1/2}\vv{A}\vv{\Sigma}^{1/2}\vv{z}_k)(\vv{z^{\T}}_r\vv{\Sigma}\vv{z}_r)\right).
\end{align*}

In what follows, we compute these different types of terms step by step.

First, by letting $\vv{A}_1=\vv{A}_2=\vv{A}_3=\vv{\Sigma}$ and $\vv{A}_4=\vv{\Sigma}^{1/2}\vv{A}\vv{\Sigma}^{1/2}$ in (\ref{eq:appapp1}), we have
\begin{align*}
&{\rm Cov}\left((\vv{z^{\T}}_r\vv{\Sigma}\vv{z}_r)^2,(\vv{z^{\T}}_r\vv{\Sigma}^{1/2}\vv{A}\vv{\Sigma}^{1/2}\vv{z}_r)(\vv{z^{\T}}_r\vv{\Sigma}\vv{z}_r)\right)\\
=&4{\rm tr}(\vv{\Sigma}^2){\rm tr}(\vv{\Sigma}){\rm tr}(\vv{\Sigma}_i\vv{A})+4{\rm tr}(\vv{\Sigma}^2\vv{A}){\rm tr}(\vv{\Sigma})^2\\
+&2(\nu_{4}-3)\bigg\{{\rm tr}(\mathcal{D}(\vv{\Sigma})^2){\rm tr}(\vv{\Sigma}){\rm tr}(\vv{\Sigma}\vv{A})+{\rm tr}(\mathcal{D}(\vv{\Sigma})\vv{\Sigma}^{1/2}\vv{A}\vv{\Sigma}^{1/2}){\rm tr}(\vv{\Sigma})^2\bigg\}+o(p^3).
\end{align*}

Next, by letting $\vv{A}_1=\vv{A}_2=\vv{\Sigma}$ and $\vv{A}_3=\vv{\Sigma}^{1/2}\vv{A}\vv{\Sigma}^{1/2}$ in (\ref{eq:appapp2}), we have
\begin{align*}
&{\rm Cov}\left((\vv{z^{\T}}_r\vv{\Sigma}_i\vv{z}_r)^2,(\vv{z^{\T}}_r\vv{\Sigma}^{1/2}\vv{A}\vv{\Sigma}^{1/2}\vv{z}_r)(\vv{z^{\T}}_s\vv{\Sigma}\vv{z}_s)\right)\\
=&{\rm tr}(\vv{\Sigma}){\rm Cov}\left((\vv{z^{\T}}_r\vv{\Sigma}\vv{z}_r)^2,(\vv{z^{\T}}_r\vv{\Sigma}^{1/2}\vv{A}\vv{\Sigma}^{1/2}\vv{z}_r)\right)\\
=&{\rm tr}(\vv{\Sigma})\bigg\{4{\rm tr}(\vv{\Sigma}^2\vv{A}){\rm tr}(\vv{\Sigma})+2(\nu_{4}-3){\rm tr}(\mathcal{D}(\vv{\Sigma}^{1/2}\vv{A}\vv{\Sigma}^{1/2})\vv{\Sigma}){\rm tr}(\vv{\Sigma})\bigg\}+o(p^3).
\end{align*}

Meanwhile, by letting $\vv{A}_1=\vv{A}_2=\vv{A}_3=\vv{\Sigma}_i$ in (\ref{eq:appapp2}), we have
\begin{align*}
&{\rm Cov}\left((\vv{z^{\T}}_r\vv{\Sigma}\vv{z}_r)^2,(\vv{z^{\T}}_s\vv{\Sigma}^{1/2}\vv{A}\vv{\Sigma}^{1/2}\vv{z}_s)(\vv{z^{\T}}_r\vv{\Sigma}\vv{z}_r)\right)\\
=&{\rm tr}(\vv{\Sigma}\vv{A})\bigg\{4{\rm tr}(\vv{\Sigma}^2){\rm tr}(\vv{\Sigma})+2(\nu_{4}-3){\rm tr}(\mathcal{D}(\vv{\Sigma})^2){\rm tr}(\vv{\Sigma})\bigg\}+o(p^3).
\end{align*}

Next, we observe 
\[
\E\left((\vv{z^{\T}}_r\vv{\Sigma}\vv{z}_s)^2|\vv{z}_r\right)=\vv{z^{\T}}_r\vv{\Sigma}^2\vv{z}_r.
\]
It then follows by letting $\vv{A}_1=\vv{\Sigma}^{1/2}\vv{A}\vv{\Sigma}^{1/2}$, $\vv{A}_2=\vv{\Sigma}$ and $\vv{A}_3=\vv{\Sigma}^2$ in (\ref{eq:appapp2}) that
\begin{align*}
&{\rm Cov}\left((\vv{z^{\T}}_r\vv{\Sigma}\vv{z}_s)^2,(\vv{z^{\T}}_r\vv{\Sigma}^{1/2}\vv{\Sigma}_0\vv{\Sigma}^{1/2}\vv{z}_r)(\vv{z^{\T}}_r\vv{\Sigma}\vv{z}_r)\right)\\
=&{\rm Cov}\left(\vv{z^{\T}}_r\vv{\Sigma}^2\vv{z}_r,(\vv{z^{\T}}_r\vv{\Sigma}^{1/2}\vv{A}\vv{\Sigma}^{1/2}\vv{z}_r)(\vv{z^{\T}}_r\vv{\Sigma}\vv{z}_r)\right)\\
=&2{\rm tr}(\vv{\Sigma}^3){\rm tr}(\vv{\Sigma}\vv{A})+2{\rm tr}(\vv{\Sigma}^3\vv{A}){\rm tr}(\vv{\Sigma})\\
+&(\nu_{4}-3)\bigg\{{\rm tr}(\mathcal{D}(\vv{\Sigma}^2)\vv{\Sigma}){\rm tr}(\vv{\Sigma}\vv{A})+{\rm tr}(\mathcal{D}(\vv{\Sigma}^2)\vv{\Sigma}^{1/2}\vv{A}\vv{\Sigma}^{1/2}){\rm tr}(\vv{\Sigma})\bigg\}+o(p^2).
\end{align*}

Finally, we observe that
\begin{align*}
&\E\left[(\vv{z^{\T}}_r\vv{\Sigma}\vv{z}_s)^2,(\vv{z^{\T}}_r\vv{\Sigma}^{1/2}\vv{A}\vv{\Sigma}^{1/2}\vv{z}_r)(\vv{z^{\T}}_s\vv{\Sigma}\vv{z}_s)|\vv{z}_r\right]\\
=&(\vv{z^{\T}}_r\vv{\Sigma}^{1/2}\vv{A}\vv{\Sigma}^{1/2}\vv{z}_r)\left[{\rm tr}(\vv{\Sigma})(\vv{z^{\T}}_r\vv{\Sigma}^2\vv{z}_r)+2(\vv{z^{\T}}_r\vv{\Sigma}^3\vv{z}_r)+(\nu_{4}-3)\vv{z^{\T}}_r\vv{\Sigma}\mathcal{D}(\vv{\Sigma})\vv{\Sigma}\vv{z}_r\right].
\end{align*}
It then follows that
\begin{align*}
&{\rm Cov}\left((\vv{z^{\T}}_r\vv{\Sigma}\vv{z}_s)^2,(\vv{z^{\T}}_r\vv{\Sigma}^{1/2}\vv{A}\vv{\Sigma}^{1/2}\vv{z}_r)(\vv{z^{\T}}_s\vv{\Sigma}\vv{z}_s)\right)\\
=&\E\left[(\vv{z^{\T}}_r\vv{\Sigma}\vv{z}_s)^2,(\vv{z^{\T}}_r\vv{\Sigma}^{1/2}\vv{A}\vv{\Sigma}^{1/2}\vv{z}_r)(\vv{z^{\T}}_s\vv{\Sigma}\vv{z}_s)\right]-{\rm tr}(\vv{\Sigma}^2){\rm tr}(\vv{\Sigma}\vv{A}){\rm tr}(\vv{\Sigma})\\
=&{\rm tr}(\vv{\Sigma})\left(2{\rm tr}(\vv{\Sigma}^3\vv{A})+(\nu_{4}-3){\rm tr}(\mathcal{D}(\vv{\Sigma}^2)\vv{\Sigma}^{1/2}\vv{A}\vv{\Sigma}^{1/2})\right)\\
+&2{\rm tr}(\vv{\Sigma}\vv{A}){\rm tr}(\vv{\Sigma}^3)+(\nu_{4}-3){\rm tr}(\vv{\Sigma}\vv{A}){\rm tr}(\mathcal{D}(\vv{\Sigma})\vv{\Sigma}^2)+o(p^2).
\end{align*}

We have
\begin{align*}
&{\rm Cov}\left((\vv{z^{\T}}_r\vv{\Sigma}\vv{z}_s)^2,(\vv{z^{\T}}_r\vv{\Sigma}^{1/2}\vv{A}\vv{\Sigma}^{1/2}\vv{z}_r)(\vv{z^{\T}}_k\vv{\Sigma}\vv{z}_k)\right)\\
=&{\rm tr}(\vv{\Sigma}){\rm cov}\left((\vv{z^{\T}}_r\vv{\Sigma}\vv{z}_s)^2,(\vv{z^{\T}}_r\vv{\Sigma}^{1/2}\vv{A}\vv{\Sigma}^{1/2}\vv{z}_r)\right)\\
=&{\rm tr}(\vv{\Sigma}){\rm cov}\left((\vv{z^{\T}}_r\vv{\Sigma}^2\vv{z}_r),(\vv{z^{\T}}_r\vv{\Sigma}^{1/2}\vv{A}\vv{\Sigma}^{1/2}\vv{z}_r)\right)\\
=&{\rm tr}(\vv{\Sigma})\bigg\{2{\rm tr}(\vv{\Sigma}^3\vv{A})+(\nu_{4}-3){\rm tr}(\mathcal{D}(\vv{\Sigma}^2)\vv{\Sigma}^{1/2}\vv{A}\vv{\Sigma}^{1/2})\bigg\}.
\end{align*}

Finally,
\begin{align*}
&{\rm Cov}\left((\vv{z^{\T}}_r\vv{\Sigma}\vv{z}_s)^2,(\vv{z^{\T}}_k\vv{\Sigma}^{1/2}\vv{A}\vv{\Sigma}^{1/2}\vv{z}_k)(\vv{z^{\T}}_r\vv{\Sigma}\vv{z}_r)\right)\\
=&{\rm tr}(\vv{\Sigma}\vv{A}){\rm cov}\left((\vv{z^{\T}}_r\vv{\Sigma}^2\vv{z}_r),(\vv{z^{\T}}_r\vv{\Sigma}\vv{z}_r)\right)\\
=&{\rm tr}(\vv{\Sigma}\vv{A})\bigg\{2{\rm tr}(\vv{\Sigma}^3)+(\nu_{4}-3){\rm tr}(\mathcal{D}(\vv{\Sigma}^2)\vv{\Sigma})\bigg\}.
\end{align*}

To recap, we have
\begin{align*}
&pn^4{\rm Cov}\left({\rm tr}(\vv{S}_{n}^2),{\rm tr}(\vv{S}_n\vv{A})\left(\frac{1}{p}{\rm tr}(\vv{S}_n)\right)\right)\\
=&\frac{1}{pn^4}\bigg\{n^2{\rm tr}(\vv{\Sigma})\bigg[4{\rm tr}(\vv{\Sigma}^2\vv{A}){\rm tr}(\vv{\Sigma})+2(\nu_{4}-3){\rm tr}(\mathcal{D}(\vv{\Sigma}^{1/2}\vv{A}\vv{\Sigma}^{1/2})\vv{\Sigma}){\rm tr}(\vv{\Sigma})\bigg]\\
+&n^2{\rm tr}(\vv{\Sigma}\vv{A})\bigg[4{\rm tr}(\vv{\Sigma}^2){\rm tr}(\vv{\Sigma})+2(\nu_{4}-3){\rm tr}(\mathcal{D}(\vv{\Sigma})^2){\rm tr}(\vv{\Sigma})\bigg]\\
+&2n^3{\rm tr}(\vv{\Sigma})\bigg[2{\rm tr}(\vv{\Sigma}^3\vv{A})+(\nu_{4}-3){\rm tr}(\mathcal{D}(\vv{\Sigma}^2)\vv{\Sigma}^{1/2}\vv{A}\vv{\Sigma}^{1/2})\bigg]\\
+&2n^3{\rm tr}(\vv{\Sigma}\vv{A})\bigg[2{\rm tr}(\vv{\Sigma}^3)+(\nu_{4}-3){\rm tr}(\mathcal{D}(\vv{\Sigma}^2)\vv{\Sigma})\bigg]\bigg\}(1+o(1))\\
=&\frac{2}{n^4}\left\{\begin{aligned}
&n^2\cdot{\rm tr}(\vv{\Sigma}){\rm tr}(\vv{\Sigma})\langle\vv{A},\vv{I}_p\rangle_{\vv{\Sigma}}\\
+&n^2\cdot{\rm tr}(\vv{\Sigma}\vv{A}){\rm tr}(\vv{\Sigma})\langle\vv{I}_p,\vv{I}_p\rangle_{\vv{\Sigma}}\\
+&n^3\cdot {\rm tr}(\vv{\Sigma}\vv{A})\langle\vv{I}_p,\vv{\Sigma}\rangle_{\vv{\Sigma}}\\
+&n^3\cdot {\rm tr}(\vv{\Sigma})\langle\vv{A},\vv{\Sigma}\rangle_{\vv{\Sigma}}
\end{aligned}\right\}(1+o(1)).
\end{align*}

\subsection{Proof of Equation (\ref{eq:lemA3.1})}

Recall that 
\begin{align*}
{\rm tr}(\vv{S}_{n}^2)&=\frac{1}{n^2}\sum_{r,s=1}^n (\vv{z^{\T}}_r\vv{\Sigma}\vv{z}_s)^2,\\
{\rm tr}(\vv{S}_{n}\vv{A})&=\frac{1}{n}\sum_{k=1}^n \vv{z}^{\T}_{k}\vv{\Sigma}^{1/2}\vv{A}\vv{\Sigma}^{1/2}\vv{z}_{k}.
\end{align*}
Then, we have 
\begin{align*}
&{\rm Cov}\left({\rm tr}(\vv{S}_{n}^2),{\rm tr}(\vv{S}_{n}\vv{A})\right)\\
=&\frac{1}{n^3}\bigg\{\sum_{j=1}^{n}{\rm Cov}((\vv{z}^{\T}_{j}\vv{\Sigma}\vv{z}_{j})^2,\vv{z}^{\T}_{j}\vv{\Sigma}^{1/2}\vv{A}\vv{\Sigma}^{1/2}\vv{z}_{j})\\
&+2\sum_{l\neq j}{\rm Cov}(\vv{z}^{\T}_{j}\vv{\Sigma}^2\vv{z}_{j},\vv{z}^{\T}_{j}\vv{\Sigma}^{1/2}\vv{A}\vv{\Sigma}^{1/2}\vv{z}_{j})\bigg\}\\
&=c_{p}^2 \frac{1}{p^2}\left[4{\rm tr}(\vv{\Sigma}^2\vv{A}){\rm tr}(\vv{\Sigma})+2(\nu_4-3){\rm tr}(\mathcal{D}(\vv{\Sigma})\vv{\Sigma}^{1/2}\vv{A}\vv{\Sigma}^{1/2}){\rm tr}(\vv{\Sigma})\right]\\
&+c_{p}\frac{1}{p}\left[2{\rm tr}(\vv{\Sigma}^3\vv{A})+(\nu_4-3){\rm tr}(\mathcal{D}(\vv{\Sigma}^2)\vv{\Sigma}^{1/2}\vv{A}\vv{\Sigma}^{1/2})\right]\\
&= 2c_p^2\left(\frac{1}{p}{\rm tr}(\vv{\Sigma})\right) \langle\vv{I}_p,\vv{A}\rangle_{\vv{\Sigma}} + c_p\langle \vv{\Sigma},\vv{A}\rangle_{\vv{\Sigma}}
\end{align*}

\subsection{Proof of Equation (\ref{eq:lemA3.2})}
Recall that
\begin{align*}
\frac{1}{n}{\rm tr}(\vv{S}_{n})^2&=\frac{1}{n^3}\sum_{r,s=1}^n (\vv{z^{\T}}_r\vv{\Sigma}\vv{z}_r)(\vv{z^{\T}}_s\vv{\Sigma}\vv{z}_s),\\
{\rm tr}(\vv{S}_{n}\vv{A})&=\frac{1}{n}\sum_{k=1}^n \vv{z}^{\T}_{k}\vv{\Sigma}^{1/2}\vv{A}\vv{\Sigma}^{1/2}\vv{z}_{k}.
\end{align*}
We then have
\begin{align*}
&{\rm Cov}\left(\frac{1}{n}{\rm tr}(\vv{S}_{n})^2,{\rm tr}(\vv{S}_n\vv{A})\right)\\
=& \frac{1}{n^4}\bigg\{\sum_{j=1}^n {\rm Cov}((\vv{z^{\T}}_r\vv{\Sigma}\vv{z}_r)^2,\vv{z}^{\T}_{j}\vv{\Sigma}^{1/2}\vv{A}\vv{\Sigma}^{1/2}\vv{z}_{j})\\
&+2\sum_{l\neq j}{\rm Cov}((\vv{z^{\T}}_j\vv{\Sigma}\vv{z}_j)(\vv{z^{\T}}_l\vv{\Sigma}\vv{z}_l),\vv{z}^{\T}_{j}\vv{\Sigma}^{1/2}\vv{A}\vv{\Sigma}^{1/2}\vv{z}_{j})\bigg\}\\
=&c_{p}\frac{1}{p^2} {\rm tr}(\vv{\Sigma})\left[2{\rm tr}(\vv{\Sigma}^2\vv{A})+{\rm tr}(\mathcal{D}(\vv{\Sigma})\vv{\Sigma}^{1/2}\vv{A}\vv{\Sigma}^{1/2})\right]+O(\frac{1}{p})\\
=& 2c_p^2\left(\frac{1}{p}{\rm tr}(\vv{\Sigma})\right)\langle\vv{I}_p,\vv{A}\rangle_{\vv{\Sigma}}
\end{align*}

\bibliographystylesupp{imsart-nameyear}
\bibliographysupp{mybibsp,mybibfile}

\end{document}